\documentclass[mnsc,nonblindrev]{oo} 
\OneAndAHalfSpacedXI 

\usepackage{natbib}
 \bibpunct[, ]{(}{)}{,}{a}{}{,}%
 \def\bibfont{\small}%
\TheoremsNumberedThrough     
\ECRepeatTheorems

\EquationsNumberedThrough    

\MANUSCRIPTNO{}

\usepackage{amsmath, amssymb, bm, mathtools}

\usepackage{pgf, tikz}
\usepackage{graphicx}

\usepackage{pgfplots}
\pgfplotsset{compat=1.17}

\usetikzlibrary{arrows, automata, shapes, decorations.pathreplacing, calligraphy}

\usepackage{algorithm}
\usepackage[noend]{algpseudocode}
\algloopdefx{RepeatUntil}[1]{\textbf{repeat}  \textbf{until} #1}

\usepackage{booktabs}
\usepackage{multirow}

\usepackage{url}
\usepackage[bookmarksnumbered=true]{hyperref}
\usepackage{soul}
\usepackage{xcolor,colortbl}

\hypersetup{
	colorlinks,
	linkcolor={red!50!black},
	citecolor={blue!50!black},
	urlcolor={blue!80!black}
}

\def\L{{\mathcal L}}

\def\R{\mathbb{R}}

\def\hat{\widehat}

\def \P{\mathcal{P}}

\def \Z{{\mathbb{Z}}}

\def\K{{\mathcal K}}
\def\L{{\mathcal L}}

\def\M{{\mathcal M}}
\def\S{{\mathcal S}}

\def\N{{\mathcal N}}

\def\T{{\mathcal T}}

\def\P{{\mathcal P}}

\def\C{{\mathcal C}}

\newcommand{\exclude}[1]{}

\newcommand{\algFalse}{\textbf{false}}
\newcommand{\algTrue}{\textbf{true}}

\newcommand{\algReturn}{\textbf{return}}

\newcommand{\algCutAdd}{\texttt{cutAdded}}
\newcommand{\algSDDP}{\texttt{SDDP}}
\newcommand{\algBC}{\texttt{SDDP-B\&C}}
\newcommand{\algConverge}{\texttt{SDDPConverge}}
\newcommand{\algExact}{\texttt{exact}}
\newcommand{\algBenders}{\texttt{Benders}}

\newcommand{\tree}{\T}
\newcommand{\stage}{\N}
\newcommand{\staget}{\N_t}
\newcommand{\patht}{\P}
\newcommand{\childs}{\C}
\newcommand{\childsub}{\C^S}
\newcommand{\periods}{T}

\newcommand{\modality}{\L}
\newcommand{\horizon}{[\periods]}
\newcommand{\markov}{\M}
\newcommand{\rootnode}{\textsf{r}}
\newcommand{\parent}{a}
\newcommand{\prob}{p}
\newcommand{\pprob}{\bar{p}}

\newcommand{\nodes}{\N}
\newcommand{\nodesnoroot}{\nodes\setminus\{\rootnode\}}

\newcommand{\shelters}{\mathcal{I}}
\newcommand{\centers}{\mathcal{J}}

\newcommand{\Qref}{Q^R}

\newcommand{\Qaggre}{Q^A}

\newcommand{\hn}{\texttt{HN}}
\newcommand{\ma}{\texttt{MA}}
\newcommand{\pa}{\texttt{PM}}
\newcommand{\mm}{\texttt{MM}}
\newcommand{\fh}{\texttt{FH}}

\newcommand{\dimmarkov}{s}
 
\newcommand{\tmarkovhistory}{m^t}
\newcommand{\maptmarkov}{\Phi_t}

\newcommand{\za}{z^A}
\newcommand{\ztM}{z_{t,\tmarkovhistory}}
\newcommand{\ztphiM}{\za_{t,\maptmarkov \tmarkovhistory}}
\newcommand{\zphi}{\za_{\phi_t(n)}}
\newcommand{\zphiparent}{\za_{\phi_{t-1}(\parent(n))}}

\newcommand{\subSet}{\S}
\newcommand{\subSett}{\subSet_t}
\newcommand{\subnodes}{\nodes^S}
\newcommand{\subnodest}{\subnodes_t}
\newcommand{\Qsub}{\overline{Q}^A}
\newcommand{\Qsubsol}{\hat{Q}^A}
\newcommand{\sub}{\sigma}

\newcommand{\coeffs}{\mathcal{B}}
\newcommand{\subn}{\sigma}
\newcommand{\subnparent}{\parent(\subn)}

\newcommand{\zcopy}{\zeta}
\newcommand{\zcopyn}{\zcopy_{\subn}}
\newcommand{\zcopyparent}{\zcopy_{\subnparent}}
\newcommand{\zcopyparentsol}{\hat{\zcopy}_{\subnparent}}
\newcommand{\zcopyphin}{\zcopy_{\subn,\phi_t(\subn)}}
\newcommand{\zcopyphiparent}{\zcopy_{\subnparent,\phi_{t-1}(\subnparent)}}
\newcommand{\zcopyphisol}{\hat{\zcopy}_{\phi_{t-1}(\subn)}}

\newcommand{\mapldr}{\Psi_t} 
\newcommand{\varldr}{\mu} 
\newcommand{\varldrsol}{\hat{\varldr}} 
\newcommand{\dimxit}{l_t}
\newcommand{\dimxitt}{l^t}

\newcommand{\Qldr}{{Q}^L}
\newcommand{\Qldrsol}{\hat{Q}^L}

\newcommand{\Qldrsub}{\overline{Q}^L}
\newcommand{\coeffsldr}{\coeffs^L}

\newcommand{\mldr}{\texttt{LDR-M}}
\newcommand{\tldr}{\texttt{LDR-T}}
\newcommand{\thldr}{\texttt{LDR-TH}}

\newcommand{\extensive}{\texttt{Ex}}
\newcommand{\sddp}{\texttt{S}}
\newcommand{\sddplb}{\texttt{S-LB}}
\newcommand{\sddpub}{\texttt{S-UB}}

\newcommand{\zn}{z_n}
\newcommand{\xn}{x_n}
\newcommand{\yn}{y_n}
\newcommand{\xnsub}{x_{\subn}}
\newcommand{\ynsub}{y_{\subn}}
\newcommand{\znparent}{z_{\parent(n)}}
\newcommand{\xnparent}{x_{\parent(n)}}
\newcommand{\xnparentsub}{x_{\subnparent}}
\newcommand{\znprime}{z_{n'}}
\newcommand{\yroot}{y_{\rootnode}}
\newcommand{\xroot}{x_{\rootnode}}
\newcommand{\zroot}{z_{\rootnode}}

\newcommand{\xrootsol}{\hat{\xroot}}
\newcommand{\zasol}{\hat{z}^A}
\newcommand{\thetasol}{\hat{\theta}}

\newcommand{\xinv}{x^{\texttt{I}}}
\newcommand{\xcap}{x^{\texttt{C}}}

\newcommand{\bzero}{\boldsymbol{0}}
\newcommand{\identity}{\mathbb{I}}

\tikzstyle{arc} = [draw,line width=0.5pt,->,gray!80!black]
\tikzstyle{arc back} = [draw,line width=1pt,->,violet!80!black]
\tikzstyle{arc bold} = [arc,
preaction={
	draw,gray!30!white,-,
	double=gray!30!white,
	double distance=3pt,
}]

\tikzstyle{point} = [circle,fill=gray!80!black,font=\scriptsize, inner sep=2pt]
\tikzstyle{node} = [circle, draw = gray!80!black, font=\scriptsize, inner sep=2.5pt]

\tikzstyle{node-square} = [draw = gray!80!black, font=\scriptsize, inner sep=3pt]
\tikzstyle{node-polygon5} = [regular polygon, regular polygon sides=5, draw = gray!80!black, font=\scriptsize, inner sep=2.2pt]
\tikzstyle{node-polygon6} = [regular polygon,regular polygon sides=6, draw = gray!80!black, font=\scriptsize, inner sep=2.2pt]

\tikzstyle{node light} = [circle, draw = gray!80!black, fill=gray!30!white, font=\scriptsize, inner sep=2.5pt]
\tikzstyle{node dark} = [circle, draw = gray!80!black, fill=gray!95!white, font=\scriptsize, inner sep=2.5pt]
\tikzstyle{text node} = [font=\scriptsize]

\newcommand{\mcDark}[1]{\scalebox{#1}{\tikz \node[node dark] {};}}
\newcommand{\mcLight}[1]{\scalebox{#1}{\tikz \node[node light] {};}}

\newcommand{\pentaDarkDark}{\tikz \node[node-polygon5, violet!80!black, fill=violet!80!black] {};}
\newcommand{\pentaLightDark}{\tikz \node[node-polygon5,violet!40!white, ,fill=violet!40!white] {};}
\newcommand{\hexaDarkDark}{\tikz \node[node-polygon6, draw=cyan!70!black, fill=cyan!70!black] {};}

\newcommand{\squareDark}[1]{\scalebox{#1}{\tikz \node[node-square, draw=blue!70!white, fill=blue!70!white] {};}}
\newcommand{\squareLight}[1]{\scalebox{#1}{\tikz \node[node-square, draw=blue!40!white, fill=blue!40!white] {};}}
\newcommand{\pentaDark}[1]{\scalebox{#1}{\tikz \node[node-polygon5, draw=violet!70!white, fill=violet!70!white] {};}}

\newcommand{\squareRed}[1]{\scalebox{#1}{\tikz \node[node-square, draw=red, fill=red] {};}}
\newcommand{\pentaGreen}[1]{\scalebox{#1}{\tikz \node[node-polygon5, draw=green!70!black, fill=green!70!black] {};}}
\newcommand{\hexaOrange}[1]{\scalebox{#1}{\tikz \node[node-polygon6, draw=orange, fill=orange] {};}}

\newcommand{\scaleTableResults}{0.85}
\newcommand{\stretchTableResults}{0.8}

\newcommand{\Song}[1]{{\color{black}#1}}
\newcommand{\Margarita}[1]{{\color{black}#1}}
\newcommand{\Merve}[1]{{\color{black}#1}}
\newcommand{\heatmap}[4]{
\begin{tikzpicture}[scale=0.55]
    \pgfplotstableread{#1}\activetable
    \pgfplotstableread{#2}\colortable

    \foreach \x in {0,...,3}
        \foreach \y in {0,...,4}
        {
        \pgfplotstablegetelem{\y}{x\x}\of{\colortable}
        \draw[white, fill=violet!\pgfplotsretval!white] (\x,4 - \y) -- (\x, 4 - \y + 1) -- (\x + 1, 4 - \y + 1) -- (\x + 1, 4 - \y) -- (\x, 4 - \y);
        }
    
    \draw[step=1,gray!50!white,thin] (0,0) grid (4,5);
    \draw[black] (0,0) -- (0,5) -- (4,5) -- (4,0) -- (0,0);
    
    \foreach \x in {1,...,4}
        {\pgfmathtruncatemacro{\label}{\x}
            \node[gray!90!black]  (x\x) at (-0.4 +\x, -0.4) {\tiny \x};}

    \foreach \y in {1,...,5}
        {\pgfmathtruncatemacro{\label}{\y}
            \node[gray!90!black]  (y\y) at (-0.4, -0.4 + \y) {\tiny \y};}

    \foreach \x in {0,...,3}
        \foreach \y in {0,...,4}
        {\pgfplotstablegetelem{\y}{x\x}\of{\activetable}
            \node (n) at (\x + 0.5, 4 - \y + 0.5) {\tiny 
                \ifthenelse{\equal{#4}{1}}{
                \ifthenelse{\equal{\pgfplotsretval}{0.0}}{}{\pgfplotsretval}
                }
                {\pgfplotsretval}
                }; }
    \node (n) at (2, 5.5) {\small #3};
    
\end{tikzpicture}
}

\newcommand{\heatmapColor}[3]{
\begin{tikzpicture}[scale=0.55]
    \pgfplotstableread{#1}\activetable
    \pgfplotstableread{#2}\colortable

    \foreach \x in {0,...,3}
        \foreach \y in {0,...,4}
        {
        \pgfplotstablegetelem{\y}{x\x}\of{\colortable}
        \draw[white, fill=violet!\pgfplotsretval!white] (\x,4 - \y) -- (\x, 4 - \y + 1) -- (\x + 1, 4 - \y + 1) -- (\x + 1, 4 - \y) -- (\x, 4 - \y);
        }
    
    \draw[step=1,gray!50!white,thin] (0,0) grid (4,5);
    \draw[black] (0,0) -- (0,5) -- (4,5) -- (4,0) -- (0,0);
    
    \foreach \x in {1,...,4}
        {\pgfmathtruncatemacro{\label}{\x}
            \node[gray!90!black]  (x\x) at (-0.4 +\x, -0.4) {\tiny \x};}

    \foreach \y in {1,...,5}
        {\pgfmathtruncatemacro{\label}{\y}
            \node[gray!90!black]  (y\y) at (-0.4, -0.4 + \y) {\tiny \y};}

    \foreach \x in {0,...,3}
        \foreach \y in {0,...,4}
        {\pgfplotstablegetelem{\y}{x\x}\of{\activetable}
            \node (n) at (\x + 0.5, 4 - \y + 0.5) {\tiny 
                \ifthenelse{\equal{\pgfplotsretval}{0.0}}{}{\pgfplotsretval}
                } ;}

    \foreach \y in {0,...,4}
        {\pgfplotstablegetelem{\y}{c}\of{\activetable}
            \ifthenelse{\equal{\pgfplotsretval}{.}}{}{
                \ifthenelse{\equal{\pgfplotsretval}{r}}{
                \node[node-square, draw=red!70!white, fill=red!70!white, scale=0.75] (n) at (4.3, 4 - \y + 0.5) {};
                }{
                \node[node-polygon6, draw=orange, fill=orange, scale=0.75] (n) at (4.3, 4 - \y + 0.5) {};
                }
            }
        }
    
    \node (n) at (2, 5.5) {\small #3};
    
\end{tikzpicture}
}

\begin{document}


\RUNAUTHOR{Castro, Bodur, and Song}

\RUNTITLE{Markov  chain-based  policies  for  multi-stage  stochastic  integer  linear  programming}

\TITLE{
Markov Chain-based Policies for Multi-stage Stochastic Integer Linear Programming with an Application to Disaster Relief Logistics}

\ARTICLEAUTHORS{%
\AUTHOR{Margarita P. Castro}
\AFF{Department of Industrial and Systems Engineering, Pontificia Universidad Católica de Chile, Santiago 7820436, Chile \EMAIL{margarita.castro@ing.puc.cl}} 
\AUTHOR{Merve Bodur}
\AFF{Department of Mechanical and Industrial Engineering, University of Toronto, Toronto, Ontario M5S 3GH, Canada, \EMAIL{bodur@mie.utoronto.ca}}
\AUTHOR{Yongjia Song}
\AFF{Department of Industrial Engineering, Clemson University, Clemson, South Carolina 29631, \EMAIL{yongjis@clemson.edu}}
} 

\ABSTRACT{%
We introduce an aggregation framework to address multi-stage stochastic programs with mixed-integer state variables and continuous local variables (MSILPs). Our aggregation framework imposes additional structure to the integer state variables by leveraging the information of the underlying stochastic process, which is modeled as a Markov chain (MC). We demonstrate that the aggregated MSILP can be solved exactly via a branch-and-cut algorithm integrated with a variant of stochastic dual dynamic programming. To improve tractability, we propose to use this approach to obtain dual bounds. Moreover, we apply two-stage linear decision rule (2SLDR) approximations, in particular a new MC-based variant that we propose, to obtain high-quality decision policies with significantly reduced computational effort. We test the proposed methodologies in an MSILP model for hurricane disaster relief logistics planning. Our empirical evaluation compares the effectiveness of the various proposed approaches and analyzes the trade-offs between policy flexibility, solution quality, and computational effort. Specifically, the 2SLDR approximation yields provable high-quality solutions for our test instances supported by the proposed bounding procedure. We also extract valuable managerial insights from the solution behaviors exhibited by the underlying decision policies.

}%



\maketitle

\section{Introduction} \label{sec:introduction}

Multi-stage stochastic integer linear programming (MSILP) problems \Merve{form an important} class of 
\Merve{optimization}
models for sequential decision-making under uncertainty. These problems consider decisions at each stage of the sequential process, which consists of local variables---
that only participate in a single stage locally---and state variables---
that link multiple stages together. 
\Merve{In this paper,} we consider MSILP problems with continuous local variables and  mixed-integer state variables, which arise in many applications such as hydro-power scheduling \citep{hjelmeland2018nonconvex}, unit commitment \citep{zou2018multistage}, and disaster relief logistics planning (see Section \ref{sec:application}).

Despite the vast applicability of MSILP, these problems are computationally prohibitive to solve due to the non-convexity of the feasibility set caused by the integer state variables. Existing solution methods utilize the \Song{stage-wise decomposition framework, such as the stochastic dual dynamic programming (SDDP) algorithm, with convexification procedures for the expected cost-to-go functions, which require certain \Merve{limiting} assumptions such as pure binary state variables~\citep{zou2019stochastic} and Lipschitz continuous cost-to-go functions~\citep{ahmed2020stochastic}.}

We propose a \Song{different perspective to handle such MSILP problems without attempting to convexify the expected cost-to-go functions \Merve{that are parametrized by} a mixture of continuous and integer state variables}. Instead, we construct a \textit{partially extended reformulation} by transforming the original problem into one that has integer variables only in the first stage. By doing so, we can employ decomposition algorithms, such as the branch-and-cut (B\&C) algorithm, to decompose the problem into a master problem that corresponds to a first-stage mixed integer linear program, and a subproblem defined for the remaining stages, which is 
a multi-stage stochastic linear program that can be handled by stage-wise decomposition algorithms such as the SDDP algorithm. However, this reformulation may lead to exponentially many first-stage integer variables, \Merve{thus can be} computationally challenging \Merve{to solve}.

To alleviate the main challenge of the partially extended reformulation, we present an aggregation framework that imposes certain structures to the reformulation by aggregating the integer state variables
. This framework could potentially leverage the structure of the underlying stochastic process to impose a restriction on the problem at hand that leads to an informative decision policy. In particular, we consider the case where the underlying stochastic process is modeled as a Markov chain (MC) and present several aggregation schemes based on the MC state information that lead to a wide range of aggregated MSILP approximations.

\Song{
\Merve{We demonstrate that} the aggregated MSILP models
\Merve{can be solved exactly by a} B\&C framework integrated with a variant of the SDDP algorithm.} The framework decomposes the problem into a master problem and a set of subproblems, where the master problem is concerned with the first-stage decisions including all the integer state variables, and the subproblems deal with the remaining stages, which involve continuous state and local variables. \Song{While this framework guarantees an optimal solution for the aggregated MSILP, solving a multi-stage stochastic program (via SDDP) at each incumbent solution encountered in the branch-and-bound process can be computationally prohibitive. We therefore settle with a relaxation bound from this approach (by relaxing the \Merve{SDDP} termination criteria
), and pursue a more tractable approach that further approximates the aggregated MSILP as a two-stage stochastic program using two-stage linear decision rules (2SLDRs) \citep{bodur2018two}. We present 2SLDR alternatives from the literature and propose a new variant that leverages the MC structure of the underlying stochastic process. The resulting problem can be solved via Benders decomposition or any of its enhanced versions.} 

We perform a case study for the proposed aggregation framework and solution methods on a hurricane disaster relief logistic planning problem (HDR), where the stochastic process is given by an MC that models the evolution of the hurricane. We show how to create several aggregated MSILPs and 2SLDR approximations for this application. Our numerical experiments compare these alternatives, find suitable transformations and 2SLDR variants to obtain high-quality solutions with provable bounds, and present managerial insights about the resulting policies. 

The remainder of the paper is organized as follows. Section \ref{sec:framework} presents our aggregation framework and a literature review on MSILP methodologies. Section \ref{sec:methodology} introduces our novel B\&C framework integrated with the SDDP algorithm and describes the 2SLDR scheme with our proposed MC-based variant. Section \ref{sec:application} presents  a case study of the proposed methodologies for HDR. Section \ref{sec:experiments} presents the numerical experiment results and we end with some concluding remarks in Section \ref{sec:conclusions}.

\section{Aggregation Framework for Integer State Variables} \label{sec:framework}

This section describes the proposed aggregation framework to generate informative policies for a class of MSILP problems, 
\Merve{with} the main idea \Merve{of aggregating}
integer state variables in the MSILP by leveraging the structure of the underlying stochastic process.

\subsection{Problem Formulation for MSILP}

We consider a class of MSILP models with continuous local variables and mixed-integer state variables. Although integer variables may also appear as local variables in a generic MSILP, the restriction here (local variables being continuous only) does not limit the model from being applicable to a broad range of problems that arise in real-world applications. For example, a typical multi-period stochastic optimization problem with local production and distribution constraints can be represented with this model, such as the one we study in Section~\ref{sec:application}, motivated by the application of disaster relief logistics planning. Other possible applications include hydro-power scheduling \citep{hjelmeland2018nonconvex} and unit commitment \citep{zou2018multistage} problems.

We assume that the underlying stochastic process has finite support and, as such, it can be represented by a scenario tree $\tree$. This assumption is not restrictive since we can construct an approximate scenario tree \Merve{by means of sampling} if the stochastic process has continuous support, a common practice in the literature~\citep{shapiro2021lectures}. Let $\periods$ be the number of stages and let $\nodes$ be the set of nodes associated with the scenario tree $\tree$. The set of nodes in each stage $t \in \horizon:=\{1,...,\periods\}$ is given by $\staget$. The root node is denoted by $\rootnode$ and it is the only node in the first stage, i.e., $\stage_1=\{\rootnode\}$. Each node $n \in \staget$ with $t>1$ has a unique parent in $\stage_{t-1}$, denoted by $\parent(n)$. Thus, there is a unique path from $\rootnode$ to any node $n\in \staget$ in stage $t>1$, and we let $\patht(n)$ represent the set of nodes on this path (including $\rootnode$ and $n$). For each non-leaf node $n$ (i.e., $n\in \staget$ for $t \in [\periods-1]$), $\childs(n)$ is the set of children of $n$, that is, the nodes whose parent is $n$. The probability that node $n$ occurs is $\prob_n$, and we have $\prob_n>0$ for all $n\in \nodes$ and $\sum_{n \in \staget} \prob_n = 1$ for all $t\in \horizon$. The transition probability from node $n \in \staget$ for $t\in [\periods-1]$ to $n' \in \childs(n)$ can then be written as $\pprob_{nn'}:= \prob_{n'}/\prob_n$. 

Next, we present a generic problem formulation for the class of MSILP problems that we consider in~\eqref{mod:scenarioTree}, which is often referred to as the nested formulation. For each node $n\in \nodes$, let $\xn \in \R^k$ and $\zn \in \Z^\ell$ be the set of continuous and integer state variables, respectively, and $\yn\in \R^r$ be the continuous local variables. The objective coefficient vectors associated with $\xn, \yn$, and $\zn$ are denoted by $d_n, h_n$, and $c_n$, respectively. The overall objective is to minimize the sum of the local cost and the expected future cost at the root node:
\vspace{-0em}
\begin{align}
    Q_{\rootnode} = \min\; & c_\rootnode^\top \zroot + d_\rootnode^\top \xroot + h_\rootnode^\top \yroot + \sum_{n \in \childs(\rootnode)}\pprob_{\rootnode n} Q_n(\xroot,\zroot ) , \tag{$P$} \label{mod:scenarioTree} \\
    \text{s.t.}\ &
    H_{\rootnode} \zroot \geq g_{\rootnode}, \;
    J_{\rootnode} \xroot \geq f_{\rootnode}, \;
    C_{\rootnode} \xroot + D_{\rootnode} \zroot + E_{\rootnode} \yroot  \geq  b_{\rootnode}, \label{eq:st_root}\\
    &\xroot \in \R^k,\; \yroot \in \R^r,\; \zroot \in \Z^\ell, \nonumber
\end{align}
where the cost-to-go function $Q_n(\cdot, \cdot)$ associated with each non-root node $n\in \nodes\setminus \{\rootnode\}$ is:
\vspace{-0em}
\begin{subequations}
\begin{align}
    Q_n(\xnparent, \znparent) = \min\; &  c_n^\top \zn + d_n^\top \xn + h_n^\top \yn + \sum_{n' \in \childs(n)}\pprob_{nn'} Q_{n'}(\xn,\zn),\nonumber \\
    \text{s.t.}\ & H_n \zn \geq  G_n \znparent + g_n, \label{eq:st_onlyz} \\
    & J_n \xn \geq  F_n \xnparent + f_n, \label{eq:st_onlyx} \\
    & C_n \xn + D_n \zn + E_n \yn  \geq A_n \xnparent + B_n  \znparent + b_n, \label{eq:st_all} \\
    &\xn \in \R^k,\; \yn \in \R^r,\; \zn \in \Z^\ell. \nonumber
\end{align}
\end{subequations}
The cost-to-go function $Q_n(\cdot, \cdot)$ associated with a leaf node $n\in \stage_\periods$ only involves the cost incurred in the terminal stage $\periods$ since its set of children nodes is empty, i.e., $\childs(n) = \emptyset$. Our proposed formulation considers, separately, one set of constraints~\eqref{eq:st_onlyz} for integer state variables, one set of constraints~\eqref{eq:st_onlyx} for continuous state variables, and one set of constraints~\eqref{eq:st_all} that link state and local variables. Although alternative formulations are possible, this form accommodates the presentation of the proposed aggregation framework, which will be described in Section~\ref{sec:methodology}.

Problem \eqref{mod:scenarioTree} is notoriously challenging to solve for two reasons. The first challenge is the nature of multi-stage stochastic programs where the underlying stochastic process is modeled as a scenario tree of an exponentially large size. Thus, it is impractical to solve these large-scale problems directly and decomposition methods can be deemed necessary. The second challenge is the existence of integer decision variables, which makes the expected cost-to-go functions (i.e., the value functions that represent the expected future cost defined at each node of the scenario tree) to be nonconvex in general. The following section reviews existing methods in the literature to handle \eqref{mod:scenarioTree} and similar problems, focusing on how 
\Merve{they attempt to overcome}
these challenges. 

\subsection{Literature Review}

In this section, we discuss recent advances in the literature that address MSILPs like formulation~\eqref{mod:scenarioTree}. 
Nested Benders decomposition~\citep{birge1985decomposition} and its sampling variant under the assumption of stage-wise independence, SDDP~\citep{pereira1991multi}, are classical decomposition algorithms for multi-stage stochastic linear programs. 
Recently,~\cite{zou2019stochastic} develop the SDDiP algorithm for MSILPs with pure binary state variables, where the nonconvex expected cost-to-go functions are lower approximated by a piecewise linear convex envelope constructed by the so-called Lagrangian cuts. These cuts are guaranteed to be exact at points where the associated state variables are binary, leading to the exactness and finite convergence of the  algorithm. The problem becomes more challenging when the set of state variables includes general mixed-integer decision variables, which is the case for formulation~\eqref{mod:scenarioTree}. One option is to perform a binarization procedure to reformulate general integer state variables and approximate continuous state variables via binary variables, respectively, and then apply the SDDiP algorithm~\citep{zou2018multistage}. However, this approach may suffer from a large number of binary state variables as a result of the binarization. 

Alternatively, when the expected cost-to-go functions are assumed to be Lipschitz continuous,~\cite{ahmed2020stochastic} propose nonlinear cuts to exactly approximate the nonconvex expected cost-to-go functions via the augmented Lagrangian method. More recently,~\cite{zhang2019stochastic} propose and study a unified theoretical framework for general multistage nonconvex stochastic mixed integer nonlinear programming problems where the expected cost-to-go functions are not necessarily Lipschitz continuous. Specifically, they propose regularized expected cost-to-go functions that can be exactly approximated by the generalized conjugacy cuts, generalizing the results of~\cite{ahmed2020stochastic}. However, despite these theoretical convergence properties, the nonlinear cut generation from augmented Lagrangian and the implementation of a variant of the nested Benders decomposition algorithm that incorporates these cuts can be computationally challenging. 

\cite{fullner2022non} extend the ideas of binarization procedure in~\cite{zou2019stochastic} and the regularization procedure in~\cite{zhang2019stochastic} to develop the so-called nonconvex nested Benders decomposition algorithm for generic multi-stage mixed integer nonconvex nonlinear program. Despite its generality and successful implementation in \emph{deterministic} multi-stage MINLP problems, it is unclear, to the best of our knowledge, whether or not applying this approach to MSILPs with a large-scale scenario tree is computationally feasible. 

In sum, all the works mentioned above aim to directly address the nonconvex (and sometimes non-Lipschitz) expected cost-to-go functions defined at each node of the scenario tree by developing exact lower-bounding techniques. In contrast, we circumvent the challenge of approximating the nonconvex expected cost-to-go functions by relocating all integer state variables (defined in the node subproblems) to the first stage (defined in the root node), so that the resulting cost-to-go functions defined at non-root nodes are convex and can be approximated (exactly) by a decomposition scheme (e.g., nested Benders or SDDP). \Song{To the best of our knowledge, our paper is one of the first that presents the performance of \Merve{stage-wise decomposition based} algorithms for MSILPs with both integer and continuous state variables. In addition, we  showcase the performance of the proposed heuristic decision policies via decision rules, \Merve{2SLDRs}, as well as the lower bounding technique constructed based on \Merve{our} exact algorithm.} We next discuss \Merve{the} proposed reformulation in detail. 

\subsection{A Partially Extended Formulation for MSILP}

We propose a reformulation for MSILP in the form of a partially extended formulation, in which all integer decisions can be considered as the first-stage variables. The advantage of this reformulation is two-fold. First, the resulting cost-to-go functions at each node $n \in \nodesnoroot$ have only continuous variables (including local and continuous state variables), making it amenable to apply standard Benders-type cutting-plane approximation to these convex functions. Second, since integer variables only appear in the first-stage problem, we could apply, e.g., a B\&C procedure by branching on these first-stage integer variables and approximating the expected cost-to-go functions using stage-wise decomposition algorithms such as the SDDP (see Section \ref{sec:methodology}). Specifically, the proposed partially extended (nested) formulation is given by:
\vspace{-0.2em}
\begin{align*}
    \Qref_{\rootnode} = \min\; & 
    \sum_{n\in\nodes}\prob_nc_n^\top \zn 
    + d_\rootnode^\top \xroot + h_\rootnode^\top \yroot + \sum_{n \in \childs(\rootnode)} \pprob_{\rootnode n} \Qref_n(\xroot,z), \tag{$P^R$} \label{mod:reformulation} \\
    \text{s.t.}\; & 
     H_{\rootnode} \zroot \geq g_{\rootnode}, \;
    J_{\rootnode} \xroot \geq f_{\rootnode}, \;
    C_{\rootnode} \xroot + D_{\rootnode} \zroot + E_{\rootnode} \yroot  \geq  b_{\rootnode},  \\
    & H_n \zn \geq  G_n \znparent + g_n, & \forall n \in \nodesnoroot,
\end{align*}
where the cost-to-go function $Q^R_n(\cdot, \cdot)$ associated with 
$n\in \nodes\setminus \{\rootnode\}$ is defined in a  nested fashion as:
\vspace{-0.2em}
\begin{equation*}
    \Qref_n(\xnparent, z) = \min_{\xn\in \R^k, \yn \in \R^r} \left\{d_n^\top \xn + h_n^\top \yn  + \sum_{n' \in \childs(n)} \pprob_{n n'} \Qref_{n'}(\xn, z) \mid \eqref{eq:st_onlyx}-\eqref{eq:st_all}\right\}.
\end{equation*}

The key difference between this reformulation~\eqref{mod:reformulation} and the original formulation~\eqref{mod:scenarioTree} is that all the integer state variables $z = \{\zn\}_{n \in \nodes}$, along with their associated constraint sets~\eqref{eq:st_onlyz}, are now ``moved'' to the problem at the root node $\rootnode$, i.e., they are treated as first-stage variables and constraints. As such, the cost-to-go function associated with each non-root node only includes the state and local continuous variables and their associated constraints. Since the continuous state variables remain to be defined in the individual node-based subproblems through the nested form, we call this reformulation a \textit{partially extended formulation}. The validity of this reformulation is clear because we only relocate the (node-based) integer state variables in order to achieve the desired property of piece-wise linear convex expected cost-to-go functions, \Merve{given its first-stage decisions,} which accommodates the application of stage-wise decomposition methods such as SDDP.

\subsection{An Aggregation Framework via MC-based Structural Policies}
While the proposed partially extended formulation~\eqref{mod:reformulation} has piece-wise linear convex expected cost-to-go functions, \Merve{given its first-stage decisions,} it may lead to exponentially many integer variables in the first stage
. We propose an aggregation framework that imposes structures to the formulation by aggregating the integer state variables to alleviate this challenge while maintaining high-quality solutions. \Merve{In that regard,} the idea is to leverage the structure of the underlying stochastic process to impose a restriction to the problem and obtain an informative decision policy.

In this paper, we consider the case where the underlying stochastic process is modeled as a Markov chain (MC). We choose this specific structure for two reasons. First, many stochastic processes in real-world applications can be modeled or \Merve{well-}approximated using MCs~\citep{lohndorf2019modeling}, including disaster relief logistics planning problems that we focus on in our numerical experiments (see Section~\ref{sec:application}). Second, it facilitates the exposition of our aggregation framework and allows us to construct a wide range of structured policies. We discuss alternative structural policies beyond those based on MC models in Remark~\ref{re:otherStochasticProcesses}.

To start, let $\markov$ be the set of MC states where each $m\in \markov$ is an $\dimmarkov$-dimensional vector (i.e., $m \in \R^\dimmarkov$). Let $\markov_t$ represent the set of MC states that are reachable in stage $t\in \horizon$, starting from a given initial state. Notation $m^t=(m_1,...,m_t)$, where $m_{t'} \in \markov_{t'}$ for $t'\in [t]$, represents a sequence of admissible MC states, that is, the probability to reach $m_{t'+1}$ given $m_{t'}$ is positive for all $t' \in [t-1]$. 

We now describe the connection between the MC states and the scenario tree model used in the MSILP formulation~\eqref{mod:scenarioTree}. Each node $n\in \staget$ ($t\geq 1$) of the scenario tree corresponds to an MC state, $m_t(n) \in \markov_t$, and is uniquely determined by the trajectory of the stochastic process from the root node (initial MC state) to MC state $m_t(n)$, that is, the sequence of MC states $m^t(n) = (m_1(\rootnode), ..., m_{t-1}(a(n)),m_{t}(n))$. To simplify the exposition, we write $m^t$ instead of $m^t(n)$ and $m_t$ instead of $m_t(n)$ when we refer to any node $n \in \staget$.

As stated above, the MC representation of the stochastic process allows us to parameterize each node in the scenario tree by means of their corresponding stage and sequence of MC states, i.e., for each $n \in \staget$, we have $n \equiv (t, m^t)$. We can then rewrite the integer state variables $\zn$ associated with node $n$ as $\zn \equiv \ztM$. Based on this parameterization, our aggregation framework relies on a linear transformation that will compress the full MC history up to stage $t$, namely the vector $m^t$ of size $\dimmarkov\cdot t$, into a vector $\maptmarkov m^t$ of size $q_t \leq \dimmarkov\cdot t$, by means of a transformation matrix $\maptmarkov \in  \Z^{ q_t \times \dimmarkov\cdot t}$. That is, it parametrizes the integer state variables of stage $t$ by the limited information $\maptmarkov m^t$:
\vspace{-0em}
\begin{equation} \label{eq:generalRestricition}
    \zn \equiv \ztM \rightarrow \ztphiM  \qquad  \qquad \forall n \in \staget, t \in \horizon.
\end{equation}
It is clear that this aggregation significantly reduces the number of integer state variables if $q_t \ll \dimmarkov\cdot t$.

Different transformations can impose different structural properties on the discrete decisions and, thus, lead to policies with different levels of aggregation. Table \ref{tab:transformations} summarizes a set of transformations, the resulting aggregated variables, and the conditions under which the original integer state variables associated with nodes $n$ and $n'$ in $\stage_t$ correspond to the same (aggregated) variables, i.e., $\zn$ and $\znprime$ have the same value given by $\ztphiM$. We explain each of these policies below.

\begin{table}[tb]
\aboverulesep=0ex 
\belowrulesep=0ex 
\setlength\extrarowheight{4pt}
    \centering
    \begin{tabular}{c|c|c|c}
        \toprule
        Name & Transformation matrix $\maptmarkov$ & Aggr. var. &  Condition for $\zn \equiv \znprime, \ \forall n,n' \in \staget$ \\
         \midrule
       \hn & $\bzero^{1 \times \dimmarkov \cdot t}$  &  $\za_{t,0}$ &  $t(n) = t(n')$ \\
        \ma & $\left[
\bzero^{\dimmarkov \times \dimmarkov \cdot (t-1)} \ | \ \identity^{\dimmarkov \times \dimmarkov}\right]$  & $\za_{t,m_t} $  & $m_t(n) = m_t(n')$ \\
        \mm & $\left[
\bzero^{2 \dimmarkov \times \dimmarkov \cdot (t-2)} \ | \ \identity^{2 \dimmarkov \times 2 \dimmarkov} 
\right]$  &  $\za_{t,(m_{t-1},m_t)}$  & $ m_t(n) = m_t(n') \ \& \ m_{t-1}(a(n)) = m_{t-1}(a(n'))$ \\
        \pa & $\left[
\bzero^{(\bar{\dimmarkov} + \dimmarkov) \times (\dimmarkov \cdot t - \bar{\dimmarkov} - \dimmarkov)} \ | \ \identity^{(\bar{\dimmarkov} + \dimmarkov) \times (\bar{\dimmarkov} + \dimmarkov)} 
\right] $ & $\za_{t,(\bar{m}_{t-1},m_t)}$ & $ m_t(n) = m_t(n') \ \& \ \bar{m}_{t-1}(a(n)) = \bar{m}_{t-1}(a(n')) $ \\
        \fh & $\identity^{\dimmarkov \cdot t \times  \dimmarkov \cdot t} $ & $\za_{t,m^t}$ & $ m^t(n) = m^t(n')$ \\*[0.1cm]
        \bottomrule
    \end{tabular}
    \caption{MC-based transformations examples, the resulting aggregated variables (Aggr. var.), and the equivalence of the original variables. $\bzero$ and $\identity$ denote the zero and identity matrices, respectively.}  
    \label{tab:transformations}
\end{table}


\Margarita{
The simplest transformation we consider is 
a \textit{Here-and-Now} (\hn) transformation, which aggregates 
integer state variables associated with nodes that are in the same stage
: $\zn \equiv \ztM \rightarrow \za_{t,0}$ for all $n \in \staget$.  Thus, the corresponding 
\Merve{integer state variables}
will be independent of the realization of the stochastic process, 
\Merve{i.e., follow a} static policy 
(see, e.g., \cite{basciftci2019adaptive}).

An alternative is to leverage the structure of the MC to create a transformation that leads to more informative policies. For example, the \textit{Markovian} (\ma) transformation aggregates integer state variables associated with nodes that share the same stage and the same MC state: $\zn \equiv \ztM \rightarrow \za_{t,m_t}$ for all $n \in \staget$. Another option is to retain information from previous stages at the price of a larger number of aggregated variables. One such example is what we refer to as the \textit{Double Markovian} (\mm) transformation, which goes one step beyond \ma\ and aggregates integer state variables associated with nodes that share the same MC states in the current and previous stage. Following this pattern, the extreme case is to consider the entire history of the MC states starting from the initial MC state, referred to as the \textit{Full History} (\fh) transformation. Note that  \fh\ leads to the original partially extended formulation for MSILP~\eqref{mod:reformulation}, since no integer state variables are effectively aggregated.

Although all the aforementioned MC-based transformations are constructed based on the complete MC state information, we may also consider partial information from MC states. Consider a subset of attributes of size $\bar{\dimmarkov} < \dimmarkov$ associated with MC state $m_t$, and assume without loss of generality that these attributes of interest correspond to the last $\bar{\dimmarkov}$ attributes in the MC state vector, i.e., we represent an MC state as $m_t = (\cdot, \bar{m}_t)\in \R^{\dimmarkov -\bar{\dimmarkov}}\times \R^{\bar{\dimmarkov}}$. Then, we can create a transformation based only on this partial vector $\bar{m}_t$ associated with the MC state $m_t$. As an example, we introduce the  \textit{Partial Markovian} (\pa) transformation that considers complete MC state information in the current stage but only partial information from the previous stage.  
}

\begin{figure}[tbp]
    \centering
    \scalebox{1}{
    \begin{tikzpicture}[->,>=stealth',shorten >=1pt,auto,node distance=1cm,
thick]       
\node[node light] (r) at (0,0) {};
\node[node dark] (n11)  at (0.8,0.8)  {};
\node[node light] (n12)  at (0.8,-0.8)  {};

\node[node dark] (n21)  at (1.6,1.2)  {};
\node[node light] (n22)  at (1.6,0.4)  {};
\node[node dark] (n23)  at (1.6,-0.4)  {};
\node[node light] (n24)  at (1.6,-1.2)  {};

\node[node dark] (n31)  at (2.4,1.4)  {};
\node[node light] (n32)  at (2.4,1)  {};
\node[node dark] (n33)  at (2.4,0.6)  {};
\node[node light] (n34)  at (2.4,0.2)  {};
\node[node dark] (n35)  at (2.4,-0.2)  {};
\node[node light] (n36)  at (2.4,-0.6)  {};
\node[node dark] (n37)  at (2.4,-1)  {};
\node[node light] (n38)  at (2.4,-1.4)  {};

\node[] (name) at (1.3, 2) {Scenario Tree};

\node[text node] (name) at (-0.1, -2) {$t=1$};
\node[text node] (name) at (0.8, -2) {$2$};
\node[text node] (name) at (1.6, -2) {$3$};
\node[text node] (name) at (2.4, -2) {$4$};

\path[every node/.style={font=\sffamily\small}]
(r) 
edge[arc] node [left] {} (n11)
edge[arc] node [left] {} (n12)
(n11)
edge[arc] node [left] {} (n21)
edge[arc] node [left] {} (n22)
(n12)
edge[arc] node [left] {} (n23)
edge[arc] node [left] {} (n24)
(n21)
edge[arc] node [left] {} (n31)
edge[arc] node [left] {} (n32)
(n22)
edge[arc] node [left] {} (n33)
edge[arc] node [left] {} (n34)
(n23)
edge[arc] node [left] {} (n35)
edge[arc] node [left] {} (n36)
(n24)
edge[arc] node [left] {} (n37)
edge[arc] node [left] {} (n38)
;

\end{tikzpicture}
    \hspace{2em}
    \begin{tikzpicture}[->,>=stealth',shorten >=1pt,auto,node distance=1cm,
thick]       
\node[node light] (r) at (0,0) {};
\node[node-square, draw=blue!70!white, fill=blue!70!white] (n11)  at (0.8,0.8)  {};
\node[node-square, draw=blue!70!white, fill=blue!70!white] (n12)  at (0.8,-0.8)  {};

\node[node-polygon5, draw=violet!70!white, fill=violet!70!white] (n21)  at (1.6,1.2)  {};
\node[node-polygon5, draw=violet!70!white, fill=violet!70!white] (n22)  at (1.6,0.4)  {};
\node[node-polygon5, draw=violet!70!white, fill=violet!70!white] (n23)  at (1.6,-0.4)  {};
\node[node-polygon5, draw=violet!70!white, fill=violet!70!white] (n24)  at (1.6,-1.2)  {};

\node[node-polygon6, draw=cyan!80!white, fill=cyan!80!white] (n31)  at (2.4,1.4)  {};
\node[node-polygon6, draw=cyan!80!white, fill=cyan!80!white] (n32)  at (2.4,1)  {};
\node[node-polygon6, draw=cyan!80!white, fill=cyan!80!white] (n33)  at (2.4,0.6)  {};
\node[node-polygon6, draw=cyan!80!white, fill=cyan!80!white] (n34)  at (2.4,0.2)  {};
\node[node-polygon6, draw=cyan!80!white, fill=cyan!80!white] (n35)  at (2.4,-0.2)  {};
\node[node-polygon6, draw=cyan!80!white, fill=cyan!80!white] (n36)  at (2.4,-0.6)  {};
\node[node-polygon6, draw=cyan!80!white, fill=cyan!80!white] (n37)  at (2.4,-1)  {};
\node[node-polygon6, draw=cyan!80!white, fill=cyan!80!white] (n38)  at (2.4,-1.4)  {};

\node[] (name) at (1.1, 2) {\hn};

\node[text node] (name) at (-0.1, -2) {$t=1$};
\node[text node] (name) at (0.8, -2) {$2$};
\node[text node] (name) at (1.6, -2) {$3$};
\node[text node] (name) at (2.4, -2) {$4$};

\path[every node/.style={font=\sffamily\small}]
(r) 
edge[arc] node [left] {} (n11)
edge[arc] node [left] {} (n12)
(n11)
edge[arc] node [left] {} (n21)
edge[arc] node [left] {} (n22)
(n12)
edge[arc] node [left] {} (n23)
edge[arc] node [left] {} (n24)
(n21)
edge[arc] node [left] {} (n31)
edge[arc] node [left] {} (n32)
(n22)
edge[arc] node [left] {} (n33)
edge[arc] node [left] {} (n34)
(n23)
edge[arc] node [left] {} (n35)
edge[arc] node [left] {} (n36)
(n24)
edge[arc] node [left] {} (n37)
edge[arc] node [left] {} (n38)
;

\end{tikzpicture}
    \hspace{2em}
    \begin{tikzpicture}[->,>=stealth',shorten >=1pt,auto,node distance=1cm,
thick]       
\node[node light] (r) at (0,0) {};
\node[node-square, draw=blue!70!white, fill=blue!70!white] (n11)  at (0.8,0.8)  {};
\node[node-square, draw=blue!40!white, fill=blue!40!white] (n12)  at (0.8,-0.8)  {};

\node[node-polygon5, draw=violet!70!white, fill=violet!70!white] (n21)  at (1.6,1.2)  {};
\node[node-polygon5, draw=violet!40!white, fill=violet!40!white] (n22)  at (1.6,0.4)  {};
\node[node-polygon5, draw=violet!70!white, fill=violet!70!white] (n23)  at (1.6,-0.4)  {};
\node[node-polygon5, draw=violet!40!white, fill=violet!40!white] (n24)  at (1.6,-1.2)  {};

\node[node-polygon6, draw=cyan!80!white, fill=cyan!80!white] (n31)  at (2.4,1.4)  {};
\node[node-polygon6, draw=cyan!40!white, fill=cyan!40!white] (n32)  at (2.4,1)  {};
\node[node-polygon6, draw=cyan!80!white, fill=cyan!80!white] (n33)  at (2.4,0.6)  {};
\node[node-polygon6, draw=cyan!40!white, fill=cyan!40!white] (n34)  at (2.4,0.2)  {};
\node[node-polygon6, draw=cyan!80!white, fill=cyan!80!white] (n35)  at (2.4,-0.2)  {};
\node[node-polygon6, draw=cyan!40!white, fill=cyan!40!white] (n36)  at (2.4,-0.6)  {};
\node[node-polygon6, draw=cyan!80!white, fill=cyan!80!white] (n37)  at (2.4,-1)  {};
\node[node-polygon6, draw=cyan!40!white, fill=cyan!40!white] (n38)  at (2.4,-1.4)  {};

\node[] (name) at (1.1, 2) {\ma};

\node[text node] (name) at (-0.1, -2) {$t=1$};
\node[text node] (name) at (0.8, -2) {$2$};
\node[text node] (name) at (1.6, -2) {$3$};
\node[text node] (name) at (2.4, -2) {$4$};

\path[every node/.style={font=\sffamily\small}]
(r) 
edge[arc] node [left] {} (n11)
edge[arc] node [left] {} (n12)
(n11)
edge[arc] node [left] {} (n21)
edge[arc] node [left] {} (n22)
(n12)
edge[arc] node [left] {} (n23)
edge[arc] node [left] {} (n24)
(n21)
edge[arc] node [left] {} (n31)
edge[arc] node [left] {} (n32)
(n22)
edge[arc] node [left] {} (n33)
edge[arc] node [left] {} (n34)
(n23)
edge[arc] node [left] {} (n35)
edge[arc] node [left] {} (n36)
(n24)
edge[arc] node [left] {} (n37)
edge[arc] node [left] {} (n38)
;

\end{tikzpicture}
    \hspace{2em}
    \begin{tikzpicture}[->,>=stealth',shorten >=1pt,auto,node distance=1cm,
thick]       
\node[node light] (r) at (0,0) {};
\node[node-square, draw=blue!70!white, fill=blue!70!white] (n11)  at (0.8,0.8)  {};
\node[node-square, draw=blue!40!white, fill=blue!40!white] (n12)  at (0.8,-0.8)  {};

\node[node-polygon5, draw=violet!80!black, fill=violet!80!black] (n21)  at (1.6,1.2)  {};
\node[node-polygon5, draw=violet!70!white, fill=violet!70!white] (n22)  at (1.6,0.4)  {};
\node[node-polygon5, draw=violet!40!white, fill=violet!40!white] (n23)  at (1.6,-0.4)  {};
\node[node-polygon5, draw=violet!20!white, fill=violet!20!white] (n24)  at (1.6,-1.2)  {};

\node[node-polygon6, draw=cyan!90!black, fill=cyan!90!black] (n31)  at (2.4,1.4)  {};
\node[node-polygon6, draw=cyan!70!white, fill=cyan!70!white] (n32)  at (2.4,1)  {};
\node[node-polygon6, draw=cyan!40!white, fill=cyan!40!white] (n33)  at (2.4,0.6)  {};
\node[node-polygon6, draw=cyan!10!white, fill=cyan!10!white] (n34)  at (2.4,0.2)  {};
\node[node-polygon6, draw=cyan!90!black, fill=cyan!90!black] (n35)  at (2.4,-0.2)  {};
\node[node-polygon6, draw=cyan!70!white, fill=cyan!70!white] (n36)  at (2.4,-0.6)  {};
\node[node-polygon6, draw=cyan!40!white, fill=cyan!40!white] (n37)  at (2.4,-1)  {};
\node[node-polygon6, draw=cyan!10!white, fill=cyan!10!white] (n38)  at (2.4,-1.4)  {};

\node[] (name) at (1.1, 2) {\mm};

\node[text node] (name) at (-0.1, -2) {$t=1$};
\node[text node] (name) at (0.8, -2) {$2$};
\node[text node] (name) at (1.6, -2) {$3$};
\node[text node] (name) at (2.4, -2) {$4$};

\path[every node/.style={font=\sffamily\small}]
(r) 
edge[arc] node [left] {} (n11)
edge[arc] node [left] {} (n12)
(n11)
edge[arc] node [left] {} (n21)
edge[arc] node [left] {} (n22)
(n12)
edge[arc] node [left] {} (n23)
edge[arc] node [left] {} (n24)
(n21)
edge[arc] node [left] {} (n31)
edge[arc] node [left] {} (n32)
(n22)
edge[arc] node [left] {} (n33)
edge[arc] node [left] {} (n34)
(n23)
edge[arc] node [left] {} (n35)
edge[arc] node [left] {} (n36)
(n24)
edge[arc] node [left] {} (n37)
edge[arc] node [left] {} (n38)
;

\end{tikzpicture}
    }
    \caption{An illustrative example showing nodes with the same aggregated integer state variables for transformations \hn, \ma, and \mm. Nodes with the same color and shape share the same set of integer state variables.}
    \label{fig:example_aggregation}
\end{figure}
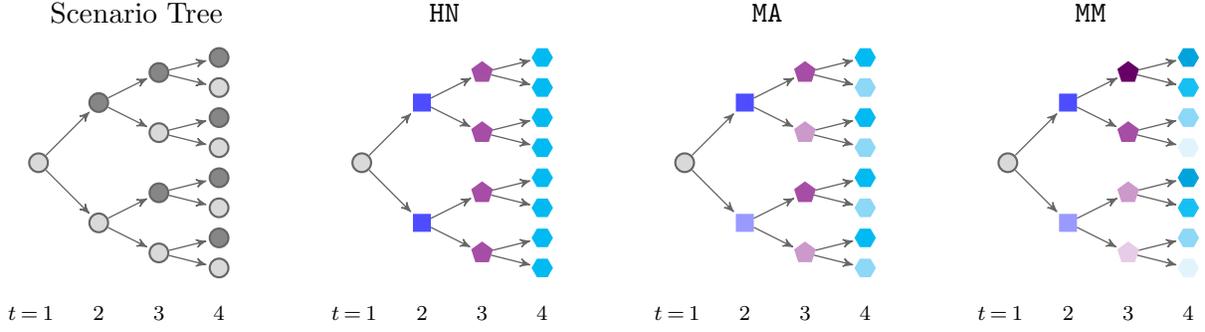

\begin{example} \label{example:scenario_tree}
Consider an MC $\markov$ with two possible states (i.e., light \mcLight{1}\ and dark \mcDark{1}) and positive transition probabilities for all transitions between states. The first drawing in Figure \ref{fig:example_aggregation} depicts a scenario tree built for $\markov$ starting at a light state and representing all possible MC states for $\periods=4$ stages, where the color of each node (light or dark) represents the corresponding MC state. The remaining drawings in Figure \ref{fig:example_aggregation} illustrate the effect of aggregation on the nodes of the scenario tree for transformations \hn, \ma, and \mm\ (we omit \pa\ because in this example the MC state vector has only one attribute). The \hn\ tree shows that all nodes in the same stage will share the same set of (aggregated) integer state variables $\za$ (i.e., they have the same color and shape). In contrast, \ma\ aggregates integer state variables in nodes with the same stage and MC state. Thus, the \ma\ tree shows nodes of two different colors in each stage, which correspond to two separate sets of aggregated integer state variables for each stage (e.g., the transformation for node $n\equiv (t,m^t) = (3,\mcLight{1}\mcDark{1}\mcLight{1})\in \nodes_3$ is $\Phi_3 m^3(\mcLight{1}\mcDark{1}\mcLight{1}) = m_3(\mcLight{1}) = \pentaLightDark$, thus,  $\za_{3,\pentaLightDark}$ represents its associated aggregated variables). 
Finally, nodes in each stage of the \mm\ tree have up to four different colors, as a result of an aggregation based on the MC states in both the current stage and the previous stage 
(e.g., $\za_{3,\pentaDarkDark}$ is associated to node (3,\mcLight{1}\mcDark{1}\mcDark{1})  because $\Phi_3 m^3(\mcLight{1}\mcDark{1}\mcDark{1}) = (m_2(\mcDark{1}), m_3(\mcDark{1})) = \pentaDarkDark$). 
Overall, \hn\ contains four sets of $\za$ variables (i.e., one for each stage), \ma \  seven, and \mm \  eleven, while the \fh\ variant contains 15 sets of $\za$ variables (i.e., one for each node of the scenario tree). 
\hfill $\square$
\end{example}

\begin{remark}
These transformations are just a few examples of a wide range of aggregations that can be constructed within our framework. One key consideration in constructing these aggregations is that there is a clear trade-off between the quality of the underlying decision policy and the corresponding computational effort. Aggregations that include (preserve) more information would naturally lead to better decision policies but also make the resulting problem harder to solve due to the increased number of integer state variables in the first stage. Choosing the right transformation depends on the specific problem type and the available computational resources. We leave the problem of adaptively finding the best transformation as a future research direction.  
\end{remark}

\begin{remark}\label{re:otherStochasticProcesses}
The main advantage of assuming that the stochastic process is defined by an MC is that most stage-dependent processes can be approximated by an MC using a discretization technique with high accuracy \citep{bally2003quantization,pages2004optimal}. Also, the work by \cite{lohndorf2019modeling} shows that such an approximation can lead to high-quality solutions compared to the ones obtained by approximating the stochastic process with autoregressive models in certain applications. We also note that our approach is applicable to stage-wise independent stochastic processes. In these cases, one can approximate the probability distribution of a random variable by partitioning its support into a set of clusters and creating aggregations based on these clusters. For example, considering an MSILP with stochastic demand parameters, we may define three levels of demand (low, medium, and high) and create aggregation schemes based on each demand level.
\end{remark}

\subsection{Aggregated MSILP}
We now present the aggregated MSILP formulation after applying a transformation to the integer state variables $z$ via matrices $\maptmarkov, t\in \horizon$. The new set of integer state variables is given by $\za_t\in \Z^{\ell \cdot q_t}$ for each $t\in \horizon$. We employ a mapping function $\phi_t:\nodes_t\rightarrow \{t\}\times[q_t]$ to relate a node $n\in \staget$, which was originally used to index variables $z$, to the index corresponding to its aggregated variables, i.e., $\phi_t(n) = (t, \maptmarkov m^t)$ for each node $n \in \staget$. Thus, the aggregated MSILP can be written as follows:
\begin{subequations}
\begin{align}
    \Qaggre_{\rootnode} = \min\; &  \sum_{n\in \nodes}\prob_n c_n^\top \zphi  + d_\rootnode^\top \xroot + h_\rootnode^\top \yroot + \sum_{n \in \childs(\rootnode)}\pprob_{\rootnode n} \Qaggre_n(\xroot,\za) \tag{$P^A$} \label{mod:aggregated} \\
    \text{s.t. }\; & 
     H_{\rootnode} \zroot \geq g_{\rootnode}, \;
    J_{\rootnode} \xroot \geq f_{\rootnode}, \;
    C_{\rootnode} \xroot + D_{\rootnode} \zroot + E_{\rootnode} \yroot  \geq  b_{\rootnode},  \label{aggre:root}\\
    & H_n \zphi \geq  G_n \zphiparent + g_n, & \forall n \in \nodesnoroot,   \label{aggre:zphi}\\
    & \xroot \in \R^k, \yroot \in \R^r, \za \in  \Z^{\ell \cdot \sum_{t \in \horizon} q_t},  \label{aggre:vars}
\end{align}
\end{subequations}
where the cost-to-go function $\Qaggre_n(\cdot,\cdot)$ associated with
$n\in \nodes\setminus \{\rootnode\}$ is defined in a  nested fashion as:
\begin{align*}
    \Qaggre_n(\xnparent, \za) = \min\; & d_n^\top \xn + h_n^\top \yn + \sum_{n' \in \childs(n)}\pprob_{n n'} \Qaggre_{n'}(\xn, \za)\\
    \text{s.t.}\; & J_n \xn \geq  F_n \xnparent + f_n, \\
    & C_n \xn + D_n \zphi + E_n \yn  \geq A_n \xnparent + B_n  \zphiparent + b_n, \\
    & \xn \in \R^k, \yn \in \R^r.
\end{align*}

We note that any feasible solution of \eqref{mod:aggregated} is also  feasible to \eqref{mod:scenarioTree} with the same objective value, because the only difference between \eqref{mod:aggregated} and \eqref{mod:scenarioTree} is that \eqref{mod:aggregated} enforces some integer variables of \eqref{mod:scenarioTree} to take the same values. However, the reverse is not true in general, and hence an optimal solution of \eqref{mod:aggregated} might be sub-optimal for \eqref{mod:scenarioTree}. 

Although the aggregated problem~\eqref{mod:aggregated} potentially has much fewer integer variables than~\eqref{mod:reformulation}, depending on the transformation $\maptmarkov$, it may remain  computationally challenging for an exact solution approach due to the existence of (aggregated) first-stage integer variables and the continuous state variables defined at each node of the scenario tree. \Margarita{We next introduce an exact solution approach based on a variant of the SDDP algorithm, and then present alternative approximation approaches to generate feasible decision policies and dual bounds that are computationally efficient.}

\section{Solution Methodology} \label{sec:methodology}

We present two solution methods for the aggregated MSILP formulation~\eqref{mod:aggregated}. Both are decomposition methods implemented as B\&C algorithms, taking advantage of the fact that all the integer variables in~\eqref{mod:aggregated} are in the first stage. The first one is an exact method that employs a variant of the SDDP algorithm adapted to our setting as a subroutine. The second method employs 2SLDR to impose additional solution structures that lead to a two-stage approximation to the 
\Merve{MSILP}.

\subsection{SDDP Integrated B\&C Algorithm}
\label{sec:sddp}

We now present an exact algorithm for solving formulation~\eqref{mod:aggregated} by integrating the SDDP algorithm within a B\&C framework. The high-level idea is to decompose the problem into a master problem and a set of subproblems, where the master problem is concerned about the first-stage decisions including all the integer variables of~\eqref{mod:aggregated}, and the subproblems deal with the remaining stages, which involve continuous state and local variables. 

The algorithm performs a branch-and-bound (B\&B) procedure for the master problem and enters a cut generation procedure every time an integer relaxation solution is encountered (e.g., done through the so-called callback functions in modern IP solvers such as CPLEX). To generate a cut, following the basic idea of Benders decomposition, we need to solve a subproblem associated with each child node of $\rootnode$ in the scenario tree, which deals with the remaining stages of~\eqref{mod:aggregated} via a variant of the SDDP algorithm for a given candidate solution $(\xrootsol, \zasol, \thetasol)$ from the master problem. \Song{Note that we adapt the standard SDDP algorithm---which is applicable when we have stage-wise independence---to our aggregated MSILP, taking advantage of the MC structure in the stochastic process \Merve{(representing stagewise-dependent uncertainty)}.} In particular, we consider the so-called multi-cut version, where a variable $\theta_n$ is introduced for each $n \in \childs(\rootnode)$ to represent the value of the 
outer approximation of the cost-to-go function at node $n$. Thus, the master problem is given by:
\begin{subequations}\label{SDDP-master}
\begin{align}
    \Qsub_{\rootnode} = \min& \sum_{n\in \nodes}\prob_n c_n^\top \zphi  + d_\rootnode^\top \xroot + h_\rootnode^\top \yroot + \sum_{n \in \childs(\rootnode)}\pprob_{\rootnode n} \theta_{n} \\
    \text{s.t.}\ & \eqref{aggre:root}-\eqref{aggre:vars}, \nonumber \\
    & \theta_n \geq \alpha_n^\top \xroot + \beta_n^\top \za + \gamma_n, \quad \forall  (\alpha_n,\beta_n,\gamma_n)\in \coeffs_n, \ n \in \childs(\rootnode),
\end{align}
\end{subequations}
where $\coeffs_n$ stores the coefficients associated with all the Benders cuts constructed to approximate $\Qaggre_n(\xroot,\za)$ during the solution procedure. 
We refer the reader to Appendix \ref{app:bc_sddp} for a detailed description of the overall B\&C procedure, including a pseudo-code shown in  Algorithm \ref{alg:bc-sddp}. 

We next explain the implementation details of this B\&C framework integrated with SDDP. We first describe the SDDP subproblem structure and then present the entire SDDP sub-routine.

\subsubsection{SDDP Subproblems.}

A key feature of the SDDP algorithm is its capability to leverage the structure of the underlying stochastic process in defining the cost-to-go functions and the corresponding SDDP subproblems. For example, we only need to create one SDDP subproblem per stage in the case of stage-wise independence \citep{pereira1991multi}, and one subproblem per stage and MC state in the case when the stage-wise dependence structure is modeled by an MC \citep{bonnans2012energy,philpott2012dynamic}. In what follows, we describe the SDDP subproblems for formulation~\eqref{mod:aggregated} based on the MC structure and the employed transformation. 

First, we observe that if the first-stage integer variables $\za$ do not exist, instead of defining one subproblem for each node $n\in \nodes$, we can have one subproblem defined for each stage $t\in \horizon$ and MC state $m_t\in \markov_t$, that is, all nodes $n \in \staget$ with $m_t(n)= m_t$ share the same cost-to-go function \citep{philpott2012dynamic}, thanks to the Markovian nature of the stochastic process. However, this might not be the case when we consider the aggregated variables $\za$ because 
cost-to-go functions with the same MC state might be associated with different aggregated variables depending on the chosen transformation. In such a case, the cost-to-go functions are defined over two different sets of variables. To see this, consider two nodes $n,n'\in \staget$ with the same MC state (i.e., $m_t(n) = m_t(n') = m_t$) but with different sets of (aggregated) integer variables $\za$ (i.e., $\phi_t(n)\neq \phi_t(n')$), then we need to define two separate SDDP subproblems for $n$ and $n'$.
In other words, we can only define the same SDDP subproblem for two nodes $n,n'\in \staget$ if they correspond to the same MC state (i.e., $m_t(n)=m_t(n')$) \emph{and} the same aggregated variables (i.e., $\phi_t(n) = \phi_t(n')$). In the following example we use a policy graph---a concept introduced by \cite{dowson2020policy}---that represents the relationship between SDDP subproblems. Note that the policy graphs for \eqref{mod:aggregated} are constructed based on the MC transition graph and the chosen transformation.

\begin{figure}[tbp]
\centering
    \scalebox{1}{
    \begin{tikzpicture}[->,>=stealth',shorten >=1pt,auto,node distance=1cm,
thick]       
\node[node light] (r) at (0,0) {};
\node[node-square, draw=blue!70!white, fill=blue!70!white] (n11)  at (1,0.8)  {};
\node[node-square, draw=blue!40!white, fill=blue!40!white] (n12)  at (1,-0.8)  {};

\node[node-polygon5, draw=violet!70!white, fill=violet!70!white] (n21)  at (2,0.8)  {};
\node[node-polygon5, draw=violet!40!white, fill=violet!40!white] (n22)  at (2,-0.8)  {};

\node[node-polygon6, draw=cyan!80!white, fill=cyan!80!white] (n31)  at (3,0.8)  {};
\node[node-polygon6, draw=cyan!40!white, fill=cyan!40!white] (n32)  at (3,-0.8)  {};

\node[] (name) at (1.5, 1.2) {\hn/\ma};

\node[text node] (name) at (-0.1, -1.2) {$t=1$};
\node[text node] (name) at (1, -1.2) {$2$};
\node[text node] (name) at (2, -1.2) {$3$};
\node[text node] (name) at (3, -1.2) {$4$};

\path[every node/.style={font=\sffamily\small}]
(r) 
edge[arc] node [left] {} (n11)
edge[arc] node [left] {} (n12)
(n11)
edge[arc] node [left] {} (n21)
edge[arc] node [left] {} (n22)
(n12)
edge[arc] node [left] {} (n21)
edge[arc] node [left] {} (n22)
(n21)
edge[arc] node [left] {} (n31)
edge[arc] node [left] {} (n32)
(n22)
edge[arc] node [left] {} (n31)
edge[arc] node [left] {} (n32)
;

\end{tikzpicture} 
    \hspace{5em}
    \begin{tikzpicture}[->,>=stealth',shorten >=1pt,auto,node distance=1cm,
thick]       
\node[node light] (r) at (0,0) {};
\node[node-square, draw=blue!70!white, fill=blue!70!white] (n11)  at (1,0.8)  {};
\node[node-square, draw=blue!40!white, fill=blue!40!white] (n12)  at (1,-0.8)  {};

\node[node-polygon5, draw=violet!80!black, fill=violet!80!black] (n21)  at (2,0.8)  {};
\node[node-polygon5, draw=violet!70!white, fill=violet!70!white] (n22)  at (2,0.25)  {};
\node[node-polygon5, draw=violet!40!white, fill=violet!40!white] (n23)  at (2,-0.25)  {};
\node[node-polygon5, draw=violet!20!white, fill=violet!20!white] (n24)  at (2,-0.8)  {};

\node[node-polygon6, draw=cyan!70!black, fill=cyan!70!black] (n31)  at (3,0.8)  {};
\node[node-polygon6, draw=cyan!80!white, fill=cyan!80!white] (n32)  at (3,0.25)  {};
\node[node-polygon6, draw=cyan!40!white, fill=cyan!40!white] (n33)  at (3,-0.25)  {};
\node[node-polygon6, draw=cyan!20!white, fill=cyan!20!white] (n34)  at (3,-0.8)  {};

\node[] (name) at (1.5, 1.2) {\mm};

\node[text node] (name) at (-0.1, -1.2) {$t=1$};
\node[text node] (name) at (1, -1.2) {$2$};
\node[text node] (name) at (2, -1.2) {$3$};
\node[text node] (name) at (3, -1.2) {$4$};

\path[every node/.style={font=\sffamily\small}]
(r) 
edge[arc] node [left] {} (n11)
edge[arc] node [left] {} (n12)
(n11)
edge[arc] node [left] {} (n21)
edge[arc] node [left] {} (n22)
(n12)
edge[arc] node [left] {} (n23)
edge[arc] node [left] {} (n24)
(n21)
edge[arc] node [left] {} (n31)
edge[arc] node [left] {} (n32)
(n22)
edge[arc] node [left] {} (n33)
edge[arc] node [left] {} (n34)
(n23)
edge[arc] node [left] {} (n31)
edge[arc] node [left] {} (n32)
(n24)
edge[arc] node [left] {} (n33)
edge[arc] node [left] {} (n34)
;

\end{tikzpicture}
    }
    \caption{SDDP policy graphs and subproblems for various types of transformations, \hn, \ma, and \mm.}
    \label{fig:sddp_sub}
\end{figure}
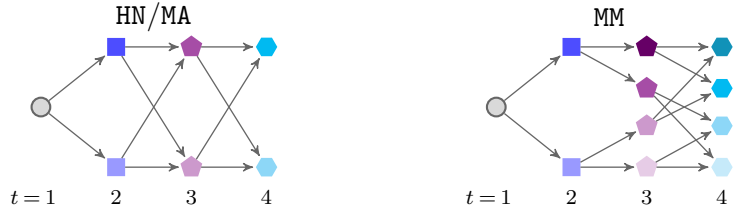

\begin{example}\label{exa:sddp_subproblems}
Figure~\ref{fig:sddp_sub} illustrates different policy graphs for the  \hn, \ma, and \mm\ transformations based on the illustrative example shown in Figure~\ref{fig:example_aggregation}. 
Recall that the underlying stochastic process is defined by an MC with two states (i.e., light \mcLight{1}\ and dark \mcDark{1}). Since at a minimum we should have one SDDP subproblem defined for each MC state, we have at least two SDDP subproblems defined in each stage $t\geq 2$. In the case of \ma, since nodes with the same MC state share the same set of aggregated integer state variables, we have exactly two SDDP subproblems per stage (i.e., one per MC state) and the policy graph only depends on the MC transitions (i.e., a Markovian policy graph). The same applies for \hn, although we only have one set of aggregated integer state variables per stage. Since their corresponding policy graphs are identical, we use \hn/\ma\ to denote the first policy graph in Figure~\ref{fig:sddp_sub}. On the other hand, the SDDP subproblem structure for \mm\ is quite different: nodes with the same MC state can be associated with different aggregated integer state variables, because the \mm\ transformation depends not only on the MC state in the current node but also its parent node. For example, nodes (3,\mcLight{1}\mcDark{1}\mcDark{1}) and (3,\mcLight{1}\mcLight{1}\mcDark{1}) are associated with the same stage and MC state ($m_{3}(\mcLight{1}\mcDark{1}\mcDark{1}) = m_3(\mcLight{1}\mcLight{1}\mcDark{1}) = \mcDark{1}$, see Figure~\ref{fig:example_aggregation}), but since the MC states of their parent nodes differ, they have different sets of aggregated integer state variables (i.e., $\za_{3,\pentaDarkDark}$ and $\za_{3,\pentaLightDark}$, respectively). Therefore, \mm\ has one SDDP subproblem for each pair of MC states, corresponding to both the current node and the parent node for stages $t\geq 3$. The policy graph for \mm\ can be seen as an aggregated version of its scenario tree (see Figure \ref{fig:example_aggregation}) with one copy for each type of nodes. \hfill $\square$
\end{example}

To formally define the SDDP subproblem for \eqref{mod:aggregated}, we represent the SDDP subproblem corresponding to node $n\in \staget$ in the scenario tree by parameterizing it with mapping $\sub(n) := (t, m_{t}(n), \phi_t(n))$. Consider $\subSett$ to be the set of all SDDP subproblems defined for stage $t\geq 2$. The number of subproblems in stage $t=2$ is always the same (i.e., $|\subSet_2| = \stage_2$)  since each node in that stage is associated to a different MC state, whereas the number of subproblems in stages $t > 2$ depends on the transformation employed. For a given SDDP subproblem $\sub\in \subSett$, we define $\subnodest(\sub)$ as the set of nodes associated with $\sub$, i.e., all nodes $n \in \staget$ such that $\sub(n) = \sub$. 

To represent the relationship between SDDP subproblems illustrated by a policy graph, let $\childsub(\sub)$ be the set of children subproblems for each subproblem $\sub$ in the policy graph. Specifically, $\sub'$ is a child of $\sub$  if $\sub'\in \subSet_{t+1}$ and there exists a node $n' \in \subnodes_{t+1}(\sub')$ such that its parent is associated with $\sub$ (e.g., the children set of \squareDark{1} in the \mm\ policy graph is $\childsub(\squareDark{1})=\{\pentaDarkDark, \pentaDark{1}\}$, as shown Figure \ref{fig:sddp_sub}). The transition probability from $\sub\in \subSett$ to one of its children $\sub' \in \childsub(\sub)$ is given by the MC transition probability, that is, $\pprob_{\sub\sub'}= \prob_{a(n')n'}$ for any node $n' \in \subnodes_{t+1}(\sub')$ with parent $a(n') \in \subnodest(\sub)$. Lastly, since there is at least one SDDP subproblem defined for every MC state, the random coefficients involved in the SDDP subproblems can be parameterized by the corresponding MC states, for example, the random cost vector in $\sub\in \subSett$ is defined as $d_\sub = d_n$ for any $n \in \subnodest(\sub)$.

To preserve the stage-wise dependency structure between SDDP subproblems in two consecutive stages, we introduce local copies of the first-stage integer state variables $\za$ in the subproblems. In other words, we copy the first-stage integer solution from an SDDP subproblem to its children subproblems to ensure that each subproblem has the needed information given by $\za$. We note that the use of local copies is common in other SDDP variants to maintain the stage-wise decomposition structure of the problem (see, e.g.,~\cite{zou2019stochastic}). Then, the SDDP subproblem $\subn:=\sigma(n)$ for any node $n\in \staget$ with $t>1$ is given by: 
\vspace{-0em}
\begin{subequations}
\begin{align}
    \Qsub_{\subn}(\xnparentsub, \zcopyparent) = \min\; & d_{\subn}^\top \xnsub + h_{\subn}^\top \ynsub + \sum_{\subn' \in \childsub(\subn)}\pprob_{\subn \subn'} \theta_{\subn'}  \nonumber \\
    \mbox{s.t.} \; & J_\sub \xnsub \geq  F_\sub \xnparentsub + f_\sub,  \\
    &\zcopyn = \zcopyparent,  \label{eq:sub_zcopy} \\ 
    & C_\sub \xnsub + D_\sub \zcopyphin + E_\sub \ynsub  \geq A_\sub \xnparentsub + B_\sub  \zcopyphiparent + b_\sub, \label{eq:sub_allvars} \\ 
    & \theta_{\subn'} \geq \alpha_{\subn'}^\top \xnsub + \beta_{\subn'}^\top \zcopyn + \gamma_{\subn'}, \qquad \forall  (\alpha_{\subn'},\beta_{\subn'},\gamma_{\subn'})\in \coeffs_{\subn'}, \ \subn' \in \childsub(\subn), \label{eq:sddp_cut}
\end{align}
\end{subequations}
where $\zcopyn$ denotes the continuous local copy of the first-stage integer state variables $\za$, which comes from its parent subproblem $\zcopyparent$ via inequality \eqref{eq:sub_zcopy}.  We note that since $\Qsub_{\rootnode}$ contains the original first-stage integer variables $\za$, the SDDP subproblems for stage $t = 2$ directly receive these variables instead of making local copies from the root (i.e., $\Qsub_{\sub}(\xnparentsub, \za)$ for each $\sub \in \childsub(\rootnode)$, see Figure \ref{fig:sddp_b&c}). A variable $\theta_{\sub}$ is introduced for each SDDP subproblem $\sub$ to represent the outer approximation of the cost-to-go function, and a set $\coeffs_{\subn}$ is introduced to store the coefficients associated with all the cutting planes constructed in this approximation during the solution procedure. 

Finally,  $\Qsub_{\subn}$ requires information from the parent of $\subn$ (i.e., $\parent(\sub)$), however, an SDDP subproblem can have multiple parents in the policy graph, (e.g., both \pentaLightDark\ and \pentaDarkDark\ are parents of \hexaDarkDark, see Figure \ref{fig:sddp_sub}). This ambiguity is not an issue for our SDDP implementation as we only consider one scenario at a time in the SDDP forward pass. So each time $\Qsub_{\subn}$ is solved, a specific node $n \in \subnodest(\sub)$ is considered as its parent node, and, thus, the parent of the SDDP subproblem would be $\parent(\subn) = \sub(\parent(n))$.

\subsubsection{SDDP Algorithm.} 
Algorithm \ref{alg:sddp} describes the SDDP sub-routine implemented in our B\&C algorithm. The procedure takes a candidate solution $(\xrootsol,\zasol,\thetasol)$ from the master problem~\eqref{SDDP-master} and one child node of the root $n'\in \childs(\rootnode)$, and runs the SDDP algorithm until we find a valid inequality that cuts off $(\xrootsol,\zasol,\thetasol)$ or show that no violated cuts can be found.  We use the multi-cut version of SDDP (i.e., one variable to approximate the objective value of each SDDP subproblem \citep{philpott2012dynamic}) and  the single-scenario iteration scheme (i.e., sample one scenario at a time from the scenario tree during the SDDP forward pass \citep{philpott2008convergence}). Lastly, our algorithm considers the quick pass version of SDDP  (i.e., we evaluate all the subproblems in a sample path before generating cuts in the backward pass \citep{morton1996enhanced}). We refer the reader to \cite{fullner2021sddpvariants} for further details on this and other SDDP variants.  

\begin{algorithm}[tb]
\small
    \linespread{1.2}\selectfont
	\caption{SDDP procedure} \label{alg:sddp}
	\begin{algorithmic}[1]
		\Procedure{\algSDDP}{$(\xrootsol,\zasol,\thetasol), n', \epsilon, K,  \algExact$}
		\Repeat
		\State Sample $K$ paths from the sub-tree rooted at $n'$ (i.e., $\tree(n')$) and store indices in $\K$. 
		\State \algCutAdd := \algFalse.
		\For{$k \in \K$}
		\For{$t \in \{2,...,\periods\}$}  {\hfill \color{gray} \%\% Forward pass sub-routines}
		\State Select node of scenario $k$ and stage $t$, $n:=n_{t,k}$, and its subproblem $\sub: =\sub(n)$.
		\State Solve subproblem $\sub$, let $\Qsubsol_{\sub} :=\Qsub_{\sub}(\hat{x}_{\parent(\sub)},\zcopyparentsol)$. Save primal and dual solutions.
		\EndFor
		\For{$t \in \{\periods,...,2\}$} {\hfill \color{gray} \%\% Backward pass sub-routines}
		\State Select node of scenario $k$ and stage $t$, $n:=n_{t,k}$, and its subproblem  $\sub: =\sub(n)$.
		\If{$|\thetasol_{\sub} - \Qsubsol_{\sub}| \geq \epsilon|\Qsubsol_{\sub}|$}
		\State Add inequality \eqref{eq:sddp_cut} to all subproblems with variables $\theta_{\sub}$, including $\Qsub_{\parent(\sub)}$.
		\State \algCutAdd := \algTrue.
		\If {$t = 2$} \algReturn\ $(\alpha_{\sub(n')}, \beta_{\sub(n')}, \gamma_{\sub(n')})$. {\hfill \color{gray} \%\% Add cut to master problem}
		\EndIf
		\EndIf
		\EndFor
		\EndFor
		\Until{\algConverge($K$, \algCutAdd, \algExact).} {\hfill \color{gray} \%\% Termination criteria}
		\EndProcedure
	\end{algorithmic} 
\end{algorithm}

We now provide a detailed explanation of the main components of the algorithm \algSDDP. We first sample (without replacement) $K$ scenarios (sample paths) from the sub-tree rooted at node $n'\in \stage_2$ to generate cuts for its associated variable $\theta_{n'}$ in the master problem. For each scenario, we then perform the forward and backward pass of the SDDP algorithm (lines 5-14). The forward pass iterates over all stages from $t=2$ to $t = \periods$ on a randomly sampled scenario (sample path), selects the associated node at each stage, and solves the corresponding SDDP subproblem (lines 6-8).  The backwards pass iterates over all stages from $t = \periods$ to $t = 2$ (lines 9-14) and checks the relative difference between the approximated value $\thetasol_{\sub}$ and $\Qsubsol_{\sub}$ (line 11). If this difference is larger than a given tolerance $\epsilon$, we generate a Benders cut as represented in \eqref{eq:sddp_cut}. In particular, the $(\alpha_{\sub'}, \beta_{\sub'}, \gamma_{\sub'})$ coefficients for any $\subn' \in \childsub(\subn)$ are given by:
\vspace{-0em}
\begin{align*}
    \alpha_{\subn'} & = {\pi^\eqref{eq:st_onlyx}}^\top F_{\subn'} + {\pi^\eqref{eq:sub_allvars}}^\top  A_{\subn'}, \qquad \beta_{\subn'} =  {\pi^\eqref{eq:sub_allvars}}^\top B_{\subn'}, \qquad
    \gamma_{\subn'}  = \Qsubsol_{\subn'} - \left( \alpha_{\subn'}^\top \hat{x}_{\subn} + \beta_{\subn'}^\top \zcopyphisol \right).
\end{align*}
where ${\pi^\eqref{eq:st_onlyx}}$ and $\pi^\eqref{eq:sub_allvars}$ are the dual solutions associated with inequalities \eqref{eq:st_onlyx} and \eqref{eq:sub_allvars} of subproblem $\sub'$, respectively. This cut is added to the parent subproblem (i.e., $\Qsub_{\sub(\parent(n))}$) but it is also valid for any subproblem in stage $t-1$ involving variable $\theta_{\sub(\parent(n))}$  \citep{philpott2012dynamic}. We note that while the forward pass randomly samples scenarios rooted at one particular child node of the root node, the backward pass can add cuts to subproblems in subtrees rooted at other children nodes of the root node. Also, we consider feasibility cuts during the forward pass if one of the subproblems is infeasible. 
Lastly, we terminate the current SDDP iterations for scenarios rooted at $n'$ if a cut is added to the first-stage (master) problem (i.e., $t=2$), otherwise, we continue iterating over the remaining scenarios until the termination criterion is met. 

Procedure \algConverge\ ensures that we terminate the SDDP algorithms when we can guarantee that no cut can be added to the master problem. To do so, we increase the number of sampled scenarios $K$ to be equal to the total number of scenarios (i.e., $K=|\stage_\periods|$)) when no cut was added using the original sample size. Once $K=|\stage_\periods|$, we continue iterating until no violated cut can be added to any subproblem (i.e., \algCutAdd=\algFalse). We note that this termination criterion can lead to several iterations over the full set of scenarios of the scenario tree, which is computationally expensive. However, this scenario increment is necessary to ensure an optimal solution to~\eqref{mod:aggregated}. 

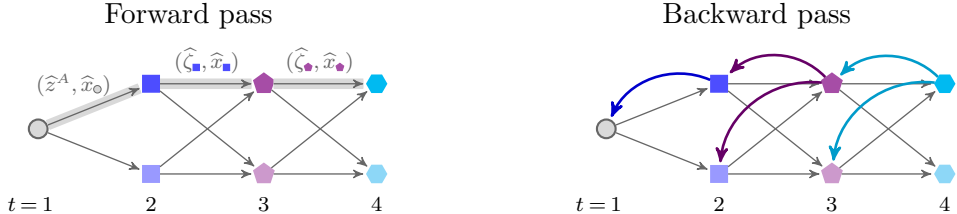
\begin{figure}[H]
    \centering
    \begin{tikzpicture}[->,>=stealth',shorten >=1pt,auto,node distance=1cm,
thick]       
\node[node light] (r) at (0,0) {};
\node[node-square, draw=blue!70!white, fill=blue!70!white] (n11)  at (1.5, 0.6)  {};
\node[node-square, draw=blue!40!white, fill=blue!40!white] (n12)  at (1.5, -0.6)  {};

\node[node-polygon5, draw=violet!70!white, fill=violet!70!white] (n21)  at (3,0.6)  {};
\node[node-polygon5, draw=violet!40!white, fill=violet!40!white] (n22)  at (3,-0.6)  {};

\node[node-polygon6, draw=cyan!80!white, fill=cyan!80!white] (n31)  at (4.5, 0.6)  {};
\node[node-polygon6, draw=cyan!40!white, fill=cyan!40!white] (n32)  at (4.5, -0.6)  {};

\node[] (name) at (2, 1.5) {Forward pass};

\node[text node] (name) at (-0.1, -1) {$t=1$};
\node[text node] (name) at (1.5, -1) {$2$};
\node[text node] (name) at (3, -1) {$3$};
\node[text node] (name) at (4.5, -1) {$4$};

\path[every node/.style={font=\sffamily\scriptsize}]
(r) 
edge[arc bold] node [above] {$(\zasol,\hat{x}_{\mcLight{0.4}})\qquad$} (n11)
edge[arc] node [left] {} (n12)
(n11)
edge[arc bold] node [above] {$(\hat{\zcopy}_{\squareDark{0.4}},\hat{x}_{\squareDark{0.4}})$} (n21)
edge[arc] node [left] {} (n22)
(n21)
edge[arc bold] node [above] {$(\hat{\zcopy}_{\pentaDark{0.4}},\hat{x}_{\pentaDark{0.4}})$} (n31)
edge[arc] node [left] {} (n32)
(n22)
edge[arc] node [left] {} (n31)
edge[arc] node [left] {} (n32)
(n12)
edge[arc] node [left] {} (n21)
edge[arc] node [left] {} (n22);

\end{tikzpicture}
    \hspace{5em}
    \begin{tikzpicture}[->,>=stealth',shorten >=1pt,auto,node distance=1cm,
thick]       
\node[node light] (r) at (0,0) {};
\node[node-square, draw=blue!70!white, fill=blue!70!white] (n11)  at (1.5, 0.6)  {};
\node[node-square, draw=blue!40!white, fill=blue!40!white] (n12)  at (1.5, -0.6)  {};

\node[node-polygon5, draw=violet!70!white, fill=violet!70!white] (n21)  at (3,0.6)  {};
\node[node-polygon5, draw=violet!40!white, fill=violet!40!white] (n22)  at (3,-0.6)  {};

\node[node-polygon6, draw=cyan!80!white, fill=cyan!80!white] (n31)  at (4.5, 0.6)  {};
\node[node-polygon6, draw=cyan!40!white, fill=cyan!40!white] (n32)  at (4.5, -0.6)  {};

\node[] (name) at (2, 1.5) {Backward pass};

\node[text node] (name) at (-0.1, -1) {$t=1$};
\node[text node] (name) at (1.5, -1) {$2$};
\node[text node] (name) at (3, -1) {$3$};
\node[text node] (name) at (4.5, -1) {$4$};

\path[every node/.style={font=\sffamily\scriptsize}]
(r) 
edge[arc] node [left] {} (n11)
edge[arc] node [left] {} (n12)
(n11)
edge[arc] node [above] {} (n21)
edge[arc] node [left] {} (n22)
edge[arc back, blue!80!black, bend right=45] node [left] {} (r)
(n21)
edge[arc] node [left] {} (n31)
edge[arc] node [left] {} (n32)
edge[arc back, violet!80!black, bend right=45] node [left] {} (n11)
edge[arc back, violet!80!black, bend right=45] node [left] {} (n12)
(n22)
edge[arc] node [left] {} (n31)
edge[arc] node [left] {} (n32)
(n12)
edge[arc] node [left] {} (n21)
edge[arc] node [left] {} (n22)
(n31)
edge[arc back, cyan!80!black, bend right=45] node [left] {} (n21)
edge[arc back, cyan!80!black, bend right=45] node [left] {} (n22);

\end{tikzpicture}
    \caption{An illustrative example of the SDDP integrated B\&C algorithm. It depicts decomposition structure and the two main SDDP subroutines (i.e., forward and backward pass).}
    \label{fig:sddp_b&c}
\end{figure}

\begin{example}\label{exa:sddp_b&c}
Figure \ref{fig:sddp_b&c} illustrates the main components of our SDDP integrated B\&C algorithm for the \ma\ transformation. Note that in this case an SDDP subproblem is defined for each MC state in each stage $t > 1$. The two illustrations show the decomposition structure over the policy graph: the root node \mcLight{1}\ corresponds to the master problem while other nodes correspond to the SDDP subproblems. The drawing on the left depicts the forward pass of Algorithm \ref{alg:sddp} for scenario $\mcLight{1} \rightarrow \mcDark{1}  \rightarrow \mcDark{1} \rightarrow \mcDark{1}$. Note that the solution of the SDDP master problem, $(\zasol, \hat{x}_{\mcLight{0.6}})$, is used as part of the input parameters to \squareDark{1}\ and, similarly, $(\hat{\zcopy}_{\squareDark{0.6}}, \hat{x}_{\squareDark{0.6}})$ is used as part of the input to \pentaDark{1}. Recall that the $\zcopy$ variables are just the copies of $\za$, so $\zasol=\hat{\zcopy}_{\squareDark{0.6}} = \hat{\zcopy}_{\pentaDark{0.6}}$. The drawing on the right illustrates the backward pass where the dual optimal solution of an SDDP subproblem is used to generate cuts to improve the outer-approximation of the cost-to-go function defined at its parent subproblems. For example, the optimal dual solution of \pentaDark{1}\ is used to generate cuts for subproblems \squareDark{1}\ and \squareLight{1}.
\hfill $\square$
\end{example}

\subsubsection{Lower Bound Computation.}\label{sec:sddp_lb}

Algorithm \ref{alg:bc-sddp} guarantees an optimal solution to~\eqref{mod:aggregated}, but the procedure can be computationally expensive because the number of SDDP subproblems depends on the underlying MC structure and the chosen transformation. In addition, we may need to call the SDDP sub-routine multiple times to correctly evaluate each candidate integer relaxation solution encountered during the B\&C procedure. With this in mind, we propose an alternative to the aforementioned procedure to considerably reduce its computation time and return valuable information about the problem in the form of a lower bound. This lower bound can be used, for example, to provide an overestimate of the optimality gap associated with a known feasible solution. We utilize this lower bound to evaluate the quality of solutions obtained by our proposed 2SLDR approximations (see Section~\ref{sec:ldr} for details).

The overall idea for this lower-bounding technique is to apply a limited number of iterations of the SDDP sub-routine. First, we  increase the cut violation tolerance parameter $\epsilon$ only to add cuts that have a significant relative difference to its true value (e.g., using $\epsilon=0.1$ instead of a conservative value $\epsilon=0.0001$). Also, we limit the number of rounds in the \algConverge\ procedure in Algorithm~\ref{alg:sddp}. Specifically, when Boolean variable $\algExact= \algFalse$, we limit the number of rounds to only three and avoid increasing the number of sampled scenarios $K$. This is in contrast with the $\algExact= \algTrue$ version (see the previous subsection), where we iterate until no more cuts are added to any subproblem and the number of sampled scenarios increases to cover all possibilities.
These modifications to the termination criteria considerably reduce the computational time of the SDDP procedure. However, we only obtain lower-approximations to the cost-to-go functions as we might be missing cuts in the SDDP master and subproblems, resulting in a lower bound for \eqref{mod:aggregated}. 

\subsection{Two-stage Linear Decision Rule Approximation} \label{sec:ldr}
In this section, we describe an approximation scheme to create feasible solutions for \eqref{mod:aggregated} by imposing additional structures to the decision policies in the form of 2SLDR  \citep{bodur2018two}. Imposing additional structures in such a way transforms the multi-stage stochastic program into a two-stage stochastic program, which requires significantly less computational effort than the B\&C framework integrated with the SDDP algorithm that we presented in the previous section. Moreover, the 2SLDR can be customized to leverage the structure of the underlying stochastic process, and we provide new schemes for the MC case.



Let $\xi_t$ be the vector of random variables at stage $t\in \horizon$, and $\xi_{t,n}$ be the realization of these random variables at node $n\in \staget$,  i.e., 
$ \xi_{t,n} = \{A_n,B_n,C_n,D_n,E_n,F_n,J_n,b_n,c_n,d_n,f_n,h_n\}. $ Note that $\xi_\rootnode$, the data associated with the root node of the scenario tree, is assumed to be a deterministic vector. In addition,  $\xi^t=\{\xi_1,...,\xi_t\}$ and $\xi^t_n=  \{\xi_{1,\rootnode},...,\xi_{t-1, \parent(n)}, \xi_{t,n}\}$ represent the trajectory of the random variables up to stage $t$ and their realizations for a particular node $n$ in the scenario tree, respectively. We consider $\xi_t$ and $\xi^t$ as vectors of random variables with dimensions $\dimxit$ and $ \dimxitt =\sum_{i = 1}^t \dimxit$ for all $t\in \horizon$, respectively, and the same applies to the realizations of these random variables.   

The general idea behind 2SLDRs is to replace the original \Merve{continuous} state variables $\xn$ with a new set of variables for each node $n \in \nodes$ by means of a linear relation of node-specific data. By doing so, we impose additional (linear relation) structure over the state variables, and in return  transform the multi-stage stochastic program into a two-stage stochastic program \citep{bodur2018two}. Specifically, 2SLDRs employ basis functions $\mapldr(\cdot)$ to map $\xi^t_n$ into values that are then used to create a linear representation for the state variables $\xn$, i.e., $\xn = \varldr_t^\top \mapldr(\xi^t_n) \in \R^k$, where $\varldr_t$ is a vector of the so-called LDR variables at stage $t$. Although the basis functions $\mapldr(\cdot)$ can be any form of functions that map the node-specific data $\xi^t_n$ to certain values, for ease of exposition, we consider simple basis functions that return the actual random variable realizations, that is, $\mapldr(\xi^t_n) = \xi^t_n$. In this case, the LDR is given by:
\vspace{-0em}
\begin{equation}\label{eq:ldr}
    \xn = \varldr_t^\top \xi^t_n, \ \forall n \in \staget, t \in \horizon.
\end{equation}
We observe that these simple basis functions return high-quality decision policies for our case study (see Section \ref{sec:experiments}), although more complex basis functions can potentially lead to better decision policies \citep{chen2008linear,bampou2011scenario}.

\subsubsection{Aggregation Framework for 2SLDR.}
Similar to the aggregation framework introduced in Section~\ref{sec:framework}, we create different aggregation  schemes for the LDR variables $\varldr$, and as a result, obtain approximate 2SLDR-based decision policies. Table \ref{tab:ldr} summaries a few LDR variable structures that we consider. The first two LDR variants (\thldr\ and \tldr) are common in the literature \citep{shapiro2005complexity,kuhn2011primal}, while \mldr\ is a new LDR variant based on the MC structure of the stochastic process.  
\begin{table}[htb]
\setlength\extrarowheight{4pt}
    \centering
    \begin{tabular}{c|l|l|l}
    \toprule
    Name & LDR variables for $t\in \horizon$ & Standard basis function & Resulting LDR  \\
    \midrule
    \thldr & $\varldr_t \in \R^{k\times \dimxitt}$ & $\mapldr(\xi^t_n) = \xi^t_n$ & $\xn =  \varldr_t^\top \xi^t_n$ \\
    \tldr & $\varldr_{t} \in \R^{k\times \dimxit}$ & $\mapldr(\xi^t_n) = \xi_{t,n}$ & $\xn = \varldr_t^\top \xi_{t,n}$ \\
    \mldr & $\varldr_{t,m}\in \R^{k\times\dimxit},  \quad \forall\; m \in \markov_t$ & $\mapldr(\xi^t_n) = \xi_{t,n}$ &  $\xn =  \varldr_{t, m(n)}^\top \xi_{t,n}$ \\
    \bottomrule
    \end{tabular}
    \caption{Examples of LDR variable structures and basis functions for each $n\in \staget$ and $t\in \horizon$.}
    \label{tab:ldr}
\end{table}

The first LDR variant, \thldr, which we refer to as a stage-history LDR, constructs the linear relation based on the trajectory of realizations of the random variables up to the current stage. The \thldr\ variant can lead to high-quality decision policies (thanks to the amount of information incorporated) but may result in a large number of LDR variables, especially for problems with a large number of stages. An alternative is \tldr\ which constructs the linear relation only based on the realization of the random vector at the current stage, which leads to a smaller number of LDR variables at the expense of lower solution quality~\citep{bodur2018two}.

In addition to these LDR variants, we  propose \mldr\ which takes advantage of the MC structure of the underlying stochastic process. The idea is to have one copy of LDR variables $\varldr_{t,m}$ for each MC state $t$ and stage $m$. It is clear that \mldr\ will lead to better policies compared to \tldr\ because both LDR variants use the same information (realizations of random vectors), but the latter has a single copy of the LDR variables for all nodes in the same stage. We note that, however, \mldr\ and \thldr\ are incomparable in general, as shown in our empirical results.

\subsubsection{2SLDR Models.}

We now present the resulting model after applying the LDR transformation \eqref{eq:ldr} to \eqref{mod:aggregated}. For ease of exposition, we only present the 2SLDR variant \thldr\ and the other two alternatives are similar. The overall idea is to replace the state variables $\xn$ with the LDR variables using linear function \eqref{eq:ldr}. The resulting model for the root node is as follows:
\vspace{-0em}
\begin{subequations}
\begin{align}
    \Qldr_{\rootnode} = \min\; & 
    \sum_{n\in \nodes}\prob_n \left( c_n^\top \zphi + d_n^\top(\varldr_{t(n)}^\top\xi^{t(n)}_n) \right) + h_\rootnode^\top \yroot + \sum_{n \in \nodesnoroot }\prob_{n} \Qldr_n(\varldr,\za) \tag{$P^L$}  \label{mod:ldr} \\ 
    \text{s.t.}\; &  H_{\rootnode} \zroot \geq g_{\rootnode}, \;
    J_{\rootnode} \varldr_{1}^\top\xi^{1}_\rootnode \geq f_{\rootnode}, \;
    C_{\rootnode} \varldr_{1}^\top\xi^{1}_\rootnode + D_{\rootnode} \zroot + E_{\rootnode} \yroot  \geq  b_{\rootnode}, \label{eq:ldr-constr0} \\
    & H_n \zphi \geq  G_n \zphiparent + g_n, & \forall n \in \nodesnoroot,  \label{eq:ldr-constr1} \\
    & J_n \varldr_{t(n)}^\top\xi^{t(n)}_n \geq F_n  \varldr_{t(a(n)) }^\top\xi^{t(a(n))}_{\parent(n)}, &  \forall n \in \nodesnoroot, \label{eq:ldr-constr2} \\
    & (\yroot, \varldr, \za)  \in  \R^r \times \R^{k\cdot \sum_{t \in \horizon}\dimxitt} \times \Z^{\ell \cdot \sum_{t \in \horizon} q_t}. \nonumber
\end{align}
\end{subequations}

Model \eqref{mod:ldr} contains all the root node variables and constraints of \eqref{mod:aggregated} and includes all the LDR variables and its associated constraints \eqref{eq:ldr-constr2} from all nodes in the scenario tree. We note that there is no need to define LDR variables for the root node since the $\xroot$ variables are already first-stage variables, but we include them in the model for ease of exposition. Lastly, the cost-to-go function for node $n \in \nodesnoroot$ only depends on the first-stage variables $\varldr$ and $\za$ and is given by:
\vspace{-0em}
\begin{equation*}
    \Qldr_n(\varldr, \za) = \min_{\yn \in \R^r} \left\{ h_n^\top \yn \mid  C_n \varldr_{t(n)}^\top\xi^{t(n)}_n + D_n \zphi + E_n \yn  \geq A_n \varldr_{t(n) -1 }^\top\xi^{t(n)-1}_{\parent(n) } + B_n  \zphiparent + b_n\right\}. 
\end{equation*}
Note that \eqref{mod:ldr} is indeed a two-stage stochastic program. Applying \eqref{eq:ldr} 
we remove the dependency between parent and child nodes, thus, all cost-to-go functions become  functions of only the first-stage decisions. This two-stage stochastic program is easier to solve than its multi-stage 
counterpart since it eliminates the nested dependency between the cost-to-go functions in different stages.  

Formulation \eqref{mod:ldr} considers one cost-to-go function for each node in the scenario tree. However, the number of cost-to-go functions can be considerably smaller depending on the LDR variant and transformation $\Phi_t$. For example, \thldr\ requires one cost-to-go function for each node in the scenario tree because it considers the entire trajectory of realizations of random variables. Therefore, the cost-to-go functions for \thldr\ are unique for each node even if they share the same stage and MC state. On the contrary, \tldr\ and \mldr\ only employ local information, so the cost-to-go functions for two nodes $n, n'\in \staget$ are the same if $\sub(n) = \sub(n')$, that is, if they share the same stage $t$, MC state $m(n)=m(n')$, and aggregated variables $\phi_t(n)= \phi_t(n')$. Therefore, we have the same number of distinct cost-to-go functions as SDDP subproblems, which can be significantly smaller than the size of the scenario tree. In conclusion, \tldr\ and \mldr\ have a significant computational advantage over \thldr\ due to the fewer LDRs variables and cost-to-go functions employed (along with all the variables and constraints necessary to define them). 

We solve the resulting two-stage stochastic programming model using Benders decomposition, a standard approach in the stochastic programming literature (see, e.g.,~\cite{zverovich2012computational}). Our specific implementation details can be found in Appendix \ref{app:2sldr_benders}.

\section{Case Study: A Hurricane Disaster Relief Planning Problem} \label{sec:application}
In this section, we present a case study to exemplify the aggregated framework and evaluate the proposed solution methodologies. We consider a class of multi-period hurricane disaster relief logistics planning (HDR) problem, where a number of contingency modality options can be activated during the planning horizon to increase the capacities at the distribution centers (DCs). 

The HDR problems are well-studied in the operations research literature because of their practical impacts on reducing economic loss and human suffering~\citep{graumann2006hurricane,seraphin2019natural}. The operational challenges in HDR problems arise from hurricanes' stochastic nature \citep{sabbaghtorkan2020prepositioning}. To address these challenges, most of the HDR literature focuses on two-stage stochastic programming models, where the first-stage model makes strategic facility location and resource pre-positioning decisions in terms of shelter preparation and resource allocation, while the second stage  deals with disaster relief logistics decisions once the damage and demand information are revealed after the hurricane's landfall (see, e.g., \cite{duran2011pre,lodree2012pre,davis2013inventory,alem2016stochastic}, and a survey paper by \cite{sabbaghtorkan2020prepositioning}). 

Recent works have also considered the multi-stage HDR variants where sequential logistics decisions are made in a rolling horizon (RH) fashion~\citep{pacheco2016forecast,siddig2022multi,yang2022optimizing}. For example, \cite{pacheco2016forecast} consider a relief supply preposition and reposition problem starting from the time when the hurricane is first detected, with updated forecast information on the hurricane's attributes every six hours. The authors proposed a forecast-driven dynamic model for the problem, and a solution approach combining scenario analysis and RH. Similarly, \cite{siddig2022multi} consider a multi-period preposition problem where the demand is realized at landfall (i.e., last period). The authors proposed an MSLP formulation and showed the value of multi-stage stochastic programming solutions compared to solutions given by an RH approach and a two-stage approximation. 

Despite their wide range of problem settings, the HDR problems have several assumptions in common. For example, the evolution of the hurricane and intensity are commonly modeled using an MC  \citep{taskin2010inventory,pacheco2016forecast} and the estimated demand at each shelter can be modeled using a deterministic mapping from the current hurricane state (i.e., hurricane location and intensity) \citep{siddig2022multi}. Similarly, the HDR problem that we study follows these assumptions as well. One unique feature of our HDR problem is that we consider a set of contingency modality options that can be activated during the planning horizon to increase the inventory capacities of the DCs. We introduce more details about the problem next.

\subsection{Problem Description and Model Formulation}

Our multi-stage HDR problem is concerned with the logistics decisions of producing and distributing relief commodities to shelters prior to the landfall of a hurricane. 
The relief commodities are used to fulfill the demand for civilians that are evacuated to shelters. As in the related problems, we consider that demand can be estimated according to the location and intensity of the hurricane, updated every six hours, which is the frequency of the hurricane forecast update by the National Hurricane Center~\citep{regnier2019hurricane}. We assume that shelters, DCs, and possible positions of the hurricane are located in a two-dimensional grid with a discrete set of $x$- and $y$-coordinates. The 
goal is to produce and transport relief commodities from the DCs to the shelters at a minimum (expected) cost, which consists of the penalty cost for unsatisfied demand and operational costs incurred in all time periods. 

In addition, we consider a set of contingency modalities that can be activated to increase the DCs' capacities. These modalities make the disaster relief logistics system more resilient to extreme hurricane events that may lead to unexpectedly high demand at the price of higher logistics costs. We make a practically relevant assumption that these modalities can be activated at most once during the planning horizon, and once activated, the modality will stay the same until the end of the planning horizon. An active modality represents an incremental capacity increase at every stage on the planning horizon. These considerations are consistent with the practical constraints associated with these large-scale critical logistics operations \Song{(e.g., the phased activation of national guards for disaster relief efforts)}, which have this ``all or nothing'' feature and have to be implemented in small increments. Lastly, we consider a one-leg delay for these modalities (i.e., the capacity starts to increase from the stage immediately after the stage when a modality activation decision is made), which partially captures the logistical challenges of modality changes during disaster relief.

We assume that the hurricane's evolution can be characterized by an MC $\markov$. Each MC state corresponds to two attributes of the hurricane: location and intensity, that is, $m=(m^x, m^y, m^i)\in \markov$ where $m^x$ and $m^y$ are the $x$- and $y$-coordinates of the hurricane's location, respectively, and $m^i$ represents the hurricane's intensity level. We consider independent probability transition matrices for each hurricane attribute. For the hurricane intensity, We use the transition probability matrix described in~\cite{pacheco2016forecast} in our test instances. For the hurricane movement, we assume that the hurricane originates from the bottom row of the grid (i.e., $m^y = 0$ for the initial MC state) and that $m^y$ increases by one in each period (i.e., the hurricane advances upwards by one step in each period) and will reach land in exactly $T$ periods (see Appendix \ref{app:hdr} for additional details). \Song{We note that our approach can be easily extended to the case when the number of stages until landfall is random, as long as the underlying stochastic process is modeled by an MC. }

Given the stochastic behavior of the hurricane, we consider an MSILP model with a stage defined for each time period when new information is realized (i.e., the hurricane moves to a new MC state). We construct a scenario tree $\tree$ by considering  an initial MC state and all possible MC transitions for up to $T$ stages (i.e., until  landfall). We consider a set of contingency modalities $\modality$, a set of shelters $\shelters$, and a set of DCs $\centers$, where each DC $j\in \centers$ has an initial capacity $C_j$ and inventory $I_j$. At each node $n\in \nodes$ of the scenario tree $\tree$, one needs to decide the amount of relief commodities $v_{j}\geq 0$ to produce for each DC $j \in \centers$, the amount of relief commodities $y_{ij}\geq 0$ to transport from each DC $j\in \centers$ to each shelter $i \in \shelters$, and the amount of unsatisfied demand $w_i \geq 0$ at each shelter $i\in \shelters$. We introduce two continuous state variables $x_n = (\xinv_{n}, \xcap_{n})$ to represent the inventory and capacity for each DC $j\in \centers$ at node $n\in \nodes$, respectively. Lastly, we introduce a binary state variable $z_{n\ell} \in \{0,1\}$ to represent whether or not a contingency modality $\ell \in \modality$ is active at node $n \in \nodes$.

Each node $n\in \nodes$ is associated with realizations of random variables according to the corresponding MC state of the hurricane. Specifically, for a given node $n \in \staget$, $t \in \horizon$, the demand of a shelter $i\in \shelters$, $d_{ni}$, the production cost at DC $j\in \centers$, $q_{nj}$, and the transportation cost from DC $j$ to shelter $i$, $f_{nij}$, depend solely on the MC state $m_t(n)$ associated to node $n$. Thus, two nodes $n,n'\in \staget$ with the same MC state at stage $t$, $m_t(n)=m_t(n')$, share the same realizations of random variables. The remaining model parameters are assumed to be deterministic, which include the unit inventory cost $g_j$ at $j\in \centers$, the contingency modality cost $c_\ell$ for $\ell\in \modality$, the penalty cost for each unit of unsatisfied demand $b_i$ for $i\in \shelters$, and the capacity increase $K_{j\ell}$ at DC $j$ in each stage when modality $\ell$ is active. Further details on the problem instance  generation can be found in Appendix \ref{app:hdr}. The MSILP model for the HDR problem is given by: 
\vspace{-0em}
\begin{subequations}
\begin{align}
     \min \quad &  \sum_{j \in \centers}\left( g_{j}\xinv_{\rootnode j} + q_{\rootnode j} v_{j} +\sum_{i \in \shelters}  f_{\rootnode ij} y_{ij} \right) + \sum_{i \in \shelters} b_{i}w_{i} + \sum_{\ell \in \modality} z_{\rootnode \ell}c_{\ell} + \sum_{n \in \childs(\rootnode)} \lefteqn{\pprob_{\rootnode n} Q_n(\xroot, \zroot) } \tag{$\textit{HDR}$} \label{dr_root} \\
    \text{s.t.} \quad & \sum_{j \in \centers} y_{ij} + w_{i} \geq d_{\rootnode i}, & \forall i \in \shelters,  \label{dr_root:demand} \\
    & \xinv_{\rootnode j} = I_{j} - \sum_{i \in \shelters} y_{ij} + v_{j}, & \forall j \in \centers, \label{dr_root:inventory} \\
    & v_{j} \leq \xcap_{\rootnode j},  & \forall j \in \centers, \label{dr_root:production}\\
    & \xcap_{\rootnode j} = C_j,  & \forall j \in \centers, \label{dr_root:capacity} \\
    & \sum_{\ell \in \modality} z_{\rootnode \ell} \leq 1, \label{dr_root:modality} \\
    & \xinv_{\rootnode j}, \xcap_{\rootnode j},   v_{j}, w_{i}, y_{ij} \geq 0, \; \zroot{}_\ell \in \{0,1\}, & \forall j \in \centers, \; i \in \shelters, \; \ell \in \modality.  \nonumber 
\end{align}
\end{subequations}
Constraints \eqref{dr_root:demand} and \eqref{dr_root:inventory} represent the demand and inventory constraints, respectively. Constraints \eqref{dr_root:production} enforce the production capacity and constraints \eqref{dr_root:capacity} represent the initial capacity of each DC. Constraint \eqref{dr_root:modality} enforces that at most one modality can be activated at the root node of the scenario tree. Lastly, the cost-to-go function for a node $n\in \nodesnoroot$ is given by:
\vspace{-0em}
\begin{subequations}
\begin{align}
    Q_n(\xnparent, \znparent) &  =  \min \sum_{j \in \centers}\left( g_{j}\xinv_{n j} + q_{n j} v_{j} + \sum_{i \in \shelters}f_{nij} y_{ij} \right) +  \sum_{i \in \shelters} b_{i}w_{i} + \lefteqn{\sum_{\ell \in \modality} z_{n \ell}c_{\ell}  
       +\!\! \sum_{n' \in \childs(n)}\!\!\pprob_{n n'} Q_{n'}(\xn, \zn) } \quad 
       \nonumber \\
   \text{s.t.} \quad  & \sum_{j \in \centers} y_{ij} + w_{i} \geq d_{ni}, & \forall i \in \shelters, \label{dr_node:demand} \\
    & \xinv_{nj} = \xinv_{a(n)j} - \sum_{i \in \shelters} y_{ij} + v_{j}, & \forall j \in \centers, \label{dr_node:inventory} \\
    & v_{j} \leq \xcap_{nj},  & \forall j \in \centers,\label{dr_node:production} \\
    & \xcap_{nj} = \xcap_{a(n)j} + \sum_{\ell \in \modality} K_{j \ell}z_{a(n) \ell},  & \forall j \in \centers, \label{dr_node:capacity}\\
    & \sum_{\ell \in \modality} z_{n\ell} \leq 1, &  \label{dr_node:modality}\\
    & z_{a(n)\ell} \leq  z_{n\ell}, & \forall \ell \in \modality,  \label{dr_node:modality2}\\
    & \xinv_{n j}, \xcap_{n j},   v_{j}, w_{i}, y_{ij} \geq 0, \; \zn{}_\ell \in \{0,1\}, & \forall j \in \centers, i \in \shelters, \ell \in \modality. \nonumber 
\end{align}
\end{subequations}
The meanings of constraints \eqref{dr_node:demand}-\eqref{dr_node:production} and \eqref{dr_node:modality} are identical to their counterparts in~\eqref{dr_root}. Constraints \eqref{dr_node:capacity} represent the capacity increase at a DC depending on whether or not a contingency modality is active in the previous stage. Constraints \eqref{dr_node:modality2} enforce that an active contingency modality in the previous stage will remain active in the current stage. 


\subsection{Aggregation Framework}
\Margarita{
We consider the four MC-based transformations introduced in Table \ref{tab:transformations} of \Merve{Section \ref{sec:framework}}, i.e., \hn, \ma, \mm, and \pa, for the aggregated model. The mathematical formulation of the aggregated model~\eqref{dr_aggre} can be found in Appendix \ref{app:hdr_aggregated}, which is a straightforward application of our proposed framework. 
  
\Merve{In the context of our HDR problem, aggregation is done for the contingency modality decisions.}
\hn\ considers the same \Merve{decisions} for all nodes in a given stage;  \ma\ aggregates variables based on the 
\Merve{hurricane} state information in a given stage, which corresponds to the location and intensity of the hurricane; and 
\mm\ incorporates the \Merve{hurricane}
states at the current node and its parent node. The 
transformation \pa\ can be customized 
depending on how one defines the partial MC state information. In our experiments, we consider the full \Merve{hurricane} state information 
(i.e., \Merve{both} location and intensity) at the current node and only the hurricane intensity information at its parent node, that is, the difference between \mm\ and \pa\ for the HDR problem is that the former incorporates the hurricane location information at the parent node into the aggregation, while the latter does not. 
}

\subsection{Two-stage LDR approximation}

As explained in Section \ref{sec:ldr}, we can apply a further restriction to the aggregated MSILP model \eqref{dr_aggre} to obtain a two-stage model by applying LDR to the continuous state variables (i.e., inventory variables $\xinv_n$). We present three LDRs variants similar to the ones introduced in Table \ref{tab:ldr}.

First, we introduce notation $\xi_{t,n}=(d_n,q_n,f_n)$ to represent the vector of realizations of random variables for each node $n \in \staget$ in the scenario tree, including the demand $d_n$, production cost $q_n$, and transportation cost $f_n$. Since the demand is one of the most important parameters for the HDR problem, in our experiment we choose to use basis functions (for constructing LDRs) that are only based on either the demand realization at the current node $\mapldr(\xi_{n}^t) = d_n$ or the entire history of demand realizations up to the current node $\mapldr(\xi_{n}^t) = (d_\rootnode,...,d_{a(n)}, d_n)$. Then, the LDR variants for the inventory variables for a given node $n\in \staget$ and DC $j \in \centers$ are defined as:
\vspace{-0em}
\[
    \thldr: \; \xinv_{n j} = \!\!\sum_{n' \in \patht(n)}\sum_{i \in \shelters}\varldr_{t(n')tji}d_{n'i} \quad\; 
    \tldr: \;  \xinv_{n j} = \sum_{i \in \shelters}\varldr_{t(n)ji}d_{ni} \quad\;
    \mldr: \; \xinv_{n j} = \sum_{i \in \shelters}\varldr_{t m_t(n) ji}d_{ni}
\]
where $\varldr_{t't} \in \R^{|\centers|\times |\shelters|}$ for $t' \in [t]$, $\varldr_{t} \in \R^{|\centers|\times |\shelters|}$, and $\varldr_{tm}\in \R^{|\centers|\times|\shelters|}$  for $m \in \markov_t$ are the LDR variables associated to each variant, respectively. 
\thldr\ considers the full demand history, while \tldr\ and \mldr\ only consider the demand realization at the current node. 
\mldr\ has one set of variables for each \Merve{hurricane state,}
while \tldr\ and \thldr\  employ one set of variables per stage. We discuss the implications of these differences in Section \ref{sec:ldr}, also in our numerical results in Section \ref{sec:experiments-approx}.

These LDR variants lead to three different two-stage approximations of \eqref{dr_aggre}, which are solved via Benders decomposition as discussed in Section \ref{sec:ldr} and Appendix \ref{app:2sldr_benders}. Note that while the MSILP problem has the relatively complete recourse property, these two-stage approximations do not necessarily have this property. For instance, LDR variables $\varldr$ can take negative values, which may lead to an infeasible solution with negative inventory values. \Margarita{We could force $\varldr\geq 0$ to void this issue, but this further limits the approximation and can lead to a worse solution, as observed during preliminary experimentation.} 
Thus, we incorporate Benders feasibility cuts \Merve{in this regard.}

\section{Numerical Experiment Results} \label{sec:experiments}

We now present the numerical results for the proposed aggregation framework 
and solution approaches.
We test the models and approaches over various instances of the described HDR problem.
In what follows, we first provide details on the test instances and the experimental setup. We then evaluate the computational performances of the proposed methodologies, focusing on the approximation methods (i.e., 2SLDR solutions with SDDP-based lower bounds) which perform best in practice. Lastly, we analyze the obtained decision policies and provide managerial insights. 

\subsection{Experimental Setup}

We generate instances for our MSILP HDR problem for two grid sizes, $4\times5$ and $5\times6$, which correspond to small-size and large-size instances, respectively (see Appendix \ref{app:model_size} for more details on the sizes of these instances). We experiment with different levels of initial capacities for the DCs (i.e., $C_j$ for $j \in \centers$) to analyze the impact of decision policies derived from the aggregation framework to systems with different capacities to address demand fluctuations. 
Specifically, we consider instances where the initial capacities across all DCs are 20\%, 25\%, or 30\% of the maximum demand across all shelters (the maximum demand is a predefined parameter for instance generation, see Appendix \ref{app:hdr} for details). We also consider two types of contingency modality options: conservative (Type-1) and aggressive (Type-2). For Type-1 modality, DC capacity expansion options are 10\%, 20\%, 30\%, or 40\% per stage, while for Type-2 modality, the options are 15\%, 30\%, 45\%, or 60\% per stage. We randomly generate a number of instances for each instance configuration and eliminate the ones where the optimal decision policies under transformations \hn\ and \fh\ are identical to make a meaningful comparison between different transformations. Overall, we obtain ten random instances for each of the six different instance configurations on two grid sizes.


The code for the decomposition algorithms was implemented in \texttt{C++} using solver IBM ILOG CPLEX 20.1 with callback functions. All experiments were run using a single core/thread and a six-hour time limit. The experiments were run in the University of Toronto SciNet server Niagara, using cores with 16GB RAM each\footnote{See \url{https://docs.scinet.utoronto.ca/index.php/Niagara_Quickstart} for further server specifications.}. We will make our code and instances available upon publication.

\subsection{Computational Performance of the Proposed Methodologies} \label{sec:experiments-approx}

We now present the empirical performance of the proposed methodologies to solve the aggregated MISLP models. We focus on the best performing procedures which could be more valuable for practitioners, that is, the  2SLDR variants described in Section \ref{sec:application} (i.e., \thldr, \tldr, and \mldr) and the SDDP-based methods to obtain a lower bound (\sddplb) and an upper bound (\sddpub). We briefly discuss the results of the exact methods at the end of this section, that is, the B\&C SDDP approach (\sddp) and the extensive model (\extensive), and refer to Appendix \ref{app:experiments-exact} for further details. For brevity, we only show results for the \pa\ transformation since it achieves the best trade-off between policy quality and computational effort. See Appendix \ref{app:approx-additional-resutls} for the results of the other transformations.

We first present results over small-size instances in Table \ref{tab:2sldr-pm-smallsize}. Columns ``Average Time (sec)'' report each approach's 
time to obtain their respective optimal solution (all techniques obtained optimal solutions within the time limit), averaged over instances with the same configuration. Columns ``Relative Difference (\%)'' compare the solution found by each approach and the (true) optimal solution obtained from \extensive, that is, $(|Obj_{\extensive} - Obj_{i}|)/Obj_{\extensive}$ for all $i \in \{\thldr,\tldr,\mldr,\sddpub,\sddplb\}$. We observe that \tldr\ and \mldr\ take the least 
time on average. Moreover, they yield the smallest relative difference (
$<$ 0.3\% on average in all configurations), indicating the high quality of their corresponding solutions. In particular, the \Merve{newly} proposed LDR variant, \mldr, provides the best approximation and takes the least time in most configurations. This shows the value of developing decision policies that leverage the structure of the underlying stochastic process. 

\begin{table}[tb]
\def\arraystretch{\stretchTableResults}
  \centering
  \caption{Solution time and  quality of 2SLDR and SDDP bounds. Results for $\pa$ over small-size instances.}
  \scalebox{\scaleTableResults}{
    \begin{tabular}{cc|rrr|r|r|rrrr|r}
    \toprule
          &       & \multicolumn{5}{c|}{Average Time (sec)} & \multicolumn{5}{c}{Relative Difference (\%)} \\
    \midrule
    Modality & \multicolumn{1}{l|}{Cap.} & \multicolumn{1}{c}{\thldr} & \multicolumn{1}{c}{\tldr} & \multicolumn{1}{c|}{\mldr} & \multicolumn{1}{c|}{\sddplb} & \multicolumn{1}{c|}{\extensive} & \multicolumn{1}{c}{\thldr} & \multicolumn{1}{c}{\tldr} & \multicolumn{1}{c}{\mldr} & \multicolumn{1}{c|}{\sddpub} & \multicolumn{1}{c}{\sddplb} \\
    \midrule
    \multirow{3}[2]{*}{Type-1} & \multicolumn{1}{l|}{20\%} & 813.0 & 164.4 & \textbf{157.9} & 549.1 & 961.3 & 0.12  & 0.26  & 0.11  & \textbf{0.01} & 0.36 \\
          & \multicolumn{1}{l|}{25\%} & 380.1 & 84.3  & \textbf{77.7} & 347.0 & 487.2 & 0.08  & 0.11  & \textbf{0.03} & 0.26  & 0.71 \\
          & \multicolumn{1}{l|}{30\%} & 408.2 & \textbf{79.0} & 85.9  & 400.7 & 222.0 & 0.12  & 0.25  & \textbf{0.00} & 0.55  & 1.24 \\
    \midrule
    \multirow{3}[2]{*}{Type-2} & \multicolumn{1}{l|}{20\%} & 1406.0 & \textbf{247.2} & 265.2 & 940.5 & 2047.3 & 0.13  & 0.26  & 0.13  & \textbf{0.00} & 0.28 \\
          & \multicolumn{1}{l|}{25\%} & 670.5 & 108.2 & \textbf{98.1} & 534.4 & 651.0 & 0.09  & 0.13  & \textbf{0.04} & 0.23  & 0.49 \\
          & \multicolumn{1}{l|}{30\%} & 532.6 & 90.0  & \textbf{89.5} & 434.8 & 270.6 & 0.12  & 0.26  & \textbf{0.01} & 0.98  & 1.29 \\
    \midrule
    \multicolumn{2}{c|}{Average} & 701.7 & \textbf{128.9} & 129.1 & 534.4 & 773.2 & 0.11  & 0.21  & \textbf{0.05} & 0.34  & 0.73 \\
    \bottomrule
    \end{tabular}%
    }
  \label{tab:2sldr-pm-smallsize}%
\end{table}%

Table \ref{tab:2sldr-pm-smallsize} also shows the performance of the SDDP lower (\sddplb) and upper (\sddpub) bound techniques. Both provide high-quality bounds with small relative differences (i.e., less than 1.5\%). Also, \sddplb\ is computationally more efficient compared to \extensive\ in most cases. Nonetheless, 
\sddpub\ (evaluating the incumbent solution from \sddplb) can \Merve{require} significant computational effort due to the additional SDDP calls:
around 1-2 hours for small-size instances and up to $12$ hours for large ones.

\begin{table}[tbp]
\def\arraystretch{\stretchTableResults}
  \centering
  \caption{Solution time and quality of 2SLDR and SDDP bounds. Results for $\pa$ over large-size instances.}
  \scalebox{\scaleTableResults}{
    \begin{tabular}{cc|rrrl|rrr}
    \toprule
          &       & \multicolumn{4}{c|}{Average Time (sec)} & \multicolumn{3}{c}{Opt. Gap (\%)} \\
    \midrule
    Modality & Cap.  & \multicolumn{1}{c}{\tldr} & \multicolumn{1}{c}{\mldr} & \multicolumn{1}{c}{\sddplb} & (opt) & \multicolumn{1}{c}{\tldr} & \multicolumn{1}{c}{\mldr} & \multicolumn{1}{c}{\sddpub} \\
    \midrule
    \multirow{3}[2]{*}{Type-1} & \multicolumn{1}{l|}{20\%} &    5,353 & \textbf{   3,732} &   13,564  & (9)   & 0.62  & 0.39  & \textbf{0.24} \\
          & \multicolumn{1}{l|}{25\%} & \textbf{   4,235} &    4,716  &   13,557  & (6)   & 2.56  & 2.06  & \textbf{1.91} \\
          & \multicolumn{1}{l|}{30\%} & \textbf{   1,555} &    1,675  &     2,067  & (10)  & 0.80  & 0.63  & \textbf{0.27} \\
    \midrule
    \multirow{3}[2]{*}{Type-2} & \multicolumn{1}{l|}{20\%} &    7,489  & \textbf{   4,855} &   18,323  & (5)   & 9.00  & \textbf{7.33} & 9.01 \\
          & \multicolumn{1}{l|}{25\%} &    6,664  & \textbf{   5,560} &   14,790  & (4)   & 4.98  & 4.38  & \textbf{4.28} \\
          & \multicolumn{1}{l|}{30\%} &    1,717  & \textbf{   1,640} &     6,971  & (10)  & 0.81  & 0.63  & \textbf{0.30} \\
    \midrule
    \multicolumn{2}{c|}{Av. (Total)} &    4,502  & \textbf{   3,696} &   11,545  & (44)  & 3.13  & \textbf{2.57} & 2.67 \\
    \bottomrule
    \end{tabular}%
    }
  \label{tab:2sldr-pm-largesize}%
\end{table}%

Table \ref{tab:2sldr-pm-largesize} presents results obtained by these approximation methods over large-size instances. The LDR variant \thldr\ is not included here as the corresponding formulation cannot be loaded into CPLEX due to the large number of first-stage LDR variables and constraints. Also, due to the lack of exact optimal solutions for large-size instances, the optimality gaps are conservative estimations computed using the best lower bound found by \sddplb, that is, we replace  $Obj_\extensive$ with $Obj_\sddplb$.

As in the small-size instances, the 2SLDR variants find high-quality solutions with relatively low computational effort. In fact, \tldr\ and \mldr\ find their respective optimal solutions within the time limit for all the instances,  while \sddplb\ 
proves optimality in 44 of the 60 tested instances (as seen in the ``(opt)" column). Nonetheless, \sddplb\ finds high-quality lower bounds and \sddpub\ yields marginally higher-quality solutions in most cases.
The results suggest that 2SLDR, specifically \mldr, is the best 
for finding high-quality solutions with limited computational effort, while the SDDP-based methods can yield high-quality lower and upper bounds given sufficient 
time.

\begin{remark}
Appendix \ref{app:experiments-exact} analyzes the results of using exact methodologies (i.e., \extensive\ and \sddp) to solve the aggregated models. To summarize, \extensive\ is faster and solves more instances than \sddp\ in small-size instances. However, \extensive\ cannot be solved for large-size problems due to memory requirements (the models use $\sim$16GB RAM after the pre-solve phase of CPLEX), while \sddp\ can solve several instances to optimality and find feasible solutions with small optimality gaps. The performance of \sddp\ is mostly explained by the large branching factor of the MC process and, consequently, the large number of SDDP subproblems (up to 593 for large-size instances). We believe that \sddp\ may be competitive in applications where the MC has fewer reachable states per stage (e.g., ten or fewer).  
\end{remark}

\subsection{Decision Policies and Managerial Insights}
\label{sec:experiments-policy}

We now analyze the quality and structure of the decision policies (in activating contingency modalities) associated with transformations \hn, \ma, \pa\ and \mm, and the original non-aggregated model (\fh). The analysis considers small-size instances (i.e., $4\times5$ grids) since its purpose is to provide insights from the provably optimal solutions, and the non-aggregated model \fh\ can only be solved by \extensive, which cannot handle the large-size instances due to computational limitations. 

Table \ref{tab:policy-quality} summarizes the performances of the considered models in terms of the resulting objective values. Columns ``\% Gap closed'' present the percentage of the gap between \hn\ and \fh\ that is closed by each transformation averaged over all the instances, where the gap closed for each instance is calculated by $(Obj_{\hn} - Obj_{i})/(Obj_{\hn}- Obj_{\fh}),  \forall i \in \{\ma, \pa, \mm\}$. Thus, $0\%$ corresponds to the same objective value given by \hn\, and 100\% corresponds to the objective value given by \fh.

\begin{table}[htbp]
    \def\arraystretch{\stretchTableResults}
  \centering
  \caption{Policy quality comparison for all transformations in our aggregation framework.}
  \scalebox{\scaleTableResults}{
    \begin{tabular}{cl|rrrrr|ccc}
    \toprule
          &       & \multicolumn{5}{c|}{\textbf{Average Objective Value}} & \multicolumn{3}{c}{\textbf{\% Gap closed} 
          } \\
    \midrule
    Modality & Cap.  & \multicolumn{1}{c}{\hn} & \multicolumn{1}{c}{\ma} & \multicolumn{1}{c}{\pa} & \multicolumn{1}{c}{\mm} & \multicolumn{1}{c|}{\fh} & \ma    & \pa    & \mm \\
    \midrule
    \multirow{3}[1]{*}{Type-1} & 20\%  &   104,162  &   102,572  &   92,924  &   92,924  &   82,193  & 7.5   & 51.9  & 51.9 \\
          & 25\%  &     73,879  &     73,442  &   66,650  &   66,238  &   62,972  & 3.2   & 69.4  & 73.1 \\
          & 30\%  &     48,970  &     48,951  &   47,425  &   47,314  &   47,117  & 0.3   & 89.5  & 96.7 \\
    \midrule
    \multirow{3}[1]{*}{Type-2} & 20\%  &   104,135  &   102,625  &   92,612  &   92,612  &   81,182  & 6.7   & 50.9  & 50.9 \\
          & 25\%  &     73,919  &     73,487  &   66,573  &   65,958  &   63,654  & 3.3   & 75.6  & 81.2 \\
          & 30\%  &     48,970  &     48,968  &   47,540  &   47,349  &   47,253  & 0.0   & 80.3  & 97.0 \\
    \bottomrule
    \end{tabular}%
    }
  \label{tab:policy-quality}%
\end{table}%

We first observe that the difference between the objective values given by \hn\ and \fh\ is smaller when the initial DC capacities are higher---about $4\%$ when the initial capacity is 30\%, and about $28\%$ when the initial capacity is $20\%$. This suggests that even the most restrictive transformation \hn\ can lead to high-quality policies when the DCs are capable of addressing random fluctuating demand without much adaptability. Conversely, less restrictive transformations, which yield more adaptive contingency modality activation plans, become valuable when the initial capacity is small. 

Table \ref{tab:policy-quality} also shows that the solution quality of \ma\ and \hn\ are similar (with $< 10\%$ difference) in all instances, indicating that information in the current stage is insufficient to obtain high-quality policies for our test instances. 
This effect is 
mitigated by in \pa\ and \mm\, since 
both can lead to significant performance improvements---closing more than 50\% of the gap between \fh\ and \hn\ in all cases. In particular, the objective values of \pa\ and \mm\ are very close to those obtained by \fh\ when the initial capacity is 30\%, but this is not the case when the initial capacity is 20\%. \pa\ and \mm\ perform similarly, suggesting that the hurricane intensity state captures most of the valuable information 
from the previous stage. In Section \ref{model-implication} we further explore the reasons behind these. 

We now take a closer look at the underlying contingency modality activation policies associated with different transformations. Specifically, we characterize the contingency modality activation by defining some key metrics and report them in Table \ref{tab:solution-structure-type1} for Type-1 instances (see Appendix \ref{app:policy-managerial} for similar results on Type-2 instances). Columns ``Nodes (\%)'' corresponds to the average percentage of nodes in the scenario tree with an active contingency modality, and columns ``\# of Contingencies'' provide the average number of contingency modalities used by a policy among all instances. Columns ``Aggressiveness (\%)'' refer to the percentage of capacity increase (between 10\% and 40\%) resulting from the contingency modality activation, and ``Intensity'' corresponds to the hurricane's intensity when a contingency modality is first activated. We note that the average values in the second and third columns are computed only over nodes with active modalities. For example,
\pa\ only actives one modality for the instance with 20\% initial capacity, which is an aggressive modality (i.e., 30\%) activated when the hurricane reaches intensity three or higher.

\begin{table}[tbp]
  \def\arraystretch{\stretchTableResults}
  \centering
  \caption{Solution structures for different policies and initial capacities for Type-1 instances.}
    \scalebox{\scaleTableResults}{
    \begin{tabular}{l|rrrrr|rrrrr|rrrrr|rrrrr}
    \toprule
          & \multicolumn{5}{c|}{Nodes (\%)			}    & \multicolumn{5}{c|}{\# of Contingencies				} & \multicolumn{5}{c|}{Aggressiveness (\%)		} & \multicolumn{5}{c}{Intensity} \\
    \midrule
    Cap.  & \multicolumn{1}{c}{\hn} & \multicolumn{1}{c}{\ma} & \multicolumn{1}{c}{\pa} & \multicolumn{1}{c}{\mm} & \multicolumn{1}{c|}{\fh} & \multicolumn{1}{c}{\hn} & \multicolumn{1}{c}{\ma} & \multicolumn{1}{c}{\pa} & \multicolumn{1}{c}{\mm} & \multicolumn{1}{c|}{\fh} & \multicolumn{1}{c}{\hn} & \multicolumn{1}{c}{\ma} & \multicolumn{1}{c}{\pa} & \multicolumn{1}{c}{\mm} & \multicolumn{1}{c|}{\fh} & \multicolumn{1}{c}{\hn} & \multicolumn{1}{c}{\ma} & \multicolumn{1}{c}{\pa} & \multicolumn{1}{c}{\mm} & \multicolumn{1}{c}{\fh} \\
    \midrule
    20\%  & 100   & 93    & 72    & 72    & 53    & 1.0   & 1.0   & 1.0   & 1.0   & 3.5   & 16    & 38    & 30    & 30    & 35    & 2.3   & 2.5   & 2.9   & 2.9   & 2.9 \\
    25\%  & 10    & 50    & 22    & 22    & 47    & 1.0   & 1.0   & 1.5   & 2.1   & 4.4   & 10    & 15    & 35    & 36    & 21    & 2.3   & 2.9   & 4.0   & 4.0   & 3.0 \\
    30\%  & 0     & 7     & 16    & 15    & 19    & -     & 1.0   & 2.0   & 4.3   & 4.6   & -     & 9     & 17    & 18    & 15    & -     & 3.0   & 4.1   & 4.1   & 4.0 \\
    \bottomrule
    \end{tabular}%
    }
  \label{tab:solution-structure-type1}%
\end{table}%

From Table \ref{tab:solution-structure-type1}, we see that 
a lower initial capacity typically leads to more contingency modality activation and modality choices that are more aggressive. Comparing \hn\ with other transformations, we see that it has more active nodes when the initial capacity is lower, and has active fewer nodes when the initial capacity is higher. This is attributed to the structure of the policy under \hn\, which resembles an ``all-or-nothing'' feature due to the lack of adaptability. The behaviors of \hn\ and \ma\ are quite similar, although the activation timing of  \ma\ is later, resulting in a different number of nodes with an active modality. We attribute this similarity to the ``propagation effect'' in MC-based aggregation policies, which we discussed in Section~\ref{model-implication}. In contrast, \pa\ and \mm\ result in more adaptive policies that resemble a ``wait-and-activate'' type of behavior, similar to \fh. In these more adaptive policies, contingency modalities tend to be activated when the hurricane reaches a higher intensity state. In addition, the more adaptive policies lead to a more diverse set of modality choices and a better balance between aggressiveness and the number of nodes actives depending on the initial capacity.
Comparing the MC-based transformations \pa\ and \mm\ with \fh, which has full adaptability, we see the value of additional adaptability in activating less nodes throughout the scenario tree with more diverse contingency modalities with appropriate aggressiveness.

\subsection{Implications of the MC-based Aggregation in HDR}\label{model-implication}

We now analyze impacts of the MC-based aggregation on the resulting HDR solution policies. 
In particular, we study why the number of possible contingency modalities is low for some of the transformations (e.g., at most one for \ma) in our test instances. First, recall that binary variables $z_{n\ell}$ represent whether contingency modality $\ell \in \modality$ is active in node $n \in \tree$. In addition, the model imposes that only one modality is active in each state \eqref{dr_node:modality} and every active modality remains active in the following stages \eqref{dr_node:modality2}. Thus, these constraints restrict the number of different modalities that can be active in any solution (including \fh), \Merve{especially }
when an aggregation is applied. 

\begin{figure}
    \centering
    \begin{tikzpicture}[->,>=stealth',shorten >=1pt,auto,node distance=1cm,
thick]       
\node[node light] (r) at (0,0) {};
\node[node dark] (n11)  at (0.8,0.8)  {};
\node[node light] (n12)  at (0.8,-0.8)  {};

\node[node dark] (n21)  at (1.6,1.2)  {};
\node[node light] (n22)  at (1.6,0.4)  {};
\node[node dark] (n23)  at (1.6,-0.4)  {};
\node[node light] (n24)  at (1.6,-1.2)  {};

\node[node dark] (n31)  at (2.4,1.4)  {};
\node[node light] (n32)  at (2.4,1)  {};
\node[node dark] (n33)  at (2.4,0.6)  {};
\node[node light] (n34)  at (2.4,0.2)  {};
\node[node dark] (n35)  at (2.4,-0.2)  {};
\node[node light] (n36)  at (2.4,-0.6)  {};
\node[node dark] (n37)  at (2.4,-1)  {};
\node[node light] (n38)  at (2.4,-1.4)  {};

\node[] (name) at (1.3, 2) {Scenario Tree};

\node[text node] (name) at (-0.1, -2) {$t=1$};
\node[text node] (name) at (0.8, -2) {$2$};
\node[text node] (name) at (1.6, -2) {$3$};
\node[text node] (name) at (2.4, -2) {$4$};

\path[every node/.style={font=\sffamily\small}]
(r) 
edge[arc] node [left] {} (n11)
edge[arc] node [left] {} (n12)
(n11)
edge[arc] node [left] {} (n21)
edge[arc] node [left] {} (n22)
(n12)
edge[arc] node [left] {} (n23)
edge[arc] node [left] {} (n24)
(n21)
edge[arc] node [left] {} (n31)
edge[arc] node [left] {} (n32)
(n22)
edge[arc] node [left] {} (n33)
edge[arc] node [left] {} (n34)
(n23)
edge[arc] node [left] {} (n35)
edge[arc] node [left] {} (n36)
(n24)
edge[arc] node [left] {} (n37)
edge[arc] node [left] {} (n38)
;

\end{tikzpicture}
    \hfill
    \begin{tikzpicture}[->,>=stealth',shorten >=1pt,auto,node distance=1cm,
thick]       
\node[node light] (r) at (0,0) {};
\node[node-square, draw=red, fill=red] (n11)  at (0.8,0.8)  {};
\node[node light] (n12)  at (0.8,-0.8)  {};

\node[node-polygon5, draw=green!70!black, fill=green!70!black] (n21)  at (1.6,1.2)  {};
\node[node-polygon5, draw=green!70!black, fill=green!70!black] (n22)  at (1.6,0.4)  {};
\node[node-polygon6, draw=orange, fill=orange] (n23)  at (1.6,-0.4)  {};
\node[node-polygon6, draw=orange, fill=orange] (n24)  at (1.6,-1.2)  {};

\node[node-polygon5, draw=green!70!black, fill=green!70!black] (n31)  at (2.4,1.4)  {};
\node[node-polygon5, draw=green!70!black, fill=green!70!black] (n32)  at (2.4,1)  {};
\node[node-polygon5, draw=green!70!black, fill=green!70!black] (n33)  at (2.4,0.6)  {};
\node[node-polygon5, draw=green!70!black, fill=green!70!black] (n34)  at (2.4,0.2)  {};
\node[node-polygon6, draw=green!70!black, fill=green!70!black] (n35)  at (2.4,-0.2)  {};
\node[node-polygon6, draw=green!70!black, fill=green!70!black] (n36)  at (2.4,-0.6)  {};
\node[node-polygon6, draw=green!70!black, fill=green!70!black] (n37)  at (2.4,-1)  {};
\node[node-polygon6, draw=green!70!black, fill=green!70!black] (n38)  at (2.4,-1.4)  {};

\node[] (name) at (1.1, 2) {\ma};

\node[text node] (name) at (-0.1, -2) {$t=1$};
\node[text node] (name) at (0.8, -2) {$2$};
\node[text node] (name) at (1.6, -2) {$3$};
\node[text node] (name) at (2.4, -2) {$4$};

\path[every node/.style={font=\sffamily\small}]
(r) 
edge[arc] node [left] {} (n11)
edge[arc] node [left] {} (n12)
(n11)
edge[arc] node [left] {} (n21)
edge[arc] node [left] {} (n22)
(n12)
edge[arc] node [left] {} (n23)
edge[arc] node [left] {} (n24)
(n21)
edge[arc] node [left] {} (n31)
edge[arc] node [left] {} (n32)
(n22)
edge[arc] node [left] {} (n33)
edge[arc] node [left] {} (n34)
(n23)
edge[arc] node [left] {} (n35)
edge[arc] node [left] {} (n36)
(n24)
edge[arc] node [left] {} (n37)
edge[arc] node [left] {} (n38)
;

\end{tikzpicture}
    \hfill
    \begin{tikzpicture}[->,>=stealth',shorten >=1pt,auto,node distance=1cm,
thick]       
\node[node light] (r) at (0,0) {};
\node[node-square, draw=red, fill=red] (n11)  at (0.8,0.8)  {};
\node[node light] (n12)  at (0.8,-0.8)  {};

\node[node-polygon5, draw=green!70!black, fill=green!70!black] (n21)  at (1.6,1.2)  {};
\node[node-polygon5, draw=green!70!black, fill=green!70!black] (n22)  at (1.6,0.4)  {};
\node[node dark] (n23)  at (1.6,-0.4)  {};
\node[node light] (n24)  at (1.6,-1.2)  {};

\node[node-polygon5, draw=green!70!black, fill=green!70!black] (n31)  at (2.4,1.4)  {};
\node[node-polygon5, draw=green!70!black, fill=green!70!black] (n32)  at (2.4,1)  {};
\node[node-polygon5, draw=green!70!black, fill=green!70!black] (n33)  at (2.4,0.6)  {};
\node[node-polygon5, draw=green!70!black, fill=green!70!black] (n34)  at (2.4,0.2)  {};
\node[node-polygon5, draw=orange, fill=orange] (n35)  at (2.4,-0.2)  {};
\node[node-polygon5, draw=orange, fill=orange] (n36)  at (2.4,-0.6)  {};
\node[node-polygon5, draw=orange, fill=orange] (n37)  at (2.4,-1)  {};
\node[node-polygon5, draw=orange, fill=orange] (n38)  at (2.4,-1.4)  {};

\node[] (name) at (1.1, 2) {\mm};

\node[text node] (name) at (-0.1, -2) {$t=1$};
\node[text node] (name) at (0.8, -2) {$2$};
\node[text node] (name) at (1.6, -2) {$3$};
\node[text node] (name) at (2.4, -2) {$4$};

\path[every node/.style={font=\sffamily\small}]
(r) 
edge[arc] node [left] {} (n11)
edge[arc] node [left] {} (n12)
(n11)
edge[arc] node [left] {} (n21)
edge[arc] node [left] {} (n22)
(n12)
edge[arc] node [left] {} (n23)
edge[arc] node [left] {} (n24)
(n21)
edge[arc] node [left] {} (n31)
edge[arc] node [left] {} (n32)
(n22)
edge[arc] node [left] {} (n33)
edge[arc] node [left] {} (n34)
(n23)
edge[arc] node [left] {} (n35)
edge[arc] node [left] {} (n36)
(n24)
edge[arc] node [left] {} (n37)
edge[arc] node [left] {} (n38)
;

\end{tikzpicture}
    \hfill
    \begin{tikzpicture}[->,>=stealth',shorten >=1pt,auto,node distance=1cm,
thick]       
\node[node light] (r) at (0,0) {};
\node[node-square, draw=red, fill=red] (n11)  at (0.8,0.8)  {};
\node[node light] (n12)  at (0.8,-0.8)  {};

\node[node-polygon5, draw=green!70!black, fill=green!70!black] (n21)  at (1.6,1.2)  {};
\node[node-polygon5, draw=green!70!black, fill=green!70!black] (n22)  at (1.6,0.4)  {};
\node[node-square, draw=red, fill=red] (n23)  at (1.6,-0.4)  {};
\node[node light] (n24)  at (1.6,-1.2)  {};

\node[node-polygon5, draw=green!70!black, fill=green!70!black] (n31)  at (2.4,1.4)  {};
\node[node-polygon5, draw=green!70!black, fill=green!70!black] (n32)  at (2.4,1)  {};
\node[node-polygon5, draw=green!70!black, fill=green!70!black] (n33)  at (2.4,0.6)  {};
\node[node-polygon5, draw=green!70!black, fill=green!70!black] (n34)  at (2.4,0.2)  {};
\node[node-polygon5, draw=green!70!black, fill=green!70!black] (n35)  at (2.4,-0.2)  {};
\node[node-polygon5, draw=green!70!black, fill=green!70!black] (n36)  at (2.4,-0.6)  {};
\node[node dark] (n37)  at (2.4,-1)  {};
\node[node light] (n38)  at (2.4,-1.4)  {};

\node[] (name) at (1.1, 2) {\fh};

\node[text node] (name) at (-0.1, -2) {$t=1$};
\node[text node] (name) at (0.8, -2) {$2$};
\node[text node] (name) at (1.6, -2) {$3$};
\node[text node] (name) at (2.4, -2) {$4$};

\path[every node/.style={font=\sffamily\small}]
(r) 
edge[arc] node [left] {} (n11)
edge[arc] node [left] {} (n12)
(n11)
edge[arc] node [left] {} (n21)
edge[arc] node [left] {} (n22)
(n12)
edge[arc] node [left] {} (n23)
edge[arc] node [left] {} (n24)
(n21)
edge[arc] node [left] {} (n31)
edge[arc] node [left] {} (n32)
(n22)
edge[arc] node [left] {} (n33)
edge[arc] node [left] {} (n34)
(n23)
edge[arc] node [left] {} (n35)
edge[arc] node [left] {} (n36)
(n24)
edge[arc] node [left] {} (n37)
edge[arc] node [left] {} (n38)
;

\end{tikzpicture}
    \caption{Activation solutions over an abstract scenario tree for different MC-based transformations.}
    \label{fig:solution_restricction}
\end{figure}
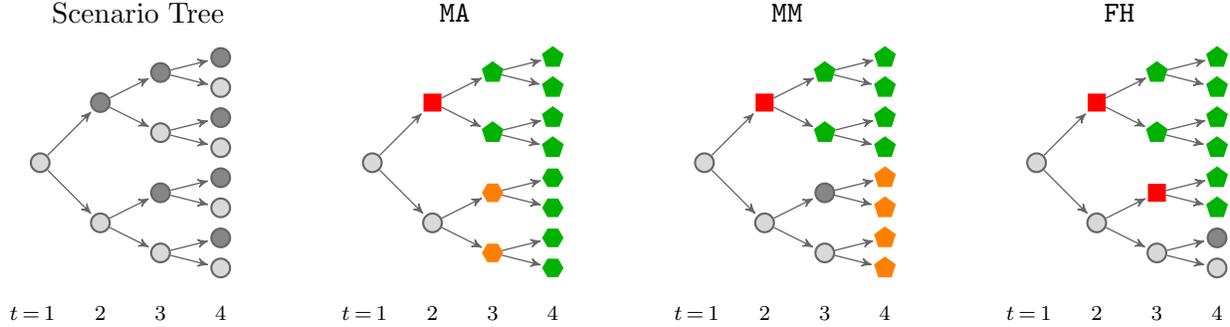

To understand the implications of an MC-based aggregation over nodes in the scenario tree with active modalities, Figure \ref{fig:solution_restricction} illustrates solutions of different transformations in an abstract scenario tree with two states 
(from Example) \ref{example:scenario_tree} 
Specifically, \squareRed{1} are nodes where a modality is \textit{initially activated} (as opposed to passively activated 
due to propagation, as explained below), \pentaGreen{1} are nodes that have active modalities because of prior activation at their parent nodes, and \hexaOrange{1} are nodes where the modality is \textit{passively activated} (i.e., the activation is a consequence of the propagation caused by the transformation that forces certain nodes to have the same $z$ values). For example, we can see that the \ma\ policy 
initially activates a modality in $t=2$ at the dark node. However, this causes the two \hexaOrange{1} nodes to be passively activated because nodes with the same MC state must have the same $z$ values according to transformation \ma. This propagation, as a result of the MC-based aggregation, may cause an initial activation to have ``unintended consequences'' of passively activating nodes that may not necessarily benefit from the modality activation. This explains the similarity between \hn\ and \ma\ that we observe in Section~\ref{sec:experiments-policy}. We can also observe such propagation behaviors with \pa\ and \mm\, but they occur less often since the aggregation depends on the MC state of the current stage and the previous one. In contrast, the activation when using \fh\ does not have any such propagation effect since no aggregation is imposed (i.e., no \hexaOrange{1} can be observed from \fh). 

\begin{figure}[tb]
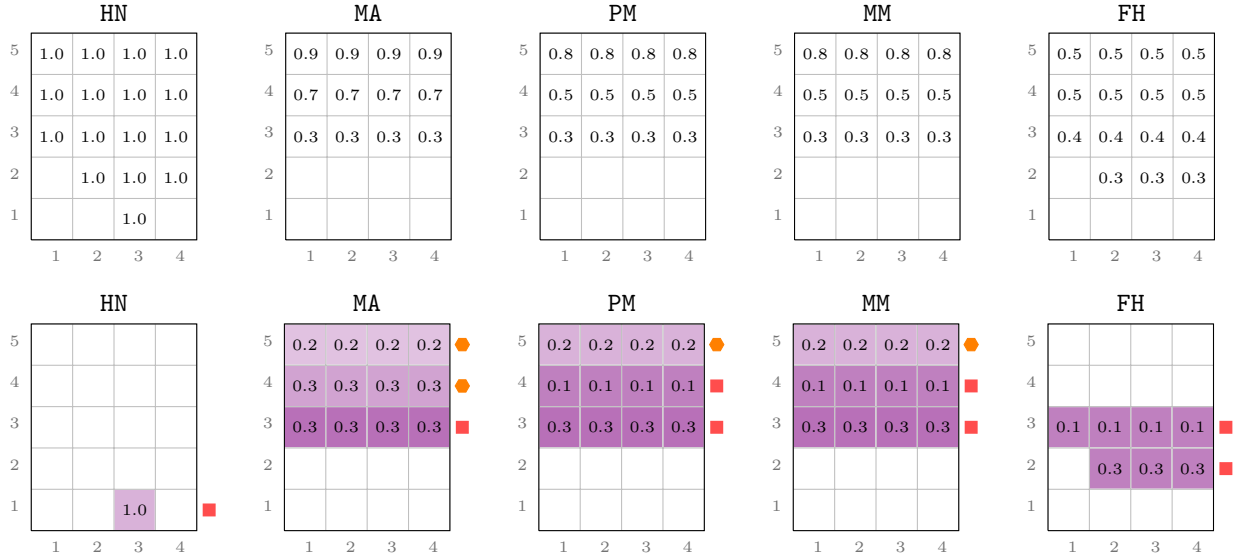

    \centering
    \heatmap{figures/heatmaps/hn-active.txt}{figures/heatmaps/none-intensity.txt}{\hn}{1}
    \hfill
    \heatmap{figures/heatmaps/ma-active.txt}{figures/heatmaps/none-intensity.txt}{\ma}{1}
    \hfill
    \heatmap{figures/heatmaps/pa-active.txt}{figures/heatmaps/none-intensity.txt}{\pa}{1}
    \hfill
    \heatmap{figures/heatmaps/mm-active.txt}{figures/heatmaps/none-intensity.txt}{\mm}{1}
    \hfill
    \heatmap{figures/heatmaps/fh-active.txt}{figures/heatmaps/none-intensity.txt}{\fh}{1}
    \vspace{0.5em}
    
    \heatmapColor{figures/heatmaps-color/hn-active.txt}{figures/heatmaps-color/hn-intensity.txt}{\hn}
    \hfill
    \heatmapColor{figures/heatmaps-color/ma-active.txt}{figures/heatmaps-color/ma-intensity.txt}{\ma}
    \hfill
    \heatmapColor{figures/heatmaps-color/pa-active.txt}{figures/heatmaps-color/pa-intensity.txt}{\pa}
    \hfill
    \heatmapColor{figures/heatmaps-color/mm-active.txt}{figures/heatmaps-color/mm-intensity.txt}{\mm}
    \hfill
    \heatmapColor{figures/heatmaps-color/fh-active.txt}{figures/heatmaps-color/fh-intensity.txt}{\fh}
    
    \caption{Graphical representation of the solutions of an instance of Type-1 with 20\% initial capacity.}
    \label{fig:heatmap}
\end{figure}

Lastly, Figure \ref{fig:heatmap} illustrates the resulting policies for a single instance in a series of images corresponding to the grid of the problem where the hurricane originates at $y=0$ (i.e., bottom row) and lands at $y=5$ (i.e., top row). The top grids show the percentage of nodes in the scenario tree with an active contingency modality for each cell. The bottom grids are heat-maps showing the proportion of nodes in a cell where a contingency modality is the first activated (either initially activated or passively activated), and the color scheme corresponds to the hurricane intensity level (i.e., darker colors for higher intensity). Each row is associated with a symbol \squareRed{1} or \hexaOrange{1}, which distinguishes if the nodes are initially  or passively activated, respectively. We observe that \ma, \pa, and \mm\ lead to a wait-and-activate type behavior by activating modalities in high-intensity nodes at later stages to reduce the number of active modalities in the last stage, a similar strategy to \fh. In contrast, \hn\ activates a modality in the first stage even though it may not be necessary for many nodes in later stages. We also observe that the intensity levels of \hexaOrange{1} nodes tend to be much smaller than the ones associated with \squareRed{1} nodes. This reflects the ``unintended consequence'' of forcing nodes of low intensity to be passively activated due to the restrictions imposed by the transformations. 


    
    

\section{Conclusions} \label{sec:conclusions}

We present an aggregation framework to address MSILP problems with mixed-integer state variables and continuous local variables. The main idea of the framework is to reformulate the original MSILP such that all integer state variables are first-stage and then apply a suitable transformation to aggregate these variables to reduce the problem size while not sacrificing the solution quality. We describe a number of transformations that rely on the underlying structure of the stochastic process, which we assume is given by an MC. We demonstrate an exact solution approach based on a B\&C framework integrated with the SDDP algorithm to solve the resulting aggregated MSILP. Also, we propose a computationally more tractable alternative based on the 2SLDR approximations and introduce new MC-based decision rules. These approaches are tested on a hurricane disaster relief logistics application with contingency modality decisions. Empirical results show the trade-off between different transformations regarding the solution quality and computational effort, and the best trade-off is yielded by transformations that consider information from current and previous MC states (i.e., \pa). The SDDP lower and upper bounding procedures have the best performance in terms of solution quality for large-scale instances, while the MC-based 2SLDR approach returns high-quality solutions with significantly less computational effort.  

As shown by our empirical results, the performance of our aggregation framework largely depends on the chosen transformation. With this in mind, one future research direction is to investigate how to choose an appropriate transformation for a given problem instance. \Song{An alternative research direction is to develop a dynamic extension of the current approach to iteratively refine the aggregation over time, e.g., in a rolling horizon procedure.}

%
%
%





{
\bibliographystyle{informs2014} 
\bibliography{main} 
}



\newpage
\begin{APPENDICES}

\section{B\&C Integrated with the SDDP algorithm} 
\label{app:bc_sddp}

Algorithm \ref{alg:bc-sddp} describes the B\&C procedure. The procedure receives as input a few parameters used by an SDDP sub-routine to be explained in detail later. Note that lines 4-7 correspond to the standard B\&B algorithm of commercial solvers.

\begin{algorithm}[H] 
\small
    \linespread{1.2}\selectfont
	\caption{A B\&C framework integrated with the SDDP algorithm} \label{alg:bc-sddp}
	\begin{algorithmic}[1]
		\Procedure{\algBC}{$\algExact, K, \epsilon$}
		\State Initialize B\&B search for $\Qsub_\rootnode$ (master problem) and the \algSDDP\ environment.
		\State Initialize $\coeffs_n :=\emptyset$ for all $n \in \childs(\rootnode)$ and $\algCutAdd := \algFalse$.
		\Repeat
		\State Choose a B\&B node from its list of open nodes. Solve the associated node relaxation problem. 
		\If {Current B\&B node is infeasible or the relaxation bound is worse than the bound given by an incumbent solution}
		\State Prune this B\&B node.  
		\EndIf
		\If {Values $\zasol$ of the current solution $(\xrootsol, \zasol, \thetasol)$ are integer} 
		\Repeat {\hfill \color{gray} \%\% Cutting plane subroutine}
		\State $\algCutAdd := \algFalse$.
		\For {$n' \in \childs(\rootnode)$} {\hfill \color{gray} \%\% Find cuts parameters for each child subproblem}
		\State Get cut coefficients $(\alpha_{n'},\beta_{n'},\gamma_{n'})$ := \algSDDP($(\xrootsol, \zasol, \thetasol), n',\epsilon,K,\algExact$).
		\If{Cut with coefficients $(\alpha_{n'},\beta_{n'},\gamma_{n'})$ is violated by  $(\xrootsol, \zasol, \thetasol)$}
		\State $\coeffs_{n'}:= \coeffs_{n'}\cup(\alpha_{n'},\beta_{n'},\gamma_{n'})$, $\algCutAdd:= \algTrue$.
		\EndIf
		\EndFor
        \Until{$\algCutAdd = \algFalse$.}
        \Else
        \State Branch on the current  B\&B node and update the list of open nodes.
		\EndIf
		\Until{The list of open B\&B nodes is empty.}
		\EndProcedure
	\end{algorithmic} 
\end{algorithm}

There are a few important points to mention about this procedure. First, we call the SDDP sub-routine for each child of the root node $\rootnode$, in which we perform the SDDP forward pass only considering sample paths that contain that specific child. Also, we look for cuts for all children nodes of the root node $\rootnode$ before re-solving the current B\&B node relaxation. An alternative would be to re-solve the B\&B node relaxation whenever we find a violated cut, but this option showed worse performance in our preliminary experiments. Lastly, all cuts added to enhance the outer approximation of the cost-to-go functions associated with the SDDP subproblems in the previous iterations are automatically carried over to the current iteration. 

\section{Benders Decomposition for the  2SLDR Model} \label{app:2sldr_benders}

In this section we present a Benders decomposition approach to solve \eqref{mod:ldr}. We consider a Benders master problem that handles the first-stage decisions and decompose the second-stage problem into one Benders subproblem per node. Because the number of subproblems can be exponentially many (e.g., when \thldr\ is used), we consider an aggregated approach where we have one variable $\theta_{t,m}^L$ that approximates all the cost-to-go functions associated with nodes in stage $t$ and with MC state $m \in \markov_t$. This variant can be seen as a hybrid between the single-cut version (i.e., using a single $\theta$ variable) and the multi-cut version (i.e., using one $\theta_{t,n}^L$ variable for each subproblem) of Benders decomposition \citep{van1969shaped,birge1988multicut}. Our preliminary experiments show that our Benders variant performs the best on our test instances, however, other variants of Benders decomposition and computational enhancements can also be considered (see, e.g., \cite{zverovich2012computational} and \cite{bodur2017strengthened}). 

In what follows we present the details of our Benders decomposition for \thldr\ (the other two cases, \tldr\ and \mldr\, are similar). The Benders master problem is given by
\begin{align}
    \Qldrsub_{\rootnode} = \min\; & 
    \sum_{n\in \nodes}\prob_n \left( c_n^\top \zphi + d_n^\top(\varldr_{t(n)}^\top\xi^{t(n)}_n)  \right) + h_\rootnode^\top \yroot  + \sum_{t = 2}^\periods\sum_{m \in \markov_t}\theta_{t,m}^L \nonumber\\ 
    \text{s.t.}\; &  \eqref{eq:ldr-constr0}-\eqref{eq:ldr-constr2}, \nonumber \\
    & \theta^L_{t,m} \geq \alpha^L_{t,m}{}^\top \mu + \beta^L_{t,m}{}^\top \za + \gamma^L_{t,m}, \quad \forall  (\alpha^L_{t,m},\beta^L_{t,m},\gamma^L_{t,m})\in \coeffsldr_{t,m}, \ m \in  \markov_t, \ t \in \{2,...,\periods\}, \label{eq:ldr_benderscut}\\
    & (\yroot, \varldr, \za)  \in  \R^r \times \R^{k\times \sum_{t\in \horizon}\dimxitt} \times \Z^{\ell \cdot \sum_{t \in \horizon} q_t}, \nonumber
\end{align}
where $\theta^L \in \R^{\sum_{t=2}^{\periods}|\markov_t|}$ is introduced to represent the cost-to-go function approximation for each Markovian state in all remaining stages (except for the first stage) via constraints~\eqref{eq:ldr_benderscut}. Here $\coeffsldr_{t,m}$ stores the coefficients associated with the Benders cuts for each MC state $m \in \markov_t$. 

\begin{algorithm}[htb]
\small
    \linespread{1.2}\selectfont
	\caption{Benders Decomposition} \label{alg:benders}
	\begin{algorithmic}[1]
		\Procedure{\algBenders}{$\epsilon$}
		\State Initialize B\&B search for $\Qldrsub_\rootnode$ (master problem)
		\Repeat
		\State Choose a node B\&B from its node list. Solve the associated node relaxation problem. 
		\If {Current B\&B node is infeasible or the relaxation bound is worse than the bound given by an incumbent solution}
		\State Prune this B\&B node.  
		\EndIf
		\If {Current solution $(\varldrsol, \zasol, \thetasol^L)$ is integer} 
		\Repeat {\hfill \color{gray} \%\% Cutting plane subroutine}
		\For {$m \in \markov_t$ and $t \in \{2,...,\periods\}$}
		\For {$n \in \staget(m)$}
		\State Solve $\Qldrsol_n:=\Qldr_n(\varldrsol, \zasol)$ and save dual solutions ${\pi^{n}}$. 
		\EndFor
		\State Save value $\Qldrsol_{t,m}:=\sum_{n \in \staget(m)} p_n\Qldrsol_n$.
		\If{$|\Qldrsol_{t,m} - \thetasol^L_{t,m}| \geq \epsilon |\Qldrsol_{t,m}|$} {\hfill \color{gray} \%\% Add cut if needed}
		\State Add inequality \eqref{eq:ldr_benderscut} and exit the for loop. 	
		\EndIf
		\EndFor
        \Until{No new inequality is added to the master problem $\Qldrsub_\rootnode$.}
        \Else
        \State Branch on the current  B\&B node and update list of nodes.
		\EndIf
		\Until{No more nodes in the B\&B node list.}
		\EndProcedure
	\end{algorithmic} 
\end{algorithm}

Algorithm \ref{alg:benders} describes our Benders decomposition for \eqref{mod:ldr}. The procedure is a B\&C algorithm similar to the one described in Algorithm \ref{alg:bc-sddp}. Once an integer solution is found during B\&B search, we enter the cutting plane subroutine. The algorithm iterates over each stage and MC state and its corresponding nodes, where $\staget(m)$ is the set of all nodes in stage $t$ associated with MC state $m \in \markov_t$. It then evaluates the Benders subproblems for all the nodes associated with a particular $m$ and $t$ and checks if the current approximation $\thetasol^L_{t,m}$ is close enough to the actual value $\Qldrsol_{t,m}$ for a given tolerance $\epsilon$. If the approximation is not good enough, then the algorithm returns a cut and stops evaluating other subproblems. Similar to our SDDP variant, we add cuts to each candidate integer solution until no more cuts are found. The complete procedure ends when there are no more nodes to explore in the B\&B search tree, i.e., we have found an optimal solution or proven infeasibility.

The cuts added to the master problem \eqref{eq:ldr_benderscut} are Benders cuts over all the nodes associated to a given stage and MC state. Specifically, the cut coefficients for $t \in \{2,...,T\}$ and $m \in \markov_t$ are 
\begin{align*}
    \alpha_m^L & = \sum_{n \in \staget(m)}\prob_n {\pi^{n}}^\top  A_{n}, \qquad \beta_m =  \sum_{n \in \staget(m)} \prob_n {\pi^{n}}^\top B_{n}, \\
    \gamma_m  &= \sum_{n \in \staget(m)} \prob_n \left( \Qldrsol_n - ({\pi^{n}}^\top  A_{n} )^\top \hat{\mu}_{t(n) -1}- ({\pi^{n}}^\top B_{n})^\top \zasol_{\phi_t(a(n))} \right),
\end{align*}
where ${\pi^{n}}$ are the dual solutions associated with the cost-to-go function $\Qldr(n)$. Note that each coefficient here corresponds to a weighted sum of the associated coefficients in the subproblems. Also, the algorithm generates Benders feasibility cuts when necessary, which is common in 2SLDR approximations \citep{bodur2018two}.

\vspace{1em}
\section{Aggregated Model for HDR Applicaton} \label{app:hdr_aggregated}

We now present how to transform the MSILP model \eqref{dr_root} with integer state variables at each node of the scenario tree to an MSILP model with integer variables only at the root node, following the ideas presented in Section~\ref{sec:framework} and Section~\ref{sec:methodology}. The general idea is to consider all integer state variables as first-stage variables and apply a suitable transformation $\maptmarkov$ to aggregate the number of variables to a manageable size. In what follows, we present the resulting aggregated model \eqref{dr_aggre} for a generic transformation and then show four alternatives given the specific structure of \eqref{dr_root}. 

Following the notation introduced in Section \ref{sec:framework}, we consider $\za$ to be the set of aggregated integer state variables for each node. There are two aspects to consider when applying this transformation. First, constraints \eqref{dr_node:modality} and \eqref{dr_node:modality2} are now first-stage constraints. Second, we can omit state variables $\xcap_n$  because these variables are only used to model the capacity dependency from one stage to the next, which is fully dictated by the first-stage variables $\za$. The resulting aggregated model is:
\vspace{-0em}
\begin{subequations}\label{aggregated-HDR}
\begin{align}
      \min \ &  \sum_{j \in \centers}\left( g_{j}\xinv_{\rootnode j} + q_{\rootnode j} v_{j} + \sum_{i \in \shelters}f_{\rootnode i j} y_{ij} \right) + \sum_{i \in \shelters} b_{i}w_{i} + \sum_{n \in \nodes}\sum_{\ell \in \modality} \lefteqn{
      \za_{\phi_t(n) \ell} c_{\ell} + \sum_{n \in \childs(\rootnode)}\pprob_{\rootnode n} \Qaggre_n(\xroot^I, \za) } \tag{$\textit{HDR}^A$} \label{dr_aggre} \\
    \text{s.t.} \ & \eqref{dr_root:demand}-\eqref{dr_root:inventory}  \nonumber \\
    & v_{j} \leq C_j,  & \forall j \in \centers, \label{dr_aggre:capacity}\\
    & \sum_{\ell \in \modality} \za_{\phi_t(n) \ell} \leq 1, & \forall n\in \nodes, \label{dr_aggre:modality1} \\
    & \za_{\phi_{t-1}(a(n)) \ell} \leq  \za_{\phi_t(n) \ell}, & \forall \ell \in \modality, n \in \nodesnoroot, \label{dr_aggre:modality2} \\ 
    & \xinv_{\rootnode j}, v_{j}, w_{i}, y_{ij} \geq 0, \; \za_{\phi_t(n) \ell} \in \{0,1\}, & \forall j \in \centers, \; i \in \shelters, \ell \in \modality, n \in \nodesnoroot, \nonumber 
\end{align}
\end{subequations}
where the cost-to-go function for a node $n\in \nodesnoroot$ is given by:
\vspace{-0em}
\begin{subequations}
\begin{align}
    \Qaggre_n(\xnparent^I, \za)  =  \min \ & \sum_{j \in \centers}\left( g_{j}\xinv_{n j} + q_{n j} v_{j} + \sum_{i \in \shelters} f_{nij}y_{ij} \right) +  \sum_{i \in \shelters} b_{i}w_{i} + \lefteqn{\sum_{n' \in \childs(n)}\pprob_{n n'} \Qaggre_{n'}(\xn^I, \za)} \nonumber \\
    \text{s.t.}\; & \eqref{dr_node:demand}-\eqref{dr_node:inventory}  \nonumber \\
    & v_{j} \leq C_j + \sum_{n' \in \patht(n)}\sum_{\ell \in \modality} K_{j \ell}\za_{\phi_t(n') \ell},  & \forall j \in \centers, \label{dr_aggre:capacity2}\\
    & \xinv_{n j}, v_{j}, w_{i}, y_{ij} \geq 0, & \forall j \in \centers, \; i \in \shelters. \nonumber 
\end{align}
\end{subequations}
Note that the right-hand side of constraint \eqref{dr_aggre:capacity2} represents the current capacity of a DC by considering the modality activation and capacity increases at previous stages. Also, although constraints \eqref{dr_aggre:modality1} and \eqref{dr_aggre:modality2} are imposed for each node, the fact is that some of them might be identical because they share the same problem parameters and aggregated variables. For example, transformation \hn\ only needs one set of constraints \eqref{dr_aggre:modality1} per stage because of the stage-based variable aggregation.

\vspace{1em}
\section{Problem Description and Instance Generation for HDR} \label{app:hdr}

We now describe the main components of our HDR problem and present details of the instance generation procedure. We  focus on the grid representation for the problem, the MC, and three main parameters (i.e., demand, modalities, and capacity). All other details about instances generation (e.g., cost parameters) can be found in our instances generation code. We will make the instance generation code and the set of instances publicly available upon publication. 

\subsection{Grid Representation and Locations of Shelters and DCs}

Similar to previous works in the literature, we use a grid network to represent the potential locations of shelters and DCs, as well as possible locations of the hurricane (see Figure~\ref{fig:hurricane_grid}). The top row on the grid represents the land, where each cell contains a subset of shelters and DCs. The remaining cells in the grid correspond to the sea and are used to represent possible locations of the hurricane at any stage during the planning horizon. For simplicity, we assume that any shelter can be supplied by any DC in the network.

\begin{figure}[tb]
    \centering
    \begin{tikzpicture}
	\draw[fill=blue!10!white] (0,0) -- (0,3) -- (7,3) -- (7,0) -- (0,0); 
	\draw[fill=brown!30!white] (0,3) -- (0,4) -- (7,4) -- (7,3) -- (0,3); 
	\draw[step=1,gray!50!white,thin] (0,0) grid (7,4);
	\draw[black] (0,0) -- (0,4) -- (7,4) -- (7,0) -- (0,0); 
    \foreach \x in {0,...,6}
        \foreach \y in {0,...,3} 
            {\pgfmathtruncatemacro{\label}{\x +  \y*3}
                \node[gray!90!black]  (\x\y) at (0.5 +\x, 0.5 + \y) {\scriptsize \x ,\y};} 
    \node (n) at (-4.0, 0.5) {$\quad$};
    \node (n) at (9.5, 0.5) {$\longrightarrow$ Hurricane starts here};
    
    \draw [ decorate, decoration = {calligraphic brace,
        raise=5pt, amplitude=5pt}] (0,0) --  (0,3)
        node[pos=0.5,left=10pt,black]{Sea};
        
    \draw [ decorate, decoration = {calligraphic brace,
        raise=5pt, amplitude=5pt}] (0,3) --  (0,4)
        node[pos=0.5,left=10pt,black]{Land};
        
    \draw[violet, fill=violet!70!white, opacity=0.5] (5.5,0.5) circle (0.3); 
    \filldraw[violet!70!white, opacity=0.2] (5.3,0.7) -- (2.5,3.5) -- (6.5, 3.5) -- (6.5, 1.5) -- (5.7,0.7);
    
\end{tikzpicture}
    \caption{A $6\times 4$ grid example that shows the initial location of the hurricane $(m^x,m^y)=(5,0)$ and all possible locations in the following stages (i.e., shaded region).}
    \label{fig:hurricane_grid}
\end{figure}
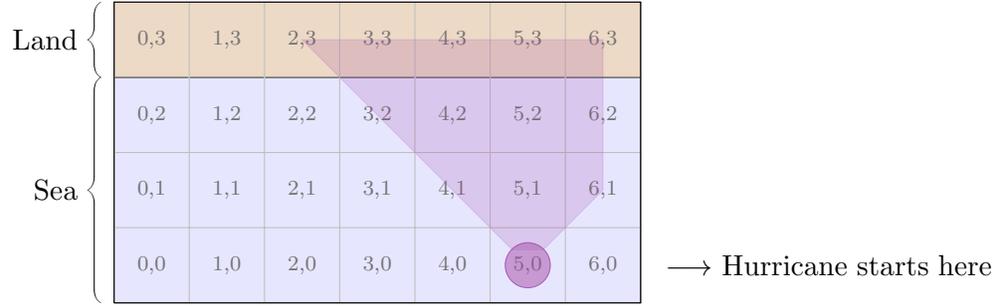

Land and sea cells have pre-defined dimensions: land cells have a size of $50\times100$, while sea cells have a size of $20\times100$. We randomly generate 3 to 7 shelters and 2 to 4 DCs in each land cell. Shelters and DCs are positioned uniformly at random inside their corresponding cells. The coordinates of the hurricane located in a cell are assumed to be the coordinates of the center of the cell.

\subsection{MC description}

As described in the main text, each MC state corresponds to two attributes of the hurricane: location and intensity, that is, $m=(m^x, m^y, m^i)\in \markov$ where $m^x$ and $m^y$ are the $x$- and $y$-coordinates of the hurricane's location, respectively, and $m^i$ represents the hurricane's intensity level. The MC is given by two independent transition probability matrices: one for the hurricane intensity level and one for its location. The intensity level transition matrix is identical to the one used by \cite{pacheco2016forecast}, which considers six levels of intensity $\{0,1,\dots,5\}$ where level 5 corresponds to the maximum hurricane intensity and level 0 corresponds to the case when the hurricane dissipates. The following matrix shows the transition probabilities where $P_{ij}$ represents the probability of transitioning from intensity level $i$ to $j$.
\[
P = \bordermatrix{
~ & 0   & 1   & 2   & 3   & 4   & 5     \cr
0 & 1   & 0   & 0   & 0   & 0   & 0     \cr
1 & 0.11& 0.83& 0.06& 0   & 0.6 & 0     \cr
2 & 0   & 0.15& 0.6 & 0.25& 0   & 0     \cr
3 & 0   & 0   & 0.04& 0.68& 0.28& 0     \cr
4 & 0   & 1   & 0   & 0.18& 0.79& 0.03  \cr
5 & 0   & 0   & 0   & 0   & 0.5 & 0.5   \cr
}
\]

For the hurricane movement, we assume that the hurricane originates at the bottom row of the grid (i.e., $m^y = 0$ for the initial MC state) and that $m^y$ increases by one in each period (i.e., the hurricane advances upwards by one step in each period) and will reach land in exactly $T$ periods. The hurricane can also move on the x-axis by either staying in the same x-coordinate, moving one cell to the left or one cell to the right. To determine the transition probabilities of each movement, we assign random weights to each movement and normalize them. In particular, we assign a uniform weight between 30 and 40 for staying in the same x-coordinate and a weight between 20 and 40 for moving to the right and left cells.  If the current x-coordinate is at one of the grid borders (e.g., $m^x=0$), then the probability of stepping out of the grid is zero. As an example, consider a $4\times 5$ grid, current hurricane location $(m^x,m^y)=(1,1)$, and weights 25, 36, and 22, for moving to locations $(0,2)$, $(1,2)$, $(2,2)$, respectively. Then, the transition probability of moving from location $(1,1)$ to $(0,2)$ is $\frac{25}{83}\approx 0.3$.

The hurricane movement transition probability matrix is generated at random for each instance. The initial $x$-coordinate is chosen uniformly at random considering all possible options. The initial intensity level is also generated with a uniform random distribution considering hurricane intensity levels of two or higher.

\subsection{Demand Generation}

The demand for each shelter depends on the current MC state of the hurricane (i.e., location and intensity). We consider a maximum demand parameter $d^{\max}$ and a maximum distance parameter $\delta^{\max}$ to model this dependency. The maximum demand $d^{\max}$ is chosen at random between values 1000 and 1500 for each instance and  represents the maximum demand for all shelters in a cell. The maximum demand at a cell is randomly distributed across all shelters in that cell. Therefore, each shelter $j\in \shelters$ has a maximum demand $d^{\max}_j$ and the total maximum demand of all shelters in the same cell is equal to $d^{\max}$.

The maximum distance $\delta^{\max}=150$ represents the distance of influence of the hurricane, i.e., all shelters which are farther than  $\delta^{\max}$ from the hurricane have zero demand. If the distance between a shelter $j \in \shelters$ and the hurricane is less than $\delta^{\max}$, then the shelters demand for MC state $m=(m^x, m^y, m^i) \in \markov$ is given by:
\[ d_{j,m} = d^{\max}_j \left(1 - \frac{\delta_{m,j}}{\delta^{\max}}\right)\cdot \left(\frac{m^i}{5}\right)^2,   \]
where $\delta_{m,j}$ is the Euclidean distance between the current location of the hurricane and the shelter. Thus, the demand generation is linear with respect to the hurricane location and quadratic with respect to the hurricane intensity. Note that the demand is zero if the intensity level is 0 and that the maximum demand is achieved when the hurricane is on top of the shelter and has  the highest intensity level (i.e., $\delta_{m,j} = 0$ and $m^i=5$).

\subsection{Initial Capacity, Contingency Modality Options, and Capacity Expansion}

The initial capacity of a DC depends on both $d^{\max}$ and the initial capacity percentage defined in Section \ref{sec:experiments}. The total initial capacity over all DCs in a cell is given by $C^\text{ini}= d^{\max}\cdot p$, where $p \in \{0.2, 0.25, 0.3\}$ is the initial capacity percentage. The total initial capacity of a cell $C^\text{ini}$ is then distributed uniformly at random across all DCs in that cell. Thus, the initial capacity $C_j$ of a particular DC $j \in \centers$ is given by a portion of $C^\text{ini}$.

We generate different contingency modality options depending on the grid size (i.e., number of land cells) and the modality Type (i.e., either Type-1 or Type-2, see Section \ref{sec:experiments}). Our modality set $\modality$ considers options that increase DCs capacity in adjacent cells and options that increase the capacity for all DCs. For example, in a $4\times5$ grid we consider modality options that increase capacity for DCs in cells 0 and 1, 1 and 2, 2 and 3, and options that increase the capacity in all DCs, that is, a total of four different contingency modality options based on location. 

In addition, we consider capacity increments of Type-1 or Type-2, each one with four different options: Type-1 has options $\{10\%, 20\%, 30\%, 40\%\}$ and Type-2 has options $\{15\%, 30\%, 45\%, 60\%\}$. Then each modality option in $\modality$ corresponds to a combination of the affected DCs based on location and the capacity increase. Considering a $4\times5$ grid example with Type-1 modalities, a possible modality option is to increase the capacity of DCs in cells 2 and 3 by 20\% increments. In this example, we have a total of 16 modality options (i.e., four location options and four increments options).

The capacity increments for each modality are based on the initial capacity of the DCs. For example, consider the  $4\times5$ grid with Type-1 modalities and active modality $\ell \in \modality$ affecting  cells 2 and 3 by 20\% increments. Then, a DC $j \in \centers$ in cell 0 has no capacity increments given the current active modality (i.e., $K_{j\ell} =0$), but a DC $j' \in \centers$ in cell 2 with initial capacity $C_{j'}=100$ has capacity increments  $K_{j'\ell}= 20$.

\vspace{1em}
\section{Model Size for HDR Instances}
\label{app:model_size}
Tables \ref{tab:model_size_small} and \ref{tab:model_size_large} show the average number of variables and constraints for both small-size ($4\times5$ grid) and large-size ($5\time6$ grid) instances.

\begin{table}[htbp]
\def\arraystretch{\stretchTableResults}
  \centering
  \caption{Average number of variables and constraints for different aggregations over small-size instances.}
  \scalebox{\scaleTableResults}{
    \begin{tabular}{l|rrrrr}
    \toprule
          & \multicolumn{1}{c}{\hn} & \multicolumn{1}{c}{\ma} & \multicolumn{1}{c}{\pa} & \multicolumn{1}{c}{\mm} & \multicolumn{1}{c}{\fh} \\
    \midrule
    Continuous Variables &        566,370  &        566,370  &        566,370  &        566,370  &        566,370  \\
    Integer Variables &                 80  &            1,120  &            6,720  &          16,528  &          41,760  \\
    Total Variables &        566,450  &        567,490  &        573,090  &        582,898  &        608,130  \\
    Constraints &        120,129  &        124,834  &        129,373  &        137,171  &        164,414  \\
    \bottomrule
    \end{tabular}%
    }
  \label{tab:model_size_small}%
\end{table}%

\begin{table}[htbp]
\def\arraystretch{\stretchTableResults}
  \centering
  \caption{Average number of variables and constraints for different aggregations over large-size instances.}
  \scalebox{\scaleTableResults}{
    \begin{tabular}{l|rrrrr}
    \toprule
          & \multicolumn{1}{c}{\hn} & \multicolumn{1}{c}{\ma} & \multicolumn{1}{c}{\pa} & \multicolumn{1}{c}{\mm} & \multicolumn{1}{c}{\fh} \\
    \midrule
    Continuous Variables &   12,281,770  &   12,281,770  &   12,281,770  &   12,281,770  &   12,281,770  \\
    Integer Variables &               120  &            2,300  &          13,800  &          46,800  &        733,240  \\
    Total Variables &   12,281,890  &   12,284,070  &   12,295,570  &   12,328,570  &   13,015,010  \\
    Constraints &     2,163,164  &     2,175,033  &     2,188,428  &     2,224,132  &     2,932,940  \\
    \bottomrule
    \end{tabular}%
    }
  \label{tab:model_size_large}%
\end{table}%

\vspace{1em}
\section{Computational Effort for Exact Methods}
\label{app:experiments-exact}

We now compare two solution methods for the aggregated models with four types of transformations, that is, solving the extensive form model directly with CPLEX (\extensive) versus the B\&C procedure integrated with the SDDP algorithm proposed in Section~\ref{sec:methodology} (\sddp). Table~\ref{tab:exact_small} shows the average solution time (for instances where optimal solutions are obtained within the time limit) and the average optimality gaps (for unsolved instances) for the small-size instances. We also include the results of solving \fh\ with \extensive\ as a point of reference.  We note that \extensive\ solves to optimality all the instances for all transformations except for 33 when solving \fh, while \sddp\ fails to solve one instance for \pa\ and solves only one instance to optimality for \mm. 

From Table~\ref{tab:exact_small} we see a clear advantage of \extensive\ over \sddp\ in both computational time and number of instances solved within the time limit across all transformations and instance configurations. These results are mostly explained by the fact that approximately 90\% of the computational time is spent inside the SDDP sub-routine. One main factor is the number of SDDP subproblems for each transformation: \hn\ and \ma\ have on average $69$ SDDP subproblems,  \pa\ has $144$, and \mm\  has $294$, which explains the relatively poor performance of \mm.

\begin{table}[tb]
  \def\arraystretch{\stretchTableResults}
  \centering
  \caption{Performance comparison between \sddp\ and \extensive\ for small-size instances ($4\times5$ grid size).}
  \scalebox{\scaleTableResults}{
    \begin{tabular}{cc|rr|rr|rr|rr|r|rr|r}
    \toprule
          &       & \multicolumn{9}{c|}{\textbf{Average Time (sec)}}                      & \multicolumn{3}{c}{\textbf{Gap (\%)}} \\
    \midrule
          &       & \multicolumn{2}{c|}{\hn} & \multicolumn{2}{c|}{\ma} & \multicolumn{2}{c|}{\pa} & \multicolumn{2}{c|}{\mm} & \multicolumn{1}{l|}{\Margarita{\fh}} & \multicolumn{1}{c}{\pa} & \multicolumn{1}{c}{\mm} & \multicolumn{1}{c}{\Margarita{\fh}} \\
    \midrule
    Modality & \multicolumn{1}{l|}{Cap.} & \multicolumn{1}{c}{\extensive} & \multicolumn{1}{c|}{\sddp} & \multicolumn{1}{c}{\extensive} & \multicolumn{1}{c|}{\sddp} & \multicolumn{1}{c}{\extensive} & \multicolumn{1}{c|}{\sddp} & \multicolumn{1}{c}{\extensive} & \multicolumn{1}{c|}{\sddp} & \multicolumn{1}{c|}{\extensive} & \multicolumn{1}{c}{\sddp} & \multicolumn{1}{c}{\sddp} & \multicolumn{1}{c}{\extensive} \\
    \midrule
    \multirow{3}[2]{*}{Type-1} & \multicolumn{1}{l|}{20\%} & \textbf{          260 } &        1,221  & \textbf{          712 } &     5,969  & \textbf{       961 } &     4,553  & \textbf{  1,853 } &  -    &   1,402  & -     & 14.9  & - \\
          & \multicolumn{1}{l|}{25\%} & \textbf{          132 } &           195  & \textbf{          361 } &     1,313  & \textbf{       487 } &     2,760  & \textbf{     626 } &  -    &   3,974  & -     & 3.5   & 0.8 \\
          & \multicolumn{1}{l|}{30\%} & \textbf{            73 } &           102  & \textbf{          109 } &        272  & \textbf{       222 } &     3,380  & \textbf{     367 } &  -    & -     & -     & 2.1   & 0.4 \\
    \midrule
    \multirow{3}[2]{*}{Type-2} & \multicolumn{1}{l|}{20\%} & \textbf{          253 } &        1,138  & \textbf{          919 } &     5,880  & \textbf{    2,047 } &   13,224  & \textbf{  2,951 } &   20,699  &   2,546  & 19.9  & 17.8  & 0.9 \\
          & \multicolumn{1}{l|}{25\%} & \textbf{          167 } &           140  & \textbf{          393 } &     1,324  & \textbf{       651 } &     3,635  & \textbf{     817 } &  -    & -     & -     & 6.2   & 0.8 \\
          & \multicolumn{1}{l|}{30\%} & \textbf{            82 } &           102  & \textbf{          119 } &        205  & \textbf{       271 } &     4,347  & \textbf{     556 } &  -    & -     & -     & 2.3   & 0.6 \\
    \midrule
    \multicolumn{2}{c|}{Average} & \textbf{          161 } &           483  & \textbf{          436 } &     2,494  & \textbf{       773 } &     5,316  & \textbf{  1,195 } &   20,699  &   2,641  &   19.9  &  7.8  &     0.7  \\
    \bottomrule
    \end{tabular}%
    }
  \label{tab:exact_small}%
\end{table}%

The previous results may suggest that \extensive\ requires less computational effort for small-size instances than \sddp. However, large-size instances cannot even be loaded into the solver due to the large number of variables and constraints in the formulation, that is,  approximately 16GB of memory after the presolved phase (see Appendix \ref{app:model_size} for further model size information), so \sddp\ is our only resort for these instances. Table \ref{tab:exact_sddp_large} shows the number of instances solved to optimality, the number of instances where the algorithm found a feasible solution (but could not prove optimality), and the average optimality gap over instances with a feasible solution for \sddp\ over large-size instances. We see that \sddp\ performs quite well for transformations \hn\ and \ma\, solving several instances to optimality; however, it fails to even find any feasible solution for \mm. As before, these results are explained mainly by a large number of SDDP subproblems, and for these large-size instances, \hn\ and \ma\ have on average $114$ subproblems, \pa\ has $249$, and \mm\ has $593$.

\begin{table}[tb]
\def\arraystretch{\stretchTableResults}
  \centering
  \caption{Performance of the SDDP integrated B\&C algorithm in large-size instances ($5\times6$ grid size).}
  \scalebox{\scaleTableResults}{
    \begin{tabular}{cc|rrrr|rrrr|rrrr}
    \toprule
          &       & \multicolumn{4}{c|}{\textbf{\# Optimal}} & \multicolumn{4}{c|}{\textbf{\# Feasible}} & \multicolumn{4}{c}{\textbf{Opt. Gaps (\%)}} \\
    \midrule
    Modality & \multicolumn{1}{l|}{Cap.} & \multicolumn{1}{c}{\hn} & \multicolumn{1}{c}{\ma} & \multicolumn{1}{c}{\pa} & \multicolumn{1}{c|}{\mm} & \multicolumn{1}{c}{\hn} & \multicolumn{1}{c}{\ma} & \multicolumn{1}{c}{\pa} & \multicolumn{1}{c|}{\mm} & \multicolumn{1}{c}{\hn} & \multicolumn{1}{c}{\ma} & \multicolumn{1}{c}{\pa} & \multicolumn{1}{c}{\mm} \\
    \midrule
    \multirow{3}[2]{*}{Type-1} & \multicolumn{1}{l|}{20\%} & 2     & 0     & 0     & 0     & 8     & 10    & 7     & 1     & 53.4  & 51.8  & 56.1  & 82.0 \\
          & \multicolumn{1}{l|}{25\%} & 5     & 2     & 0     & 0     & 5     & 8     & 6     & 0     & 35.7  & 26.2  & 27.9  & - \\
          & \multicolumn{1}{l|}{30\%} & 10    & 8     & 0     & 0     & 0     & 1     & 4     & 0     & -     & 15.4  & 4.3   & - \\
    \midrule
    \multirow{3}[2]{*}{Type-2} & \multicolumn{1}{l|}{20\%} & 3     & 0     & 0     & 0     & 7     & 10    & 7     & 1     & 43.5  & 52.9  & 61.4  & 86.0 \\
          & \multicolumn{1}{l|}{25\%} & 6     & 3     & 0     & 0     & 4     & 7     & 6     & 0     & 37.2  & 24.9  & 32.8  & - \\
          & \multicolumn{1}{l|}{30\%} & 9     & 8     & 0     & 0     & 1     & 1     & 6     & 0     & 12.0  & 13.3  & 6.4   & - \\
    \midrule
    \multicolumn{2}{c|}{Total/Av.} & 35    & 21    & 0     & 0     & 25    & 37    & 36    & 2     & 36.4  & 30.7  & 31.5  & 84.0 \\
    \bottomrule
    \end{tabular}%
    }
  \label{tab:exact_sddp_large}%
\end{table}%

\vspace{1em}
\section{Additional Results for Approximation Methods}
\label{app:approx-additional-resutls}

Tables \ref{tab:2sldr-hn-smallsize} to \ref{tab:2sldr-mm-largesize} show the performance of the approximated methods based on 2SLDR and SDDP for small-size and large-size instances, and transformations \hn, \pa, and \mm.

\begin{table}[htbp]
\def\arraystretch{\stretchTableResults}
  \centering
  \caption{Solution time and quality of 2SLDR and SDDP bounds.  Results for $\hn$ over small-size instances.}
  \scalebox{\scaleTableResults}{
    \begin{tabular}{cc|rrr|r|rrrr|r}
    \toprule
          &       & \multicolumn{4}{c|}{Average Time (sec)} & \multicolumn{5}{c}{Relative Difference (\%)} \\
    \midrule
    Modality & \multicolumn{1}{l|}{Cap.} & \multicolumn{1}{c}{\thldr} & \multicolumn{1}{l}{\tldr} & \multicolumn{1}{l|}{\mldr} & \multicolumn{1}{c|}{\sddplb} & \multicolumn{1}{c}{\thldr} & \multicolumn{1}{l}{\tldr} & \multicolumn{1}{l}{\mldr} & \multicolumn{1}{c|}{\sddpub} & \multicolumn{1}{c}{\sddplb} \\
    \midrule
    \multirow{3}[2]{*}{Type-1} & \multicolumn{1}{l|}{20\%} & 193.7 & 103.6 & \textbf{71.7} & 76.5  & 0.04  & 0.05  & \textbf{0.00} & 0.06  & 0.44 \\
          & \multicolumn{1}{l|}{25\%} & 132.2 & \textbf{31.3} & 33.1  & 15.2  & 0.05  & 0.07  & \textbf{0.00} & \textbf{0.00} & 0.63 \\
          & \multicolumn{1}{l|}{30\%} & 82.5  & 20.0  & \textbf{19.0} & 9.9   & 0.10  & 0.21  & \textbf{0.00} & \textbf{0.00} & 1.36 \\
    \midrule
    \multirow{3}[2]{*}{Type-2} & \multicolumn{1}{l|}{20\%} & 262.8 & \textbf{64.3} & 76.0  & 61.6  & 0.03  & 0.05  & \textbf{0.00} & 0.03  & 0.46 \\
          & \multicolumn{1}{l|}{25\%} & 131.7 & \textbf{29.3} & 35.8  & 12.4  & 0.05  & 0.07  & \textbf{0.00} & 0.01  & 0.62 \\
          & \multicolumn{1}{l|}{30\%} & 80.4  & 19.7  & \textbf{19.0} & 10.0  & 0.10  & 0.21  & \textbf{0.00} & \textbf{0.00}  & 1.36 \\
    \midrule
    \multicolumn{2}{c|}{Average} & 147.2 & 44.7  & \textbf{42.4} & 30.9  & 0.06  & 0.11  & \textbf{0.00} & 0.02  & 0.81 \\
    \bottomrule
    \end{tabular}%
    }
  \label{tab:2sldr-hn-smallsize}%
\end{table}%

\begin{table}[htbp]
\def\arraystretch{\stretchTableResults}
  \centering
  \caption{Solution time and quality of 2SLDR and SDDP bounds. Results for $\hn$ over large-size instances.}
  \scalebox{\scaleTableResults}{
    \begin{tabular}{cc|rrrl|rrr}
    \toprule
          &       & \multicolumn{4}{c|}{Average Time (sec)} & \multicolumn{3}{c}{Opt. Gap (\%)} \\
    \midrule
    Modality & Cap.  & \multicolumn{1}{c}{\tldr} & \multicolumn{1}{c}{\mldr} & \multicolumn{1}{c}{\sddplb} & (opt) & \multicolumn{1}{c}{\tldr} & \multicolumn{1}{c}{\mldr} & \multicolumn{1}{c}{\sddpub} \\
    \midrule
    \multirow{3}[2]{*}{Type-1} & \multicolumn{1}{l|}{20\%} &     2,100.9  &     1,961.8  & \textbf{       449.3} & (10)  & 0.51  & 0.33  & \textbf{0.13} \\
          & \multicolumn{1}{l|}{25\%} &     1,679.4  &     1,628.6  & \textbf{       288.0} & (10)  & 0.78  & 0.39  & \textbf{0.15} \\
          & \multicolumn{1}{l|}{30\%} &        965.1  &        864.2  & \textbf{         56.6} & (10)  & 0.61  & 0.48  & \textbf{0.18} \\
    \midrule
    \multirow{3}[2]{*}{Type-2} & \multicolumn{1}{l|}{20\%} &     2,501.3  &     1,794.3  & \textbf{       730.9} & (10)  & 0.56  & 0.31  & \textbf{0.12} \\
          & \multicolumn{1}{l|}{25\%} &     1,660.5  &     1,764.6  & \textbf{       203.7} & (10)  & 0.78  & 0.40  & \textbf{0.13} \\
          & \multicolumn{1}{l|}{30\%} &        950.1  &     1,035.3  & \textbf{         57.9} & (10)  & 0.62  & 0.49  & \textbf{0.18} \\
    \midrule
    \multicolumn{2}{c|}{Av. (Total)} &     1,642.9  &     1,508.1  & \textbf{       297.7} & (60)  & 0.64  & 0.40  & \textbf{0.15} \\
    \bottomrule
    \end{tabular}%
    }
  \label{tab:2sldr-hn-largesize}%
\end{table}%

\begin{table}[htbp]
\def\arraystretch{\stretchTableResults}
  \centering
  \caption{Solution time and  quality of 2SLDR and SDDP bounds. Results for $\ma$ over small-size instances.}
  \scalebox{\scaleTableResults}{
    \begin{tabular}{cc|rrr|r|rrrr|r}
    \toprule
          &       & \multicolumn{4}{c|}{Average Time (sec)} & \multicolumn{5}{c}{Relative Difference (\%)} \\
    \midrule
    Modality & \multicolumn{1}{l|}{Cap.} & \multicolumn{1}{c}{\thldr} & \multicolumn{1}{l}{\tldr} & \multicolumn{1}{l|}{\mldr} & \multicolumn{1}{c|}{\sddplb} & \multicolumn{1}{c}{\thldr} & \multicolumn{1}{l}{\tldr} & \multicolumn{1}{l}{\mldr} & \multicolumn{1}{c|}{\sddpub} & \multicolumn{1}{c}{\sddplb} \\
    \midrule
    \multirow{3}[2]{*}{Type-1} & \multicolumn{1}{l|}{20\%} & 431.1 & 104.2 & \textbf{95.2} & 367.3 & 0.26  & 0.29  & 0.23 & \textbf{0.08}  & 0.19 \\
          & \multicolumn{1}{l|}{25\%} & 368.3 & 58.4 & \textbf{49.2}  & 167.4 & 0.14  & 0.18  & \textbf{0.11} & 0.22 & 1.11 \\
          & \multicolumn{1}{l|}{30\%} & 119.3 & 25.7  & \textbf{25.4} & 17.1  & 0.11  & 0.22  & 0.02 & \textbf{0.00} & 1.56 \\
    \midrule
    \multirow{3}[2]{*}{Type-2} & \multicolumn{1}{l|}{20\%} & 630.6 & 119.4 & \textbf{117.1} & 439.4 & 0.12  & 0.14  & 0.10 & \textbf{0.09}  & 0.26 \\
          & \multicolumn{1}{l|}{25\%} & 464.3 & 49.3 & \textbf{46.9}  & 111.2 & 0.15  & 0.18  & \textbf{0.11} & 0.18  & 1.15 \\
          & \multicolumn{1}{l|}{30\%} & 113.8 & \textbf{22.7}  & 23.1 & 16.4  & 0.10  & 0.21  & 0.01 & \textbf{0.00}  & 1.58 \\
    \midrule
    \multicolumn{2}{c|}{Average} & 354.5 & 63.3  & \textbf{59.5} & 186.5 & 0.15  & 0.21  & \textbf{0.09} & 0.10  & 0.97 \\
    \bottomrule
    \end{tabular}%
    }
  \label{tab:2sldr-ma-smallsize}%
\end{table}%

\begin{table}[htbp]
\def\arraystretch{\stretchTableResults}
  \centering
  \caption{Solution time and quality of 2SLDR and SDDP bounds. Results for $\ma$ over large-size instances.}
  \scalebox{\scaleTableResults}{
    \begin{tabular}{cc|rrrl|rrr}
    \toprule
          &       & \multicolumn{4}{c|}{Average Time (sec)} & \multicolumn{3}{c}{Opt. Gap (\%)} \\
    \midrule
    Modality & Cap.  & \multicolumn{1}{c}{\tldr} & \multicolumn{1}{c}{\mldr} & \multicolumn{1}{c}{\sddplb} & (opt) & \multicolumn{1}{c}{\tldr} & \multicolumn{1}{c}{\mldr} & \multicolumn{1}{c}{\sddpub} \\
    \midrule
    \multirow{3}[2]{*}{Type-1} & \multicolumn{1}{l|}{20\%} &     2,266.5  &     2,855.8  & \textbf{    1,913.0} & (10)  & 0.50  & 0.32  & \textbf{0.13} \\
          & \multicolumn{1}{l|}{25\%} &     3,067.1  &     2,458.6  & \textbf{    1,295.8} & (10)  & 0.84  & 0.46  & \textbf{0.19} \\
          & \multicolumn{1}{l|}{30\%} &     1,051.1  &     1,028.4  & \textbf{       133.1} & (10)  & 0.62  & 0.49  & \textbf{0.19} \\
    \midrule
    \multirow{3}[2]{*}{Type-2} & \multicolumn{1}{l|}{20\%} &     4,065.5  & \textbf{    3,488.9} &     3,500.0  & (10)  & 0.59  & 0.34  & \textbf{0.14} \\
          & \multicolumn{1}{l|}{25\%} &     4,295.9  &     4,107.4  & \textbf{    1,960.6} & (10)  & 0.89  & 0.40  & \textbf{0.18} \\
          & \multicolumn{1}{l|}{30\%} &     1,272.7  &     1,252.9  & \textbf{       150.3} & (10)  & 0.65  & 0.52  & \textbf{0.20} \\
    \midrule
    \multicolumn{2}{c|}{Av. (Total)} &     2,669.8  &     2,532.0  & \textbf{    1,492.1} & (60)  & 0.68  & 0.42  & \textbf{0.17} \\
    \bottomrule
    \end{tabular}%
    }
  \label{tab:2sldr-ma-largesize}%
\end{table}%

\begin{table}[htbp]
\def\arraystretch{\stretchTableResults}
  \centering
  \caption{Solution time and quality of 2SLDR and SDDP bounds. Results for $\mm$ over small-size instances.}
  \scalebox{\scaleTableResults}{
    \begin{tabular}{cc|rrr|r|rrrr|r}
    \toprule
          &       & \multicolumn{4}{c|}{Average Time (sec)} & \multicolumn{5}{c}{Relative Difference (\%)} \\
    \midrule
    Modality & \multicolumn{1}{l|}{Cap.} & \multicolumn{1}{c}{\thldr} & \multicolumn{1}{l}{\tldr} & \multicolumn{1}{l|}{\mldr} & \multicolumn{1}{c|}{\sddplb} & \multicolumn{1}{c}{\thldr} & \multicolumn{1}{l}{\tldr} & \multicolumn{1}{l}{\mldr} & \multicolumn{1}{c|}{\sddpub} & \multicolumn{1}{c}{\sddplb} \\
    \midrule
    \multirow{3}[2]{*}{Type-1} & \multicolumn{1}{l|}{20\%} & 2766.2 & \textbf{262.2} & 277.0 & 3521.5 & 0.12  & 0.26  & 0.11  & \textbf{0.00} & 0.39 \\
          & \multicolumn{1}{l|}{25\%} & 1272.7 & 176.6 & \textbf{154.9} & 12447.8 & 0.09  & 0.13  & \textbf{0.04} & 0.25  & 0.98 \\
          & \multicolumn{1}{l|}{30\%} & 4049.5 & \textbf{235.3} & 310.3 & 14746.8 & 0.13  & 0.25  & \textbf{0.00} & 1.37  & 1.77 \\
    \midrule
    \multirow{3}[2]{*}{Type-2} & \multicolumn{1}{l|}{20\%} & 3916.3 & 504.6 & \textbf{396.0} & 11703.3 & 0.13  & 0.26  & 0.13  & \textbf{0.02} & 0.22 \\
          & \multicolumn{1}{l|}{25\%} & 1486.7 & 226.3 & \textbf{164.9} & 21226.9 & 0.09  & 0.12  & \textbf{0.01} & 0.38  & 2.30 \\
          & \multicolumn{1}{l|}{30\%} & 8042.7 & \textbf{431.0} & 720.4 & 15323.7 & 0.13  & 0.24  & \textbf{0.01} & 0.95  & 1.80 \\
    \midrule
    \multicolumn{2}{c|}{Average} & 3589.0 & \textbf{306.0} & 337.2 & 13161.7 & 0.11  & 0.21  & \textbf{0.05} & 0.50  & 1.24 \\
    \bottomrule
    \end{tabular}%
    }
  \label{tab:2sldr-mm-smallsize}%
\end{table}%

\begin{table}[htbp]
\def\arraystretch{\stretchTableResults}
  \centering
  \caption{Solution time and quality of 2SLDR and SDDP bounds. Results for $\mm$ over large-size instances.}
  \scalebox{\scaleTableResults}{
    \begin{tabular}{cc|rrrl|rrr}
    \toprule
          &       & \multicolumn{4}{c|}{Average Time (sec)} & \multicolumn{3}{c}{Opt. Gap (\%)} \\
    \midrule
    Modality & Cap.  & \multicolumn{1}{c}{\tldr} & \multicolumn{1}{c}{\mldr} & \multicolumn{1}{c}{\sddplb} & (opt) & \multicolumn{1}{c}{\tldr} & \multicolumn{1}{c}{\mldr} & \multicolumn{1}{c}{\sddpub} \\
    \midrule
    \multirow{3}[2]{*}{Type-1} & \multicolumn{1}{l|}{20\%} &     8,288.1  & \textbf{    4,776.9} &   20,125.0  & (1)   & 39.71 & \textbf{33.00} & 39.39 \\
          & \multicolumn{1}{l|}{25\%} &   14,398.1  & \textbf{  12,203.6} &   21,601.1  & (0)   & 15.73 & \textbf{15.08} & 17.53 \\
          & \multicolumn{1}{l|}{30\%} & \textbf{    1,911.4} &     2,076.6  &   21,601.6  & (0)   & 2.41  & \textbf{0.67} & 2.09 \\
    \midrule
    \multirow{3}[2]{*}{Type-2} & \multicolumn{1}{l|}{20\%} &   14,668.7  & \textbf{  13,135.3} &   21,600.7  & (0)   & 39.55 & \textbf{37.42} & 43.68 \\
          & \multicolumn{1}{l|}{25\%} &   16,280.7  & \textbf{  16,051.0} &   21,601.6  & (0)   & 23.98 & \textbf{23.22} & 26.70 \\
          & \multicolumn{1}{l|}{30\%} &     2,391.2  & \textbf{    1,925.4} &   21,601.4  & (0)   & 3.63  & \textbf{1.51} & 3.50 \\
    \midrule
    \multicolumn{2}{c|}{Av. (Total)} &     9,656.4  & \textbf{    8,361.5} &   21,355.2  & (1)   & 20.84 & \textbf{18.48} & 22.15 \\
    \bottomrule
    \end{tabular}%
    }
   \label{tab:2sldr-mm-largesize}%
\end{table}%

\vspace{1em}
\section{Additional Results for Policy Structures and Managerial Insights}
\label{app:policy-managerial}

Table \ref{tab:solution-structure-type2} presents the same metrics shown in Table \ref{tab:solution-structure-type1} but for Type-2 small-size instances.

\begin{table}[htbp]
  \centering
  \caption{Solution structures for different policies and initial capacities for Type-2 instances.}
    \scalebox{\scaleTableResults}{
    \begin{tabular}{l|rrrrr|rrrrr|rrrrr|rrrrr}
    \toprule
          & \multicolumn{5}{c|}{Nodes (\%)			}    & \multicolumn{5}{c|}{\# of Contingencies				} & \multicolumn{5}{c|}{Aggressiveness (\%)		} & \multicolumn{5}{c}{Intensity} \\
    \midrule
    Cap.  & \multicolumn{1}{c}{\hn} & \multicolumn{1}{c}{\ma} & \multicolumn{1}{c}{\pa} & \multicolumn{1}{c}{\mm} & \multicolumn{1}{c|}{\fh} & \multicolumn{1}{c}{\hn} & \multicolumn{1}{c}{\ma} & \multicolumn{1}{c}{\pa} & \multicolumn{1}{c}{\mm} & \multicolumn{1}{c|}{\fh} & \multicolumn{1}{c}{\hn} & \multicolumn{1}{c}{\ma} & \multicolumn{1}{c}{\pa} & \multicolumn{1}{c}{\mm} & \multicolumn{1}{c|}{\fh} & \multicolumn{1}{c}{\hn} & \multicolumn{1}{c}{\ma} & \multicolumn{1}{c}{\pa} & \multicolumn{1}{c}{\mm} & \multicolumn{1}{c}{\fh} \\
    \midrule
    20\%  & 100   & 92    & 72    & 72    & 53    & 1.0   & 1.0   & 1.0   & 1.0   & 4.1   & 14    & 39    & 52    & 52    & 46    & 2.3   & 2.5   & 2.9   & 2.9   & 2.9 \\
    25\%  & 10    & 38    & 22    & 16    & 31    & 1.0   & 1.0   & 1.9   & 3.8   & 5.3   & 15    & 20    & 51    & 55    & 33    & 2.3   & 2.9   & 4.0   & 4.1   & 3.5 \\
    30\%  & 0     & 7     & 15    & 14    & 17    & -     & 1.0   & 1.9   & 4.8   & 5.7   & -     & 14    & 29    & 27    & 23    & -     & 3.0   & 4.2   & 4.1   & 4.0 \\
    \bottomrule
    \end{tabular}%
    }
  \label{tab:solution-structure-type2}%
\end{table}%

\end{APPENDICES}

\end{document}